\documentclass[3p]{elsarticle}

\usepackage{bm}
\usepackage[parfill]{parskip}
\usepackage[utf8]{inputenc}
\usepackage{url}
\usepackage{amsmath,amssymb,amsfonts,amsthm}
\usepackage{color}
\usepackage{bm}
\usepackage{subcaption}
\usepackage[linesnumbered,ruled]{algorithm2e}

\usepackage[section]{placeins}
\makeatother

\makeatletter
\renewcommand{\fnum@figure}{Fig. \thefigure}
\makeatother

\makeatletter
\renewcommand{\fnum@table}{Tab. \thetable}
\makeatother

\journal{arXiv.org}

\newcommand{\TheTitle}{Ray Effect Mitigation for the Discrete Ordinates Method 
through Quadrature Rotation} 

\newcommand{\SN}{S$_N$\,}
\newcommand{\PN}{P$_N$\,}
\newcommand{\rSN}{rS$_N$\,}

\newcommand{\bOmega}{{\bm{\Omega}}}
\newcommand{\bx}{{\bm{x}}}
\newcommand{\bn}{{\bm{n}}}
\newcommand{\bP}{{\bm{P}}}
\newcommand{\bxi}{{\bm{\xi}}}

\begin{document}

\begin{frontmatter}

\title{\TheTitle}

\author[adressThomas]{Thomas Camminady}
\author[adressMartin]{Martin Frank}
\author[adressKerstin]{Kerstin K\"upper}
\author[adressJonas]{Jonas Kusch}

\address[adressThomas]{Karlsruhe Institute of Technology, Karlsruhe,
    thomas.camminady@kit.edu}
    \address[adressMartin]{Karlsruhe Institute of Technology, Karlsruhe,
    martin.frank@kit.edu}
\address[adressKerstin]{RWTH Aachen University, Aachen,  
kuepper@mathcces.rwth-aachen.de}
\address[adressJonas]{Karlsruhe Institute of Technology, Karlsruhe,
    jonas.kusch@kit.edu}

\begin{abstract}
Solving the radiation transport equation is a challenging task, due to the high 
dimensionality of the solution's phase space. The commonly used discrete 
ordinates (S$_N$) method 
suffers from ray effects which result from a break in rotational symmetry from 
the finite set of directions chosen by S$_N$.
The spherical harmonics (P$_N$) equations, on the other hand, preserve 
rotational symmetry, but can produce negative particle densities. 
The discrete ordinates (S$_N$) method, in turn, by construction ensures 
non-negative particle densities.

In this paper we present a modified version of the \SN method, the rotated \SN 
(rS$_N$) method. Compared to \SN, we add a rotation and interpolation step for 
the angular quadrature points and the respective function values after every time step. Thereby, the number of directions 
on which the solution evolves is effectively increased and ray effects are 
mitigated. Solution values on rotated ordinates are computed by an 
interpolation step.
%, which is facilitated by the use of mesh-based quadrature rules.
Implementation details are provided and in our experiments the rotation and 
interpolation step 
only adds $5\%$ to $10\%$ to the runtime of the \SN method. 
We apply the \rSN method to the line-source and a lattice test case, both being 
prone to ray effects. Ray effects are reduced significantly, even for small 
numbers of quadrature 
points. The \rSN method yields qualitatively similar solutions to the \SN method with less than a third of the number of 
quadrature points, both for the line-source and the lattice problem. 
The code used to produce our results is freely available and can be downloaded \cite{JuliaWN}.

%\mf{Publish the code???}
\end{abstract}


\end{frontmatter}

\section{Introduction}

Many physics applications such as high-energy astrophysics, supernovae \cite{fryer2006snsph,swesty2009numerical} and fusion \cite{matzen2005pulsed,marinak2001three} require accurate solutions of the radiation transport equation. 
Numerically solving this equation is challenging, since the phase space on which the solution (the angular flux) is defined is at least six dimensional, consisting of three spatial dimensions, two directional (angular) parameters and time. In many applications, there is additional frequency dependence.

Various angular discretization strategies exist and they all come with certain advantages and shortcomings (cf.\ \cite{brunner2002forms} for a comparison): 
The spherical harmonics (P$_N$) method 
\cite{case1967linear,pomraning1973equations,lewis1984computational}  expands the solution in terms of 
angular variables with finitely many spherical harmonics basis functions. A 
hyperbolic system of equations for the expansion coefficients is then obtained 
by testing with the spherical harmonics basis. The convergence of P$_N$ is 
spectral and the solution preserves the property of rotational invariance. 
However, the P$_N$ method suffers from oscillations in non-smooth 
regimes, leading to non-physical, negative values of the angular flux and the 
density.

Stochastic Methods such as Monte Carlo (MC), e.g.~\cite{fleck1971implicit}, 
preserve positivity and are considered to yield accurate solutions without ray 
effects. However, MC suffers from noise due to the finite number of stochastic 
samples. 

A third method is the discrete ordinates (S$_{N}$) method  
\cite{lewis1984computational}, which preserves positivity of the angular 
flux. In order to eliminate the dependency on 
angular parameters, the key idea of \SN is to solve the radiation transfer 
equation on a fixed angular grid, which results in a set of equations that 
couple through a collision term. Due to the fact that a finite number of 
possible directions is imposed on the solution, the resulting angular flux 
suffers from so called ray effects 
\cite{lathrop1968ray,morel2003analysis,mathews1999propagation}. The solution exhibits rays 
that correspond to the discrete set of directions chosen for the \SN grid. 
Consequently, one obtains a solution with poor accuracy, violating the 
rotational invariance property. Sufficiently increasing the number of ordinates 
solves this problem, but at the cost of a heavily increased run-time. 

More sophisticated strategies to mitigate ray effects make use of biased quadrature sets, which reflect the importance 
of certain ordinates \cite{abu2001angular}. In \cite{tencer2016ray}, the 
angular flux is computed for several differently oriented quadrature sets and 
ray effects are then mitigated by averaging over all solutions. Moreover, a 
method combining \SN and P$_N$ to reduce rays has been introduced in 
\cite{lathrop1971remedies} with further refinements given in 
\cite{jung1972discrete,reed1972spherical,miller1977ray}. The idea is to use a 
mixture of collocation points as well as basis functions to represent the 
solution's angular dependency. Consequently, a system for the angular expansion 
coefficients with an increased coupling of the individual equations needs to be 
solved. The accuracy of these methods has been studied in the review paper 
\cite{morel2003analysis} and it turns out that all methods still suffer from 
ray effects for a line-source in a void.

Comparisons of S$_N$, \PN and MC methods have for example been studied in 
\cite{mcclarren2010robust} for the line-source, lattice and hohlraum problem. 
These test cases are designed to show the respective disadvantages of
the \SN and \PN method. That is, solutions tend to show ray effects (for S$_{N}$) 
and become negative (for P$_N$), respectively.

In this paper, we propose a new strategy to mitigate the formation of rays when using 
the \SN method. The key idea is to allow for an effectively larger number of directions along 
which particles can travel. This is done by rotating the set of ordinates after 
each time step around a random axis, meaning that the solution is evolved on an 
enlarged set of ordinates. The solution at the rotated ordinates is then 
obtained via interpolation. To guarantee an efficient interpolation procedure, 
we make use of a quadrature set similar to \cite{thurgood1995tn}, which is based on a triangulation of the unit sphere. The resulting connectivity is used to 
efficiently find relevant 
interpolation points for the new ordinates. Since the interpolation step preserves positivity of the solution values, the resulting method inherits positivity from \SN. We demonstrate the effectiveness of 
our method, which we call \rSN method, by studying the line-source and lattice problem, 
where we observe that the method reduces ray effects while leading to positive solution values 
at affordable computational overhead.

This paper is structured as follows: In Section \ref{sec:background}, the 
radiation transport equation as well as its \SN discetization is presented. A 
simple rotated \SN version is then introduced in Section 
\ref{sec:RotationSimple}, for which we quantify the smoothing effect 
of the rotations and show numerical results. We extend the rotational 
axis of the quadrature set to arbitrary directions in Section 
\ref{sec:RotationArbitrary}. The straight-forward extension of the \SN method 
to the \rSN method is discussed in Section 
\ref{sec:Implementation}. Numerical examples are shown in Section 
\ref{sec:Results}.

\section{Background}
\label{sec:background}
In this section, we present the governing equations as well as their numerical discretization using \SN.
\subsection{The radiation transport equation}
The radiation transport equation is a linear integro-differential equation and describes the evolution of particles traveling through a background medium. It is given by
\begin{align}\label{eq:kineticEquation}
\partial_t \psi(t,\bx,\bOmega) + \bOmega \cdot \nabla_\bx \psi(t,\bx,\bOmega) + 
\sigma_t(\bx) 
\psi(t,\bx,\bOmega)= \sigma_s(\bx) \int_{\mathbb{S}^2} \psi(t,\bx,\bOmega') \, 
d\bOmega'.
\end{align}

%\mf{Die Strahlungstransportgleichung war im gesamten Paper falsch. Es fehlte 
%der Verlustterm durch Streuung. Bitte checken ob er jetzt ueberall steht. 
%Bitte 
%Rechenergebnisse checken.}

In this equation, $\psi$ is the angular flux and depends on time 
$t\in\mathbb{R}^+$, spatial position $\bm{x}=(x,y,z)^T \in\mathbb{R}^3$ and 
direction of travel $\bOmega\in\mathbb{S}^2$. The units are chosen so that particles travel with unit speed. The first two terms of 
\eqref{eq:kineticEquation} describe streaming, i.e.\ particles move in the 
direction $\bOmega$ without any interaction with the background material. The 
function $\sigma_t(\bx)= \sigma_s(\bx)+\sigma_a(\bx)$ is the total cross section which
describes the loss due to absorption and scattering by the material. The right 
hand side 
describes the gain of particles with direction $\bOmega$ due to incoming scattering. For simplicity, we consider isotropic scattering, 
but all methods here can easily be applied for anisotropic scattering.

\subsection{Discrete ordinates method}
A first step when deriving numerical schemes to calculate the solution $\psi$ 
is to discretize the direction of travel $\bOmega$. The idea of the discrete 
ordinates method is to describe the direction of travel by a fixed set of 
finitely many directions (or ordinates), meaning that the solution is computed on the set $\{ 
\bOmega_1,\cdots,\bOmega_{N_q} \}\subset \mathbb{S}^2$. The solution is now 
described by
\begin{align*}
\psi_q(t,\bx):=\psi(t,\bx,\bOmega_q) \enskip \text{ with } q = 1,\cdots,N_q.
\end{align*}
Defining a quadrature rule
\begin{align*}
\int_{\mathcal{S}^2} \psi(t,\bx,\bOmega) \,d\bOmega \approx \sum_{q = 1}^{N_q} 
w_q 
\psi_q(t,\bx),
\end{align*}
allows a discretization of the radiation transport equation \eqref{eq:kineticEquation} by
\begin{align*}
\partial_t \psi_q(t,\bx) + \bOmega_q \cdot \nabla_\bx \psi_q(t,\bx)+ 
\sigma_t(\bx)  
\psi_q(t,\bx) = \sigma_s(\bx)  \sum_{p = 1}^{N_q} w_p \psi_p(t,\bx) \enskip 
\text{ 
with }\, 
q = 
1,\cdots,N_q.
\end{align*}
The resulting system of $N_q$ partial differential equations only depends on 
$t$ and $\bm{x}$, meaning that it can be solved by a standard finite volume 
scheme. 

\subsection{Finite volume discretization}
\label{sec:FVM}
For the purpose of presenting the algorithm, we describe a first-order finite volume discretization. For the numerical experiments, 
we use a second-oder method in both space and time (a re-implementation of \cite{garrett2013comparison}).

We start by dividing the spatial domain into cells 
\begin{align*}
V_{ijl}: = [x_i,x_{i+1}]\times [y_j,y_{j+1}]\times [z_l,z_{l+1}],
\end{align*}
with volume $\vert V_{ijl} \vert$. Furthermore, the time domain 
is decomposed into equidistant time steps $t_0,\cdots,t_{N_t}$, where $\Delta t 
:= t_{n+1}-t_n$. In every spatial cell at every time step, the averaged 
solution is
\begin{align*}
\psi_{q,ijl}^n \simeq \frac{1}{\vert V_{ijl} \vert}\int_{V_{ijl}} 
\psi_q(t_n,\bx) \,d\bm{x}.
\end{align*}
The finite volume method is then given by
\begin{align}\label{eq:FVDiscretization}
\psi_{q,ijl}^{n+1} =& \psi_{q,ijl}^n - \frac{\Delta t}{\vert V_{ijl} \vert} \left( g_{i+1/2,j,l} - g_{i-1/2,j,l} + g_{i,j+1/2,l}- g_{i,j-1/2,l}+ g_{i,j,l+1/2}- g_{i,j,l-1/2} \right)\nonumber \\
&+\Delta t\, \sigma_{s,ijl}\sum_{p=1}^{N_q} w_p \psi_{p,ijl}^n- 
\sigma_{t,ijl}\, \psi_{q,ijl}^n\;, 
\end{align}
where the numerical flux at the interface between cells $V_{i,j,l}$ and $V_{i,j+1,l}$ with unit normal $\bm{n}$ is given by
\begin{align*}
g_{i,j+1/2,l} = 
\begin{cases}
\bm{n}^T\bm{\bOmega}_q \psi_{q,ijl}^n & \text{ if }\bm{n}^T\bm{\bOmega}_q>0\\
\bm{n}^T\bm{\bOmega}_q \psi_{q,i,j+1,l}^n & \text{ else}
\end{cases}\;.
\end{align*}
The remaining numerical fluxes are chosen analogously. A typical S$_N$ solution 
is depicted in Fig.~\ref{fig:SNsolutionLS}. In the following, we propose a 
method to mitigate the non-physical ray-effects.

\section{Two-dimensional case}
\label{sec:RotationSimple}
\subsection{Rotation and interpolation}
In this section, we introduce our idea in a two-dimensional setting. In addition to evolving the angular flux in time by repeatedly calling the finite volume 
update \eqref{eq:FVDiscretization}, we rotate the set of ordinates after each 
timestep around the z-axis. This simplified setting is considered, since rotation around the z-axis as well as the corresponding interpolation step are straight forward when using a standard tensorized quadrature. The case when using an arbitrary rotation will be considered in Sec.~\ref{sec:RotationArbitrary}.

First, we present the commonly used tensorized quadrature set on the sphere: To discretize the direction $\bOmega\in\mathbb{S}^2$, we define the azimuthal 
and the polar angles $\theta$ and $\phi$, so that
\begin{equation}\label{eq:spherical-coordinates}
\bOmega = (\cos\phi \sin \theta, \sin\phi \sin \theta, \cos \theta)^T\;.
\end{equation}
We use a product quadrature on the sphere with some arbitrary quadrature for 
$\mu = \cos(\theta) \in[-1,1]$ (e.g.\ Gauss quadrature for $\mu$) and equally 
weighted, equally spaced points for $\phi$. Thus, let 
\begin{equation}
\phi_i = i \Delta\phi \quad \text{for} \quad i=1,\ldots,N_q \quad \text{and} 
\quad \Delta\phi = \frac{2\pi}{N_q}\;,
\end{equation}
where $N_q$ is the number of quadrature points. Due to the alignment of $\mu = 
\cos(\theta)$ with the z-axis, a rotation around the z-axis only affects the 
azimuthal angle $\theta$. If we rotate the quadrature set by an angle 
$\delta\in(0,\Delta\phi)$, we can approximate the solution on the rotated 
points $\phi_i^\delta = \phi_i+\delta$, $i=1,\ldots,N_q$ using linear 
interpolation
\begin{equation}\label{eq:interpolation}
\psi(t_n,\bx,\theta,\phi_i^\delta) = (1-a) \psi(t_n,\bx,\theta,\phi_i) + 
a\psi(t_n,\bx,\theta,\phi_{i+1})\;,
\end{equation}
where $a = \frac{\delta}{\Delta\phi}\in(0,1)$ and $\phi_{N_q+1} = \phi_1$. 
We call this a \emph{rotation} and \emph{interpolation} step. After having 
interpolated the solution at the new ordinates $\bOmega^{\delta} = 
(\cos\phi_i^\delta \sin \theta, \sin\phi_i^\delta \sin \theta, \cos \theta)^T$, 
a finite volume step is performed on the new ordinates, i.e.
\begin{align}\label{eq:FVDiscretization2}
\psi_{q,ijl}^{\delta,n+1} =& \psi_{q,ijl}^{\delta,n} - \frac{\Delta t}{\vert V_{ijl} \vert} \left( g_{i+1/2,j,l}^{\delta} - g_{i-1/2,j,l}^{\delta} + g_{i,j+1/2,l}^{\delta}- g_{i,j-1/2,l}^{\delta}+ g_{i,j,l+1/2}^{\delta}- g_{i,j,l-1/2}^{\delta} \right)\nonumber \\
&+\Delta t  \sigma_{s,ijl}\sum_{p=1}^{N_q} w_q \psi_{p,ijl}^{\delta,n}-  \sigma_{t,ijl} \psi_{q,ijl}^{\delta,n}\;, 
\end{align}
with
\begin{align*}
g_{i+1/2,j,l}^{\delta} = 
\begin{cases}
\bm{n}^T\bOmega_q^{\delta} \psi_{q,ijl}^{\delta,n} & \text{ if 
}\bm{n}^T\bOmega_q^{\delta}>0\\
\bm{n}^T\bOmega_q^{\delta} \psi_{q,i,j+1,l}^{\delta,n} & \text{ else}
\end{cases}\;.
\end{align*}
This process is repeated until a specified final time is reached. Conservation 
of the zeroth order moment is guaranteed due to linear interpolation.

\subsection{Modified equation analysis}
In this section we analyze the effect of the interpolation. This analysis is based on modified equations, which are a common technique to determine the dispersion or diffusion of the numerical discretization of a hyperbolic balance law \cite[Chapter~11.1]{levequenumerical}. In the following, we assume that $\delta$ is fixed and we rotate back and forth between the original and rotated quadrature set. For this setting, we analyze the effects resulting from the combination of rotation/interpolation and update steps. 

For simplicity, we only consider the advection operator (which is responsible 
for the ray effects) 
\begin{equation}\label{eq:advection}
\partial_t \psi + \bOmega\cdot \nabla_\bx \psi = 0\;,
\end{equation}
i.e., collisions and absorption are omitted. The update in time is performed by the 
explicit Euler method with time step $\Delta t$, that is
\begin{equation}\label{eq:update}
\psi(t_{n+1},\bx,\bOmega) = \psi(t_n,\bx,\bOmega)-\Delta t \bOmega \cdot 
\nabla_\bx 
\psi(t_n,\bx,\bOmega) \quad \text{with} \quad t_n = n\Delta t \quad \text{and} 
\quad n = 0,1,2,\ldots\;.
\end{equation}
In our scheme, the update step~\eqref{eq:update} and the rotation 
step~\eqref{eq:interpolation} are alternating. In the following, we want to 
analyze the concatenation of a rotation around the z-axis by an angle $\delta$, 
an update from $t_n$ to $t_{n+1}$ for some $n$, and another rotation around the 
z-axis by an angle $-\delta$, so that we return to the original set of 
quadrature points. Note that since we are only interested in 
investigating how an interpolation and rotation step affects the standard \SN time 
update, we are not considering a full cycle of our scheme 
as this would include another time update. Furthermore, the spatial variable 
$\bx$ and the polar angle $\theta$ are continuous variables for now.

First, we apply the interpolation
\begin{equation}
\psi(t_{n+1},\bx,\theta,\phi_i) = (1-a) 
\psi(t_{n+1},\bx,\theta,\phi_i^\delta) +a \psi(t_{n+1}, \bx, \theta, 
\phi_{i-1}^\delta)\;.
\end{equation}
Second, we perform the update in time
\begin{equation}
\begin{aligned}
\psi(t_{n+1},\bx,\theta,\phi_i) =
(1-a)&\left( \psi(t_n,\bx,\theta, \phi_i^\delta)-\Delta t \bOmega_i^\delta 
\cdot \nabla_\bx  \psi(t_n,\bx,\theta, \phi_i^\delta) \right) \\
+ a &\left(\psi(t_n,\bx,\theta, \phi_{i-1}^\delta)-\Delta t 
\bOmega_{i-1}^\delta \cdot \nabla_\bx  \psi(t_n,\bx,\theta, 
\phi_{i-1}^\delta)\right)\;,
\end{aligned}
\end{equation}
where $\bOmega_i^\delta$ is defined according~\eqref{eq:spherical-coordinates} 
using $\theta_i^\delta$. Finally, we apply the interpolation again
\begin{equation}
\begin{aligned}
\psi(t_{n+1},\phi_i) =
(1-a) \biggl(&(1-a)\psi(t_n,\phi_i) + a \psi(t_n,\phi_{i+1})\\
& -\Delta t \bOmega_i^\delta \cdot \nabla_\bx \left((1-a) 
\psi(t_n,\phi_i)+a\psi(t_n,\phi_{i+1})\right)\biggr) \\
+ a \biggl(&(1-a) 
\psi(t_n,\bx,\phi_{i-1})+a\psi(t_n,\bx,\phi_{i})\\
&-\Delta t \bOmega_{i-1}^\delta \cdot \nabla_\bx  \left((1-a)\psi(t_n, 
\phi_{i-1})+a\psi(t_n, \phi_{i})\right)\biggr)\;.
\end{aligned}
\end{equation}
Here, we omitted the arguments $\bx$ and $\theta$ of the solution $\psi$ to 
shorten the notation. The above equation can be rewritten as
\begin{equation}\label{eq:discrete-concatenation}
\begin{aligned}
\frac{\psi(t_{n+1},\phi_i)-\psi(t_n,\phi_i)}{\Delta t}+ \bOmega_i\cdot 
\nabla_\bx 
\psi(t_n,\phi_i) 
= \frac{a(1-a)}{\Delta t} \biggl(\psi(t_n,\phi_{i+1}) -2\psi(t_n,\phi_i)+\psi(t_n,\phi_{i-1})\biggr)\\
-a(1-a)\biggl( \bOmega_i^\delta \cdot \nabla_\bx \psi(t_n,\phi_{i+1}) 
+\bOmega_{i-1}^\delta \cdot \nabla_\bx \psi(t_n,\phi_{i-1} ) \biggr)\\
-\biggl((1-a)^2 \bOmega_i^\delta+a^2 \bOmega_{i-1}^\delta 
-\bOmega_i\biggr) \cdot \nabla_\bx \psi(t_n,\phi_i)\;,
\end{aligned}
\end{equation}
so that the left-hand side is a discretization of the advection 
equation~\eqref{eq:advection} and the right-hand side is the result of the 
rotation and interpolation. 

Now we require that the scheme has a non-trivial limit for $\Delta t\to 0$. Because of the term $\frac{a(1-a)}{\Delta t}$, we have to choose
 $a = c\Delta t$ for some constant $c$. When $\Delta t \to 0$, we then get the following limiting equation
\begin{equation}
\label{eq2dmod}
\partial_t \psi(t,\phi_i) + \bOmega_i\cdot \nabla_\bx \psi(t,\phi_i) = c\Delta 
\phi^2 \frac{\psi(t,\phi_{i+1}) -2\psi(t,\phi_i)+\psi(t,\phi_{i-1})}{\Delta 
\phi^2}\;.
\end{equation} 
This is a semi-discretized advection equation with a discrete second-order 
derivative in the azimuthal angle on the right-hand side, i.e. 
$\frac{\partial^2\psi}{\partial\phi^2}$. However, the right-hand side scales 
with $\Delta \phi^2$, so that the diffusive effect of the second-order 
derivative vanishes with increasing angular resolution~$\Delta \phi\to 0$. On 
the other hand, for fixed $\Delta t$ and $\Delta\phi\to0$, the above 
Eq.~\eqref{eq:discrete-concatenation} becomes
\begin{equation}
\frac{\psi(t_{n+1},\phi)-\psi(t_n,\phi)}{\Delta t}+ \bOmega\cdot \nabla_\bx 
\psi(t_n,\phi) = 0\;.
\label{eq:modequation}
\end{equation}
This means that the effect of the rotation vanishes when the angular discretization is refined.

The important point of this analysis is that the rotation introduces diffusion 
(in the angular variable) into the system and we have to choose $a$ 
proportional to $\Delta t$, i.e.\ $a=c\Delta t\in[0,1]$ for some 
constant~$c$. In particular, the angle of the rotation $\delta = a \Delta 
\phi = c\Delta t \Delta \phi$ is then proportional to the timestep $\Delta t$ and the angular discretization $\Delta 
\phi$. 

%For an arbitrary rotation in three dimensions, we expect to obtain a discrete 
%Laplace-Beltrami operator, but this requires a different 
%quadrature/interpolation rule.

\subsection{Numerical results for S$_N$ with rotation}
We briefly discuss numerical results for the S$_N$ solution with and without rotation 
around the z-axis. The solution of S$_N$ with $N=8$ is computed, i.e. we make 
use of $2\cdot 
N$ equidistant discretization points for the angle $\Phi$ and $N/2$ Gauss 
quadrature points for $\mu$, i.e.\ the total number of quadrature points is $N_q=N^2=64$. For the S$_N$ method with 
rotation, we rotate the quadrature set back and forth by an angle of 
$\delta=10\, 
\Delta t \Delta \phi$. 

%\mf{Darueber steht dass alpha in [0,1] sein muss???}

\begin{figure}[h!]
\begin{subfigure}{0.32\linewidth}
	\centering
	\includegraphics[scale=0.24]{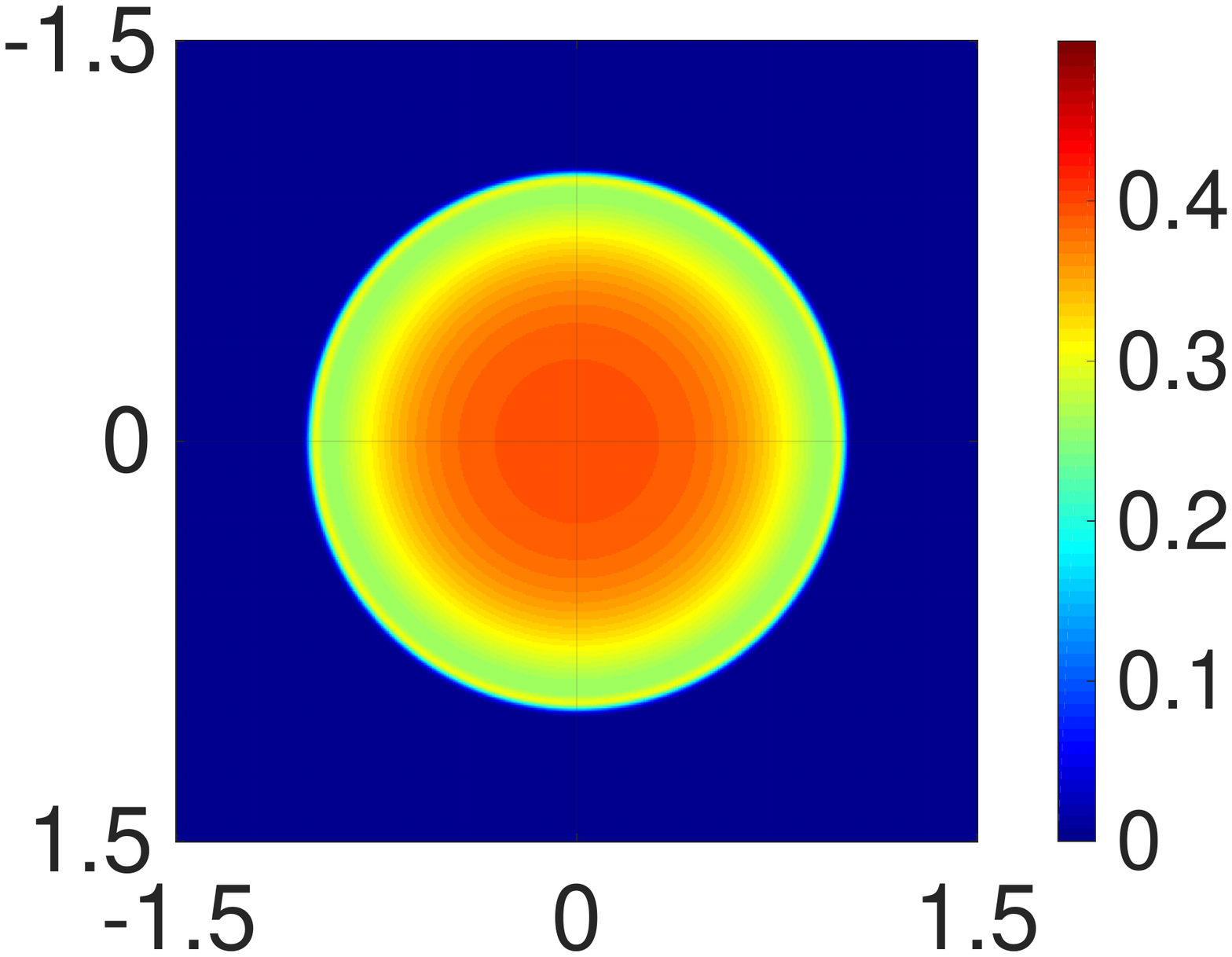}	
	\label{fig:referencesolutionLS}
	\vspace*{-20mm} 
	\caption{Analytical solution.}
\end{subfigure}% gure}% 
\begin{subfigure}{0.32\linewidth}
	\centering
	\includegraphics[scale=0.24]{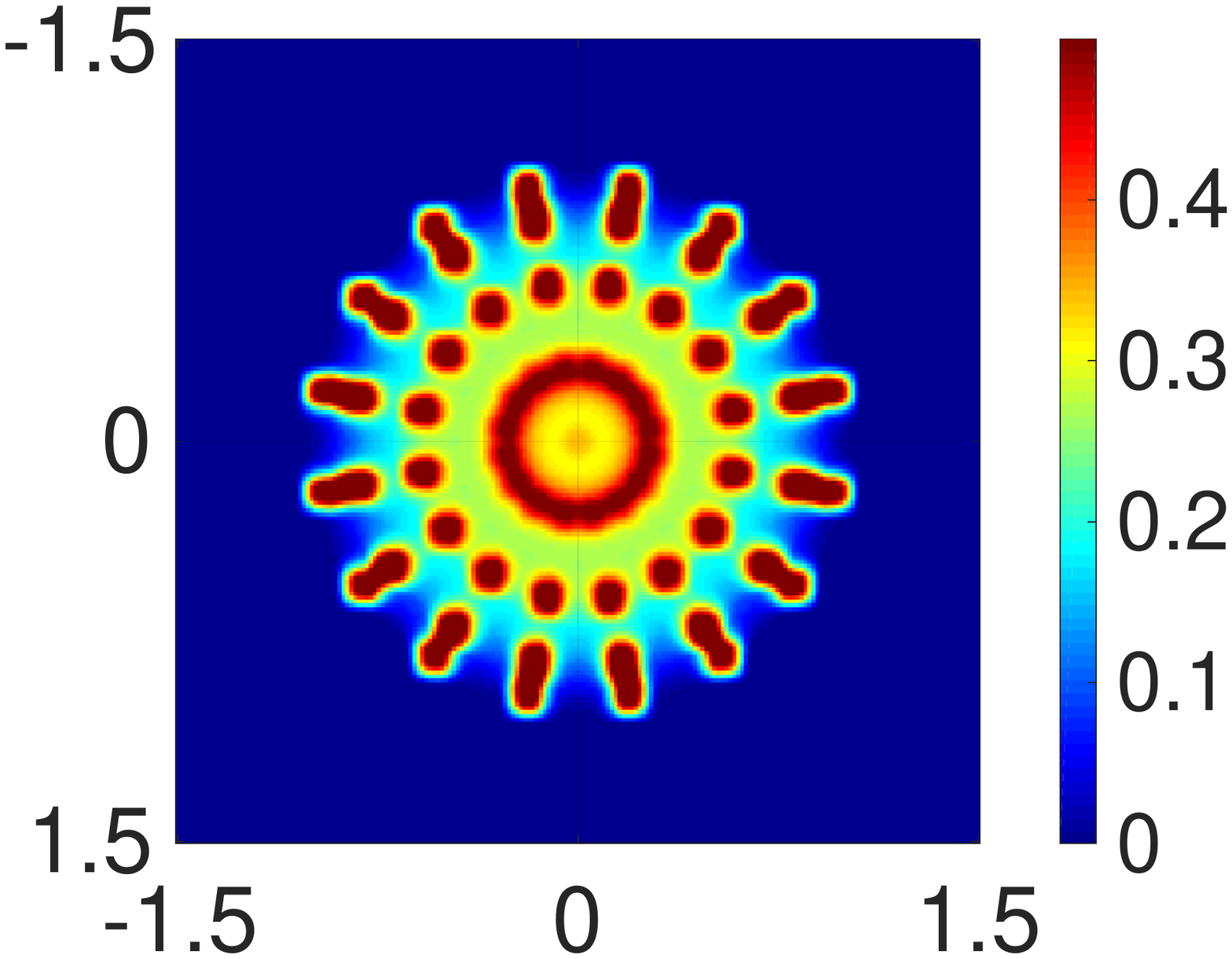}	
		\vspace*{-20mm} 
	\caption{S$_8$ solution with ray effects.}
		\label{fig:SNsolutionLS}
\end{subfigure}% 
\begin{subfigure}{0.32\linewidth}
	\centering
	\includegraphics[scale=0.24]{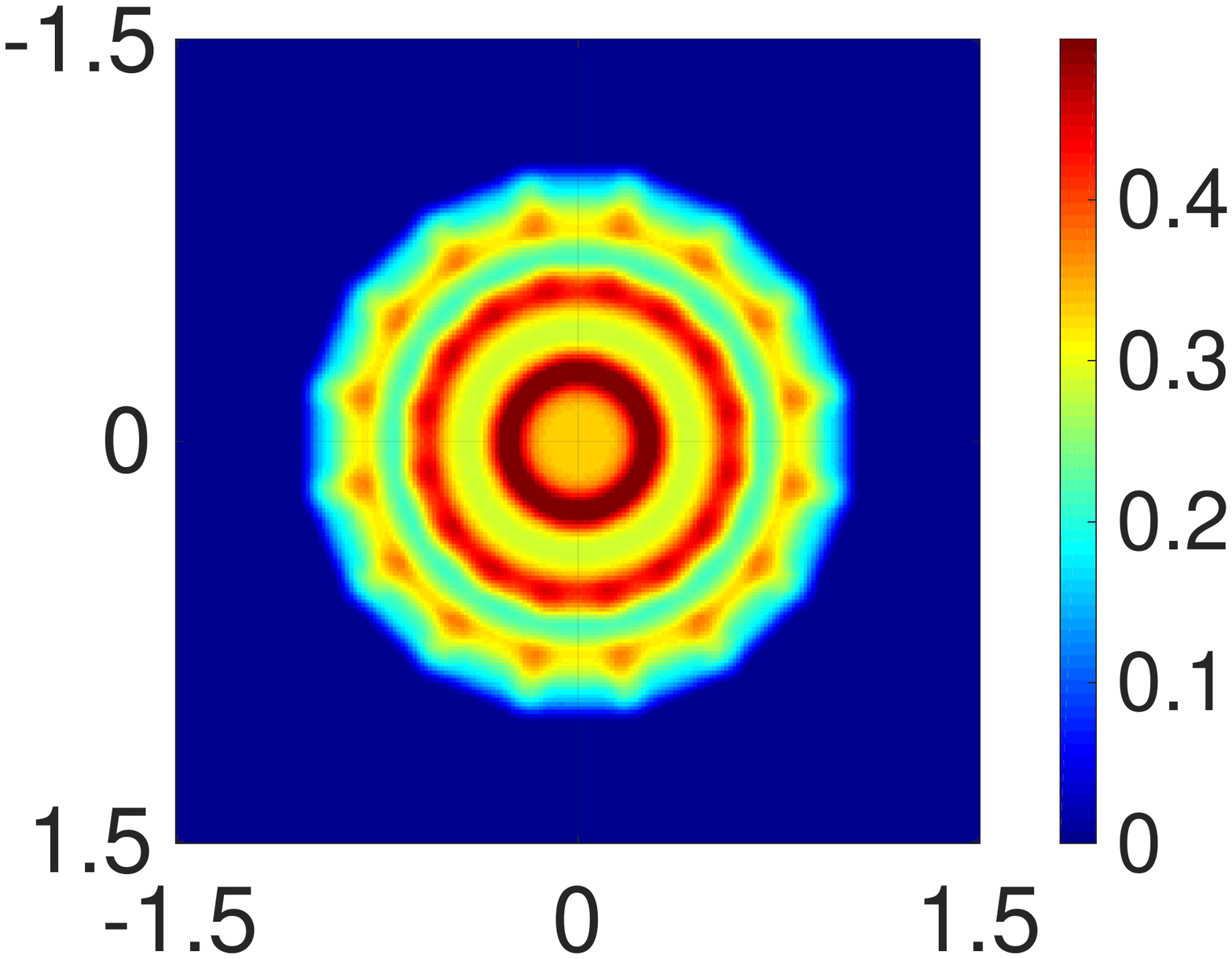}	
		\vspace*{-20mm} 
	\caption{S$_8$ solution with rotation.}
\end{subfigure}%
\caption{Line-source test case.}
\label{fig:LineSourceSetting}
\end{figure}

Further parameters of the 
computation as well as 
details on the line-source test case 
can be found in Section \ref{sec:LineSourceResults}. Results of the density
\begin{align}\label{eq:density}
\rho(t_{\text{end}},\bx) = \int_{\mathbb{S}^2} 
\psi(t_{\text{end}},\bx,\bOmega)\,d\bOmega
\end{align}
plotted on the physical domain $\bm{x}\in\mathbb{R}^2$ are shown in 
Fig.~\ref{fig:LineSourceSetting}.
The exact solution has the property of being rotationally symmetric, which is 
especially violated by the \SN method, since the solution suffers from ray 
effects. By rotating the ordinates of the \SN method around the z-axis, one 
mitigates ray effects. However, oscillations are still present in the radial 
dimension, because we only rotate around the $z$-axis.

%https://www.researchgate.net/profile/A_Pollard/publication/245361815_The_TN_Quadrature_Set_for_the_Discrete_Ordinates_Method/links/563bb55308aec6f17dd4e968.pdf

\section{Three-dimensional case}
\label{sec:RotationArbitrary}
To generalize the procedure explained in Section~\ref{sec:RotationSimple} to rotations around an 
arbitrary axis, we need a quadrature set that allows for easy interpolation of 
the rotated points in every spatial cell. For this purpose, a quadrature set 
that is the result of an underlying triangulation of the unit sphere is chosen. 
Given this triangulation, function values at rotated points can be interpolated 
via barycentric interpolation on the sphere. The quadrature points and weights, as well as the triangulation, result from 
projecting a triangulation of planar triangles onto the sphere.
% We want to 
%obtain a quadrature with a small variance in the distribution of quadrature 
%weights. Therefore we present different possible ways to construct this 
%quadrature set based on different \textit{base geometries} and different 
%triangulations of these geometries.

\subsection{Quadrature points}
\label{sec:quadpoints}

\begin{figure}[h!]
\begin{center}
		\includegraphics[scale = 
		0.6]{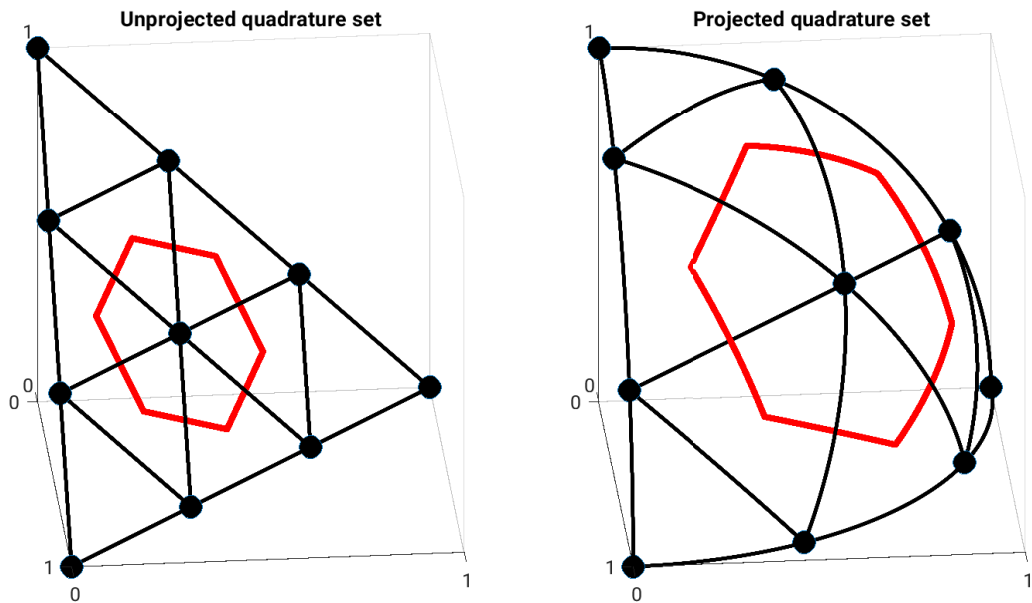}
		\caption{Triangulation for $N=4$ with the corresponding quadrature 
		points in the plane and on the surface of the unit sphere, together 
		with the connectivity and the quadrature weight for the center point in 
		red. Quadrature points and weights for the seven other octants result 
		from symmetry.}
		\label{fig:construct_quad_points1}
	\end{center}
\end{figure}

For the quadrature points, consider one face of the standard octahedron, i.e.\ the triangle 
with nodes $(1,0,0), (0,1,0)$ and $(0,0,1)$. This triangle is now 
triangulated in an equidistant 
way as shown in Fig. \ref{fig:construct_quad_points1} on the left. All vertices 
of the triangulation are projected onto the unit sphere $\mathbb{S}^2$ in a next step, 
presented Fig.\ \ref{fig:construct_quad_points1} on the right.

%\tc{this is different from the TN quadrature, see: The TN	Quadrature  Set  
%for  the  	Discrete  Ordinates  Method}

With $N$ being the number of points on the line segment between the points $A$ 
and $B$, and therefore between $B$ and $C$, as well as between $C$ and $A$, the 
total number of quadrature points on the whole sphere, denoted by $N_q$, is 
given 
by the relation $N_q=4N^2-8N+6$.
Each vertex belongs to six distinct triangles, except for the six vertices at 
the poles which only belong to four distinct triangles.

%To obtain a smaller variance in the distribution of the quadrature weight, we 
%can consider a face of the regular icosahedron as a base geometry from which 
%we 
%then project the quadrature points. In this case, the number of quadrature 
%points is $N_q = 10N^2-20N+12$.

%Another possible modification of the described procedure can be used to reduce 
%the variance of the weight distribution even further. Instead of equidistantly 
%placing the points in the planar geometry followed by a projection onto the 
%sphere, the points an be placed equidistantly with respect to their distance 
%on 
%the sphere. This is commonly referred to as spherical linear interpolation, or 
%\textit{slerp}. 
Our construction of the quadrature set is very similar to the idea of the T$_N$ 
quadrature \cite{thurgood1995tn}. However, instead of taking the midpoints of 
the resulting triangles in Fig. \ref{fig:construct_quad_points1}, we take the 
surrounding vertices. This changes the number of quadrature points 
that are being generated and the corresponding weights. More importantly 
however, it directly yields a connectivity between quadrature points which we 
use for the interpolation later on. For the T$_N$ quadrature, it is not clear 
how to connect vertices from different octants. Our proposed quadrature does 
not suffer from this ambiguity since vertices fall on the connecting edges 
between the octants, thus allowing us to keep a triangulation for the whole 
unit sphere.

\subsection{Quadrature weights}
\label{sec:quadweights}
For each quadrature point, the quadrature weight corresponds to the area 
associated with that given point.
%This area is defined by connecting the midpoints of all surrounding
%triangles. 
%For the poles, the resulting shape is a quadrilateral and a hexagon for all 
%other points. The associated quadrature weight is then given by the associated 
%area when projected onto the surface of the unit ball. 
This area is defined by first connecting the midpoints of all surrounding
triangles on the unprojected grid and then projecting this shape onto the surface of the unit ball. 
For the poles, the resulting shape is a quadrilateral and a hexagon for all 
other points.

%When choosing the \textit{slerp} method to generate the triangulation, the 
%associated area is given by connecting the midpoints of the surrounding 
%triangles on the unit sphere.

\begin{figure}[h!]
	\centering
	\includegraphics[width=1.\linewidth]{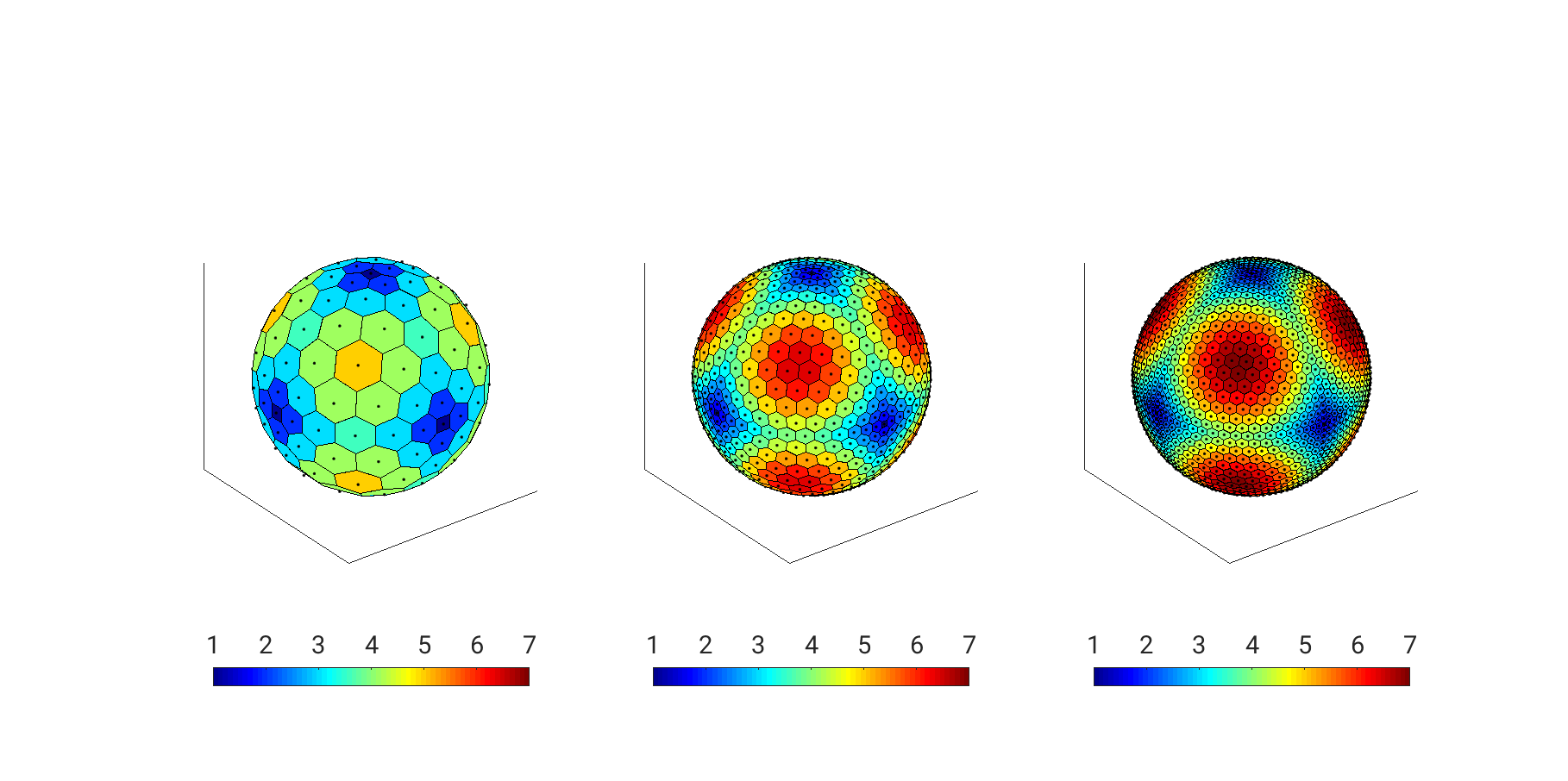}
	\caption{Weight ratio distribution for the quadrature set. Color coded is 
	the ratio to the minimal quadrature weight found. The maximal ratio 
	converges 
	towards $9\sqrt{3}/2\approx 7.7942$. The 
		number of quadrature points is $N_q=146$, $N_q=678$ and $N_q=1602$, 
		respectively.}
	\label{fig:weightscolorcoded}
\end{figure}

The complete set of quadrature weights is shown in 
Fig.~\ref{fig:weightscolorcoded}, together with the corresponding hexagons or 
quadrilaterals. 
Due to symmetry, it is sufficient to compute the quadrature weights on a single 
octant and then copy them to all seven other octancts.
% For the icosahedron, we 
%proceed analogously.
%A boxplot for the distribution of the quadrature weights is shown in 
%Fig. \ref{fig:weightdistribution}. We observe the smallest variance for a 
%quadrature set that is generated with the \textit{slerp} method with base 
%geometry being the icosahedron. The outliers that can be observed in the 
We observe the smallest quadrature weights for the poles. 
This is due to two effects. Firstly, the increase of the corresponding area when projecting 
it onto the unit sphere is smaller than for other points. Secondly, only four 
neighboring triangles contribute to the quadrature weights at the poles, as 
opposed to six for all other points.
%/home/thomas/PhD/Projects/QuadratureOnSphere/src/matlab/mainplots.m
%146         678        1602
%\begin{figure}[h!]
%	\centering
%	\includegraphics[scale = 0.3]{figs/weightdistribution.png}
%	\caption{Boxplot for different types of quadratures with the 
%	number of quadrature points being approximately $N_q\approx100.000$. The 
%	outliers correspond to the quadrature weights for the poles. }
%	\label{fig:weightdistribution}
%\end{figure}
%\newpage

\subsection{Accuracy}
We now integrate certain functions on the sphere with the described quadrature 
to check the implementation and investigate its accuracy. A more detailed analysis of the quadrature rule, and possible improvements, will 
be the topic of a folllow-up paper.
%As an 
%example, the following results are based on the octahedron version, using 
%linearly spaced 
%points.
We consider functions, mapping from $\mathbb{S}^2$ to $\mathbb{R}$ and compute 
the error for different numbers of quadrature points against the analytic 
solution. The two functions used to test the accuracy of the chosen quadrature 
are
\begin{align*}
%&f(x,y,z) = xyz \text{ and 
%}I_{f} = \int_{\mathbb{S}^2} f\, d\bOmega=0, \\
&g(x,y,z) = x^4y^2 \text{ and 
}I_{g} = \int_{\mathbb{S}^2} g\, d{\bOmega}=\frac{4}{35} \pi,\\
&h(x,y,z) = \cos(x)+\sin(y)+z^6 \text{ and 
}I_{h}= \int_{\mathbb{S}^2} h\, d\bOmega=\frac{4}{7} 
\pi(1+7\sin(1)).
\end{align*}

In Table \ref{tab:tab1} we show the errors resulting from the numerical 
integration with the quadrature set of the specified order $N$. Integrated were 
two different functions $g$ and $h$. The error ratio is the absolute 
value of the ratio of two consecutive errors.
The results show, that the order of convergence is two with respect 
to $N$, which implies a first order convergence with respect to the number of 
quadrature points $N_q$. First order convergence was expected from the construction.
We did not aim for a high order quadrature rule, as we are primarily using the described quadrature due to its low 
variance in quadrature weights and the naturally arising interpolation 
capabilities. 

	%The table lists the results for applying the quadrature with different 
	%number of points.
	\begin{center}
		\begin{tabular}{ll|ll|ll|ll}
			$N$ & $N_q$ &  Error for 
			$g$ & Error ratio & Error for $h$ & Error ratio\\
			\hline
			2 & 6 &  -0.359039  & - & 	2.46015  & - \\
			4 & 38 &  -0.012968 & 27.6865 & 	0.073617 & 33.4183  \\
			8 & 198 &  -0.00234195 &5.53729 & 	0.0265397 & 
			2.77384\\
			16 & 902 &-0.000530132 &4.41767 & 	0.00607712 
			& 4.36715\\
			32 & 3846 & -0.000125148 & 4.23606& 0.00143802 
			&4.22602 \\
			64 & 15878 & -3.03595e-05 & 4.12219 & 	0.000349035 
			& 4.11999
		\end{tabular}
	\label{tab:tab1}
	\captionof{table}{Error for integration of $g$ and $h$. The results 
	indicate a first order convergence with respect to the number of quadrature 
	points used. }
	\end{center}

\subsection{Rotation and interpolation}
\label{sec:quadinterp}
Rotation of the quadrature set is straight forward. A rotation is defined by an 
axis $\bn=(n_x,n_y,n_z)^T\in\mathbb{R}^3$ with 
$||\bn||_2=1$ and a rotation magnitude $\delta$. For a quadrature point 
$q\in\mathbb{S}^2$, $R_\delta^{\bn} q$ rotates $q$ around $\bn$ by 
an amount $\delta$, with the rotation matrix
\begin{align}
R_\delta^\bn = \begin{pmatrix} 
n_x^2 \left(1-\cos(\delta)\right) +\cos(\delta) & u_x u_y \left(1-(\cos 
\delta)\right) - n_z \sin(\delta) & n_x n_z \left(1-\cos(\delta)\right) + n_y 
\sin(\delta) \\ 
n_y n_x \left(1-\cos(\delta)\right) + n_z \sin(\delta) &   
n_y^2\left(1-\cos(\delta)\right)+\cos(\delta) & n_y n_z \left(1-\cos (
\delta)\right) 
- n_x 
\sin(\delta) \\ 
n_z n_x \left(1-\cos(\delta)\right) - n_y \sin(\delta) & n_z n_y \left(1-\cos( 
\delta)\right) + n_x \sin(\delta) &   n_z^2\left(1-\cos( 
\delta)\right)+\cos(\delta)
\end{pmatrix}.
\label{eq:rotmatrix}
\end{align}
When performing the rotation of the quadrature points, the associated 
quadrature weights are kept. The same holds true for the connectivity 
between the vertices that define the triangulation.

\begin{figure}[h!]
	\begin{center}
		\includegraphics[scale = 0.60]{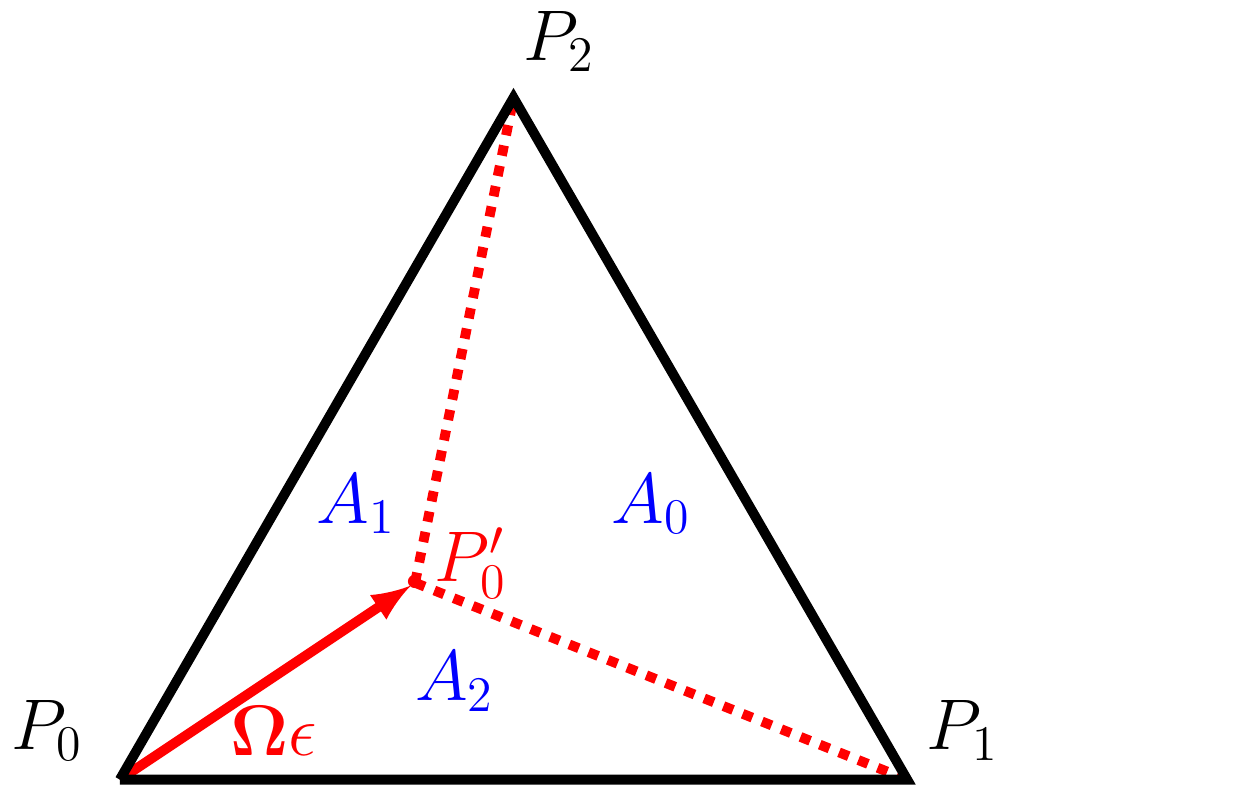}
		\caption{Interpolation weights are given proportional to the covered 
		area 
		inside the triangle.}
		\label{fig:weights_in_triangle}
	\end{center}
\end{figure}

In a next step, we interpolate function values on the 
rotated quadrature point set, given function values on the original quadrature 
point set.
To do so, we use the triangulation that was set up to create the quadrature 
points originally.
Each rotated point falls into one triangle of the original triangulation.
We can then interpolate a function value for any rotated point by the three 
function values at the vertices belonging to the triangle that the rotated 
point lies in. Interpolation is then performed via barycentric interpolation. 
The barycentric interpolation is visualized in Fig. 
\ref{fig:weights_in_triangle} for the planar case. When interpolating a new 
function value at $\bP_0'$, we sum the function values at $\bP_i$ with weights 
$w_i 
= A_i/A$ for $i=0,1,2$ and $A$ being the area of the triangle. For the 
spherical case, all areas are computed as areas on the unit sphere.

Before computing the relevant interpolation weights, we need to find the 
corresponding triangle that a quadrature point has been rotated into. As we 
will restrict ourselves to small rotations later on, these triangles are the 
neighboring triangles for any quadrature point. However, the interpolation 
procedure works for any rotation magnitude. 
Since each spatial cell has the same set of quadrature points, the 
interpolation weights have to be computed only once. Afterwards, the 
interpolation can be performed in each spatial cell separately.

%\subsubsection{Different ways to interpolate like shepard? }

%\tc{Diese Subsection kann weg?! }

\subsection{Modified equation for the planar case}
\label{sec:modeq}
We will now consider a simplified setting to explore the effects of the 
rotation and interpolation step analytically. 
The fact that the quadrature is still anisotropic makes the analysis 
on the sphere difficult. We postpone a more detailed discussion to the end of this subsection. 

\begin{figure}[h!]
	\begin{center}
		\includegraphics[scale = 0.4]{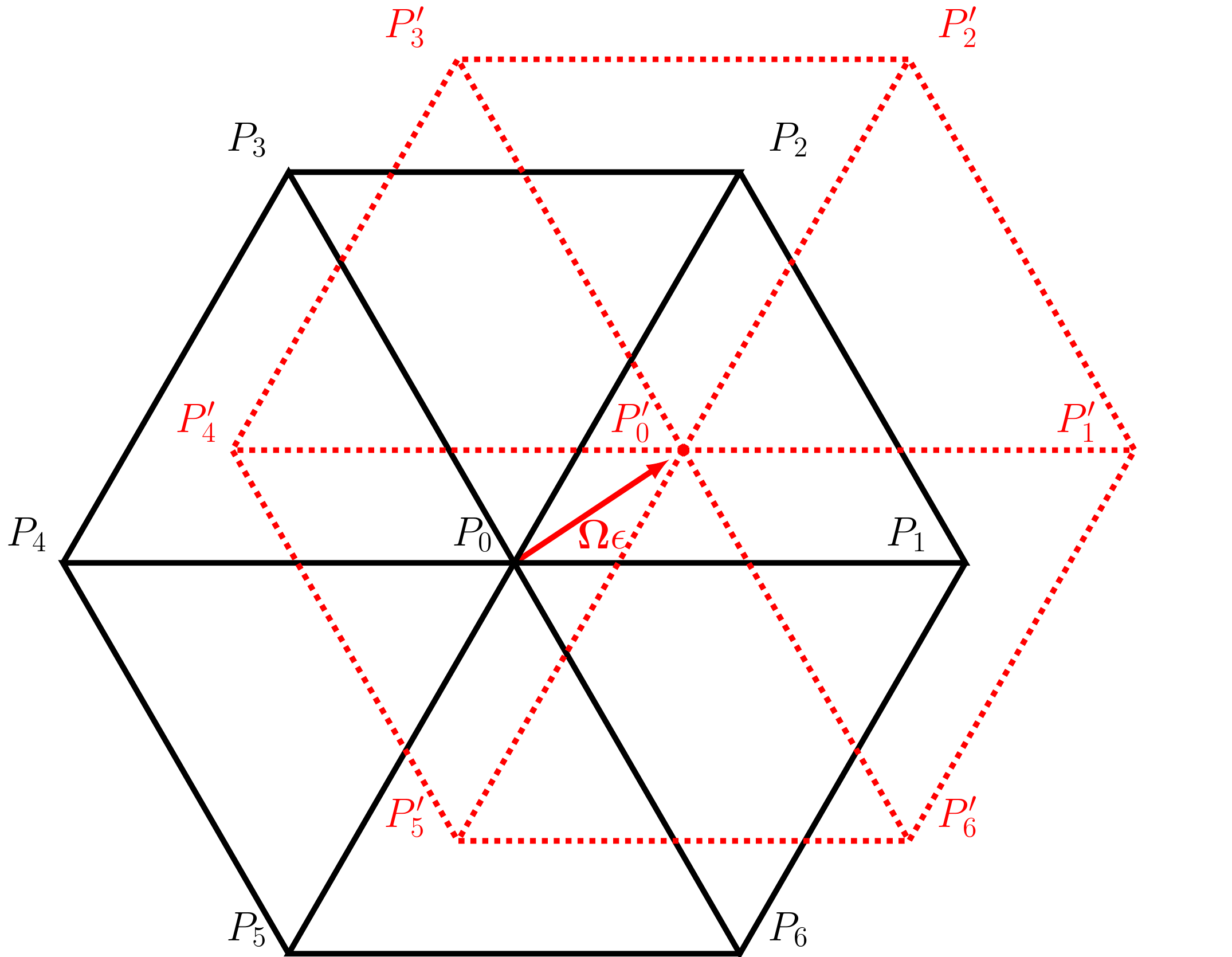}
		\caption{Original set of quadrature points (black) and the translated points (red), together with the respective triangulation.}
	\label{fig:hexagonshift2}	
\end{center}
\end{figure}

Assume therefore a triangulation in planar geometry with equilateral triangles 
that is 
being translated by a vector $\bOmega \, \epsilon$, where $ 
{\bOmega}=(\cos(\alpha),\sin(\alpha))^T$. An excerpt of the original points, 
together 
with the surrounding triangles is shown in black and the shifted points with 
the 
corresponding triangles in dashed red in Fig. \ref{fig:hexagonshift2}. Each 
point $\bP_i$ in the original set of points is being shifted to a new point 
$\bP_i'=\bP_i+\bOmega\, \epsilon$.

The interpolation weights are again the barycentric weights, shown in Fig. 
\ref{fig:weights_in_triangle}, that is $w_1 = A_1/A, w_2 = A_2/A$ and 
$w_0=1-w_1-w_2$, with $A=A_0+A_1+A_2$ being the area of the equilateral 
triangle.
The interpolation weights can be computed analytically and expressed in terms 
of $\epsilon$ and $\alpha$ as
\begin{align*}
	w_1(\alpha) 
	&=\frac{2}{\sqrt{3}}\sin\left(\pi/3-\alpha\right)\epsilon=c_1(\alpha) 
	\epsilon,\\
	w_2(\alpha) &= \frac{2}{\sqrt{3}}\sin(\alpha)\epsilon=c_2(\alpha) 
	\epsilon,\\
	w_0(\alpha) &= 1-(c_1(\alpha)+c_2(\alpha))\epsilon.
\end{align*}

\begin{figure}[h!]
\centering
\includegraphics[width=0.5\linewidth]{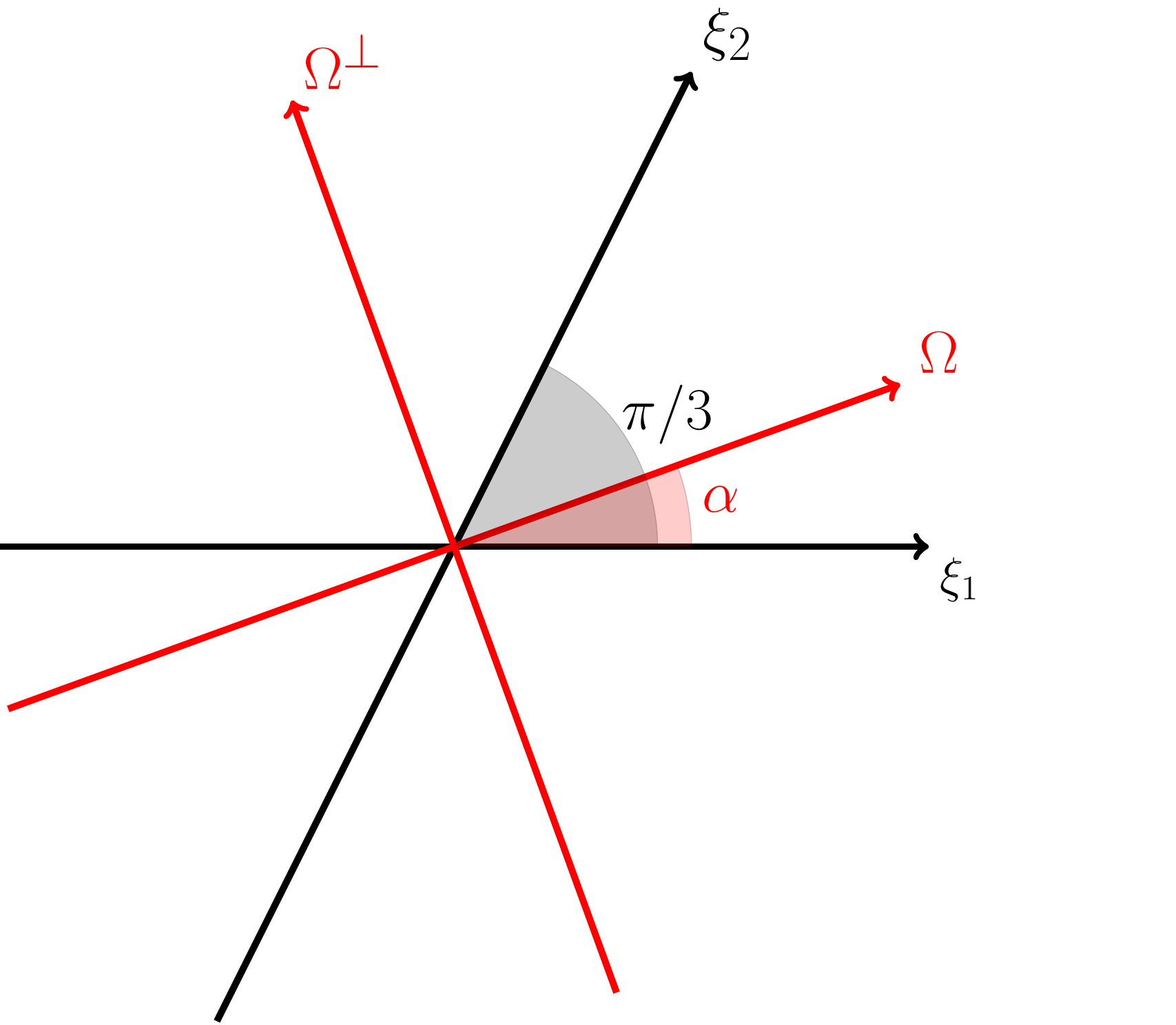}
\caption{Transforming from the $(\bxi_1,\bxi_2)$ coordinate system to 
$(\bOmega,\bOmega^\perp)$.}
\label{fig:transformation}
\end{figure}

The interpolated function value at $\bP_0'$ is then given by
\begin{align*}
	\tilde{f}(\bP_0') = (1-\epsilon \,(c_1(\alpha)+c_2(\alpha))) 
	f(\bP_0)+\epsilon\, 
	c_1(\alpha)  f(\bP_1) + 
	\epsilon\, c_2(\alpha) f(\bP_2).
\end{align*}
%where we omit the dependency of $c_1$ and $c_2$ on $\alpha$ for the sake of 
%notation.
Similarly we compute expressions for $\tilde{f}(\bP_4')$ and 
$\tilde{f}(\bP_5')$.
%\begin{align*}
%f(\bP_4') &= (1-(c_1+c_2)\epsilon) f(\bP_4)+c_1 \epsilon f(\bP_0) + 
%c_2\epsilon 
%%%f(\bP_3), \\
%f(\bP_5') &= (1-(c_1+c_2)\epsilon) f(\bP_5)+c_1 \epsilon f(\bP_6) + 
%c_2\epsilon 
%%%f(\bP_0).
%\end{align*}
If we now reverse the shift by moving all points into the direction $-\bOmega 
\, \epsilon$, we can interpolate a new value for the point $\bP_0$ that is
\begin{align}
	\tilde{f}(\bP_0) =& (1-\epsilon\,(c_1(\alpha)+c_2(\alpha))) f(\bP_0') + 
	 \epsilon \, c_1(\alpha) f(\bP_4') + \epsilon \, c_2(\alpha) f(\bP_5') 
	 \notag  
	 \\ 
%	=&
%	 (1-(c_1+c_2)\epsilon)\left(
%	 (1-(c_1+c_2)\epsilon) f(\bP_0)+c_1 \epsilon f(\bP_1) + c_2\epsilon f(\bP_2)
%	 \right) \\ \notag
%	 &+
%	 c_1\epsilon \left(
%	 (1-(c_1+c_2)\epsilon) f(\bP_4)+c_1 \epsilon f(\bP_0) + c_2\epsilon f(\bP_3)
%	 \right) \\ \notag
%	 &+
%	 c_2\epsilon \left(
%	 (1-(c_1+c_2)\epsilon) f(\bP_5)+c_1 \epsilon f(\bP_6) + c_2\epsilon f(\bP_0)
%	 \right) \\
	 =&f(\bP_0)+
	f(\bP_0) + \epsilon \,\left[
	c_1(\alpha) \left(
	f(\bP_1)-2f(\bP_0)+f(\bP_4)
	\right) +
	c_2(\alpha) \left(
	f(\bP_2)-2f(\bP_0)+f(\bP_5)
	\right)
	\right]%+\mathcal{O}(\epsilon^2),
	\label{eq:stencil}
\end{align}
%or equivalently
%\begin{align*}\label{eq:stencil}%
%		\tilde{f}(\bP_0)  = f(\bP_0) + \epsilon \left(
%		 c_1(\alpha) \left[
%		 f(\bP_1)-2f(\bP_0)+f(\bP_4)
%		 \right] +
%		 c_2(\alpha) \left[
%		 f(\bP_2)-2f(\bP_0)+f(\bP_5)
%		 \right]
%		\right).
%\end{align*}
Note that the exact same result would hold in the case of an equilateral 
triangle with sides of length $\Delta\xi$ when the shift is performed by 
${\bOmega}\, \Delta\xi \, \epsilon$. In order to identify the differential 
operator that the derived stencil approximates, we perform a Taylor expansion 
around 
$\bP_0$, where we use the coordinate axis $\bxi_1$ and $\bxi_2$.
These are the axes that run along the hexagonal grid, centered in $\bP_0$, 
show 
in Fig. \ref{fig:transformation}.
%, which are defined by
%\begin{align*}
%\bm{x} = \bP_0 + (\bP_1-\bP_0)\bxi_1 +  (\bP_2-\bP_0)\bxi_2.
%\end{align*}
From this, we obtain
\begin{align*}
f(\bP_1) &= f(\bP_0) + \frac{\partial}{\partial \bxi_1}f(\bP_0)\Delta 
x +\frac{\partial^2}{\partial \bxi_1^2}f(\bP_0)\frac{\Delta\xi^2}{2}+O(\Delta\xi^3),\\
f(\bP_2) &= f(\bP_0) + \frac{\partial}{\partial 
\bxi_2}f(\bP_0)\Delta\xi 
+\frac{\partial^2}{\partial\bxi_2^2}f(\bP_0)\frac{\Delta\xi^2}{2}+O(\Delta 
\xi^3).
\end{align*}
Plugging this into the derived stencil \eqref{eq:stencil} gives
\begin{align*}
		\tilde{f}(\bP_0)  = f(\bP_0) + \epsilon \Delta\xi^2\left(
		 c_1(\alpha) \left[
		 \frac{\partial^2}{\partial \bxi_1^2}f(\bP_0)
		 \right] +
		 c_2(\alpha) \left[
		 \frac{\partial^2}{\partial \bxi_2^2}f(\bP_0)
		 \right]
		\right)+O(\Delta 
		\xi^3).
\end{align*}

Instead of writing the derivatives in dependency of $\bxi_1$ and $\bxi_2$ we 
transform the derivatives to only rely on the direction $\bOmega$ and the 
direction perpendicular to $\bOmega$, namely $\bOmega^\perp$.

\begin{figure}
\centering
\includegraphics[width=0.5\linewidth]{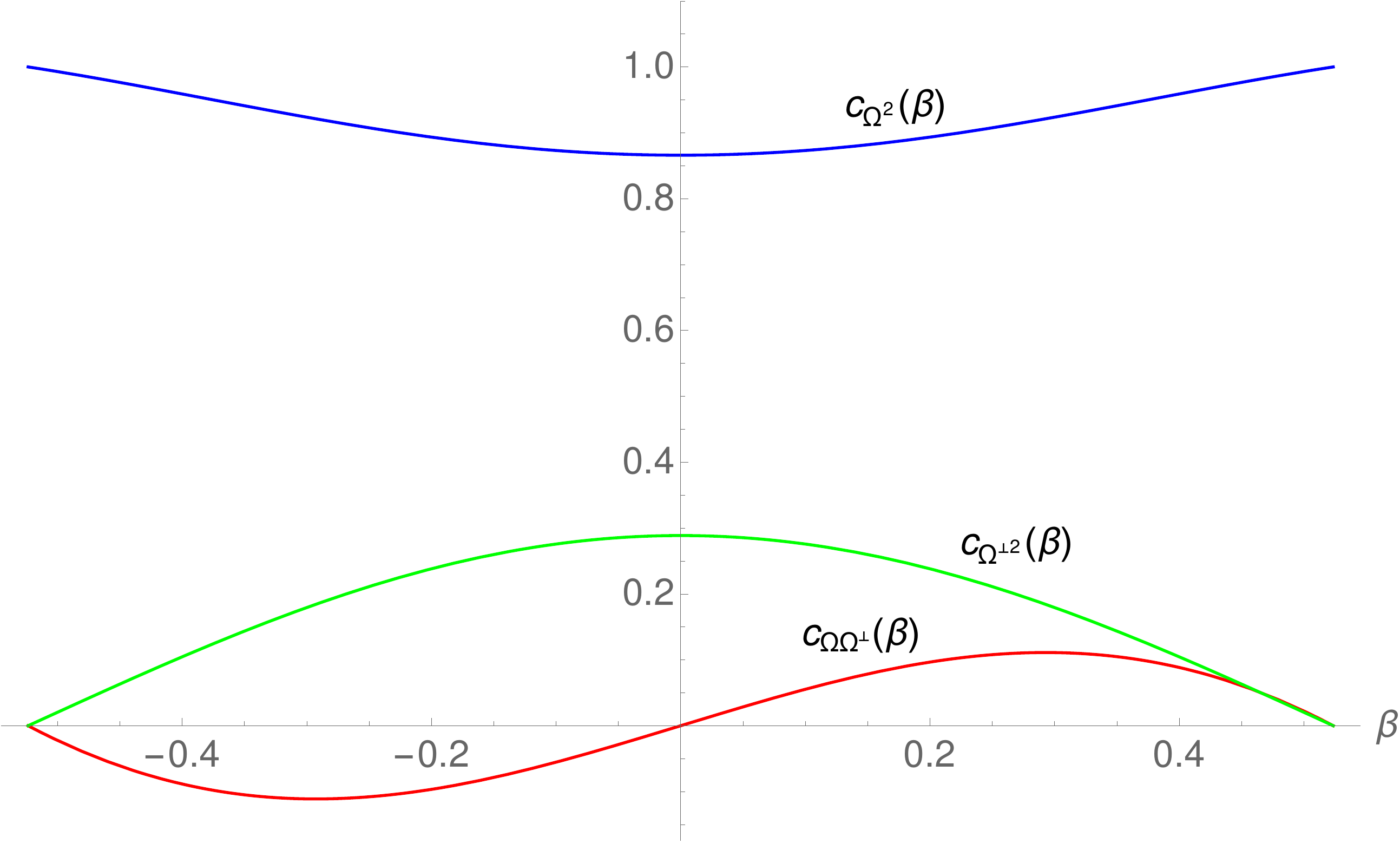}
\caption{Magnitudes of the different second derivative operators in dependency of $\beta$.}
\label{fig:TransformationMagnitude}
\end{figure}

From the geometry of Fig. \ref{fig:transformation} follows
\begin{align*}
&\langle \bOmega,\bxi_1\rangle = \cos(\alpha), \enskip \langle \bOmega,\bxi_2\rangle = \cos(\pi/3-\alpha),\\
&\langle \bOmega^\perp,\bxi_1\rangle =-\sin(\alpha), \enskip \langle \bOmega^\perp,\bxi_2\rangle = \sin(\pi/3-\alpha).
\end{align*}
%\begin{align*}
%	\bOmega &= \langle \bOmega,\bxi_1\rangle \bxi_1 + \langle 
%	\bOmega,\bxi_2\rangle \bxi_2
%			= \cos(\alpha) \bxi_1 + \cos(\pi/3-\alpha)\bxi_2,
%	\\	\bOmega^\perp &= \langle \bOmega^\perp,\bxi_1\rangle \bxi_1 + \langle 
%	\bOmega^\perp,\bxi_2\rangle \bxi_2
%		= -\sin(\alpha)\bxi_1	+\sin(\pi/3-\alpha)\bxi_2.
%\end{align*}
Together with
\begin{align*}
	\frac{\partial^2}{\partial \bxi_i^2}
	=
 \langle \bOmega,\bxi_i\rangle^2\frac{\partial^2}{\partial \bOmega^2}
	+
	2 
	\langle\bOmega,\bxi_i\rangle\langle\bOmega^\perp,\bxi_i\rangle\frac{\partial^2}{\partial
	 \bOmega \partial \bOmega^\perp}
	+
	 \langle \bOmega^\perp,\bxi_i\rangle^2\frac{\partial^2}{\partial 
	 {\bOmega^\perp}^2}
\end{align*}
we obtain
\begin{align*}
c_1(\alpha)\frac{\partial^2}{\partial \bxi_1^2}
&=
\frac{2}{\sqrt{3}}\sin\left(\pi/3-\alpha\right)
\left[
\cos(\alpha)^2\frac{\partial^2}{\partial \bOmega^2}
-
2 \cos(\alpha)\sin(\alpha)\frac{\partial^2}{\partial \bOmega \partial 
\bOmega^\perp}
+
\sin(\alpha)^2\frac{\partial^2}{\partial {\bOmega^\perp}^2}
\right],\\
c_2(\alpha)\frac{\partial^2}{\partial \bxi_2^2}
&=
\frac{2}{\sqrt{3}}\sin\left(\alpha\right)
\left[
\cos(\pi/3-\alpha)^2\frac{\partial^2}{\partial \bOmega^2}
+
2 \cos(\pi/3-\alpha)\sin(\pi/3-\alpha)\frac{\partial^2}{\partial \bOmega 
\partial \bOmega^\perp}
+
\sin(\pi/3-\alpha)^2\frac{\partial^2}{\partial {\bOmega^\perp}^2}
\right].
\end{align*}
Next, we substitute $\alpha = \beta+\pi/6$ with $\beta \in [-\pi/6,\pi/6]$ 
to obtain
\begin{align*}
&c_1(\alpha)\frac{\partial^2}{\partial \bxi_1^2}  + 
c_2(\alpha)\frac{\partial^2}{\partial \bxi_2^2}  
=c_1(\pi/6+\beta)\frac{\partial^2}{\partial \bxi_1^2}  + 
c_2(\pi/6+\beta)\frac{\partial^2}{\partial \bxi_2^2} 
%&=
%\frac{2}{\sqrt{3}}\Bigg(
%{\color{blue}
%\left[\sin(\pi/6-\beta)\cos(\pi/6+\beta)^2 + 
%%%\sin(\pi/6+\beta)\cos(\pi/6-\beta)^2\right]}\frac{\partial^2}{\partial 
%%%%\bOmega^2} \notag\\
%&
%+{\color{red}2\left[
%-\sin(\pi/6-\beta)\cos(\pi/6+\beta)\sin(\pi/6+\beta)
%+\sin(\pi/6+\beta)\cos(\pi/6-\beta)\sin(\pi/6-\beta)
%\right]}\frac{\partial^2}{\partial \bOmega \partial \bOmega^\perp} \notag \\
%&+
%{\color{green}\left[\sin(\pi/6-\beta)\sin(\pi/6+\beta)^2+\sin(\pi/6+\beta)\sin(\pi/6-\beta)^2\right]}
%\frac{\partial^2}{\partial {\bOmega^\perp}^2}\Bigg).
\\&=
c_{\bOmega^2}(\beta) \frac{\partial^2}{\partial \bOmega^2}
+
c_{\bOmega\bOmega^\perp}(\beta)\frac{\partial^2}{\partial \bOmega 
\partial \bOmega^\perp} 
+
c_{{\bOmega^\perp}^2}(\beta)\frac{\partial^2}{\partial 
{\bOmega^\perp}^2}
.
\end{align*}
Here, we defined the following constants
\begin{align*}
c_{\bOmega^2}(\beta) &=
\frac{1}{2\sqrt{3}}(4 \cos (\beta )-\cos (3 \beta )), \\
c_{\bOmega\bOmega^\perp} (\beta)&
=\frac{1}{2\sqrt{3}} \left( -4\sin(\beta)+2 \sin(3\beta)\right),\\
c_{{\bOmega^\perp}^2}(\beta) &
=\frac{1}{2\sqrt{3}} \cos(3\beta),
\end{align*} 
which are visualized in Fig. \ref{fig:TransformationMagnitude}.

%\end{align}

%The expressions can be simplified and finally yield
%\begin{align}
%&c_1(\alpha)\frac{\partial^2}{\partial \bxi_1^2}  + 
%%%c_1(\alpha)\frac{\partial^2}{\partial \bxi_1^2}  
%%%%=c_1(\pi/6+\beta)\frac{\partial^2}{\partial \bxi_1^2}  + 
%%%%%c_1(\pi/6+\beta)\frac{\partial^2}{\partial \bxi_1^2} \notag\\
%&=
%\frac{2}{\sqrt{3}}\Bigg(
%{\color{blue}
%\left[
%\frac{1}{4} (4 \cos (\beta )-\cos (3 \beta ))\right]	}
%\frac{\partial^2}{\partial \bOmega^2} \notag\\
%&
%+{\color{red}2\left[\frac{1}{4} (\sin (3 \beta )-2 \sin (\beta ))	\right]}
%\frac{\partial^2}{\partial \bOmega \partial \bOmega^\perp} \notag \\
%&+
%{\color{green}\left[%
%	\frac{1}{4} \cos (3 \beta )\right]}
%\frac{\partial^2}{\partial {\bOmega^\perp}^2}\Bigg)
%\\
%&=
%{\color{blue}c_{\bOmega^2}(\beta)} \frac{\partial^2}{\partial \bOmega^2}
%	+
%	{\color{red}c_{\bOmega\bOmega^\perp}(\beta)}\frac{\partial^2}{\partial 
%\bOmega 
%\partial \bOmega^\perp} (\beta)
%		+
		%{\color{green}c_{{\bOmega^\perp}^2}(\beta)}\frac{\partial^2}{\partial 
		%{\bOmega^\perp}^2}
%		.
%\end{align}
%where we defined the following constants
%\begin{align}
%	{\color{blue}c_{\bOmega^2}(\beta)} &{\color{blue}= }
%	{\color{blue}\frac{1}{2\sqrt{3}}(4 \cos (\beta )-\cos (3 \beta ))}, \\
%		{\color{red}c_{\bOmega\bOmega^\perp} (\beta)}&
%	{\color{red}=\frac{1}{2\sqrt{3}} \left( -4\sin(\beta)+2 
%\sin(3\beta)\right)},\\
%		{\color{green}c_{{\bOmega^\perp}^2}(\beta) }&
%	{\color{green}=\frac{1}{2\sqrt{3}} \cos(3\beta)},
%\end{align}
%corresponding to the colored expressions in the equations above and visualize 
%them in Fig. \ref{fig:TransformationMagnitude}.

Finalizing the change of coordinate systems, we derive
\begin{align}
\tilde{f}(\bP_0)  &= f(\bP_0) + \epsilon \left(
c_1(\alpha) \left[
\frac{\partial^2}{\partial \bxi_1^2}f(\bP_0)\Delta\xi^2
\right] +
c_2(\alpha) \left[
\frac{\partial^2}{\partial \bxi_2^2}f(\bP_0)\Delta\xi^2
\right]
\right)+O(\Delta 
\xi^3) \notag \\
&=
f(\bP_0) + \epsilon \Delta\xi^2 \left(c_{\bOmega^2}(\beta) 
\frac{\partial^2}{\partial \bOmega^2}
+
c_{\bOmega\bOmega^\perp}(\beta)\frac{\partial^2}{\partial \bOmega 
\partial \bOmega^\perp} 
+
c_{{\bOmega^\perp}^2}(\beta)\frac{\partial^2}{\partial 
{\bOmega^\perp}^2}
\right)f(\bP_0) +O(\Delta 
\xi^3).
%\\
%&=
%f(\bP_0) + \frac{\epsilon \Delta\xi^2}{2\sqrt{3}} \left((4 \cos 
%(\beta )-\cos (3 \beta )) \frac{\partial^2}{\partial \bOmega^2}
%+
%\left( -4\sin(\beta)+2 
%\sin(3\beta)\right)\frac{\partial^2}{\partial \bOmega \partial \bOmega^\perp} 
%(\beta)
%+
%\cos(3\beta)\frac{\partial^2}{\partial {\bOmega^\perp}^2}
%\right)f(\bP_0).
\label{eqmod}
\end{align}

We are now going to draw some conclusions from this derivation. First, we 
observe a diffusive behavior, where diffusion is strongest along the direction 
of shifting. Furthermore, when the shift is performed in alignment with the 
lattice, diffusion only occurs along that direction. 

Comparing equation \eqref{eqmod} with \eqref{eq2dmod} we observe a similar 
scaling behavior of the diffusion term. Choosing $\epsilon=\delta \, \Delta  t 
/ \Delta  \xi$
again results in vanishing diffusion for a refinement of the angular 
discretization ($\Delta \xi\rightarrow 0$). Since $\Delta \xi$ scales like the 
number of angular 
quadrature points $N_q$, we will perform rotations by $\delta \Delta t / N_q$ 
instead of $\delta \Delta t / \Delta \xi$. This is due to the fact that $\Delta 
\xi$ is not constant for the different triangles on the sphere, it does however 
scale like $N_q$.

%Compared to \eqref{eq2dmod}, \dots \mf{how to choose magnitude of rot angle 
%depending on time step and angular discretization? Check notation for 
%consistency: $\alpha$ is used twice in different meaning. Relationship 
%$\alpha,\epsilon,c,\dots$}

The implementation of the \rSN method differs from this simplified 
analysis. There, we are moving the quadrature points on the sphere and not in 
planar geometry. Furthermore, not all triangles have the same size. 
Additionally, determining the corresponding points from which we interpolate 
the function values is not trivial. It might happen, that two rotated points 
fall into the same triangle of the old quadrature set.
In the implementation of the \rSN method, we also rotate randomly 
around different axes and perform the usual \SN update, i.e.\ stream and 
collide, in between interpolating. All these aspects make the theoretical 
analysis significantly more difficult than for the simplified planar geometry.
However, the numerical results indicate a diffusive behavior which is equally 
strong in any direction. We believe that randomly choosing a rotation axis has 
an averaging effect. That is, diffusion will be equally strong along all 
directions since we rotate differently in each time step.
%The observable diffusive behavior indicates a similar modified equations as 
%derived in Eq. \eqref{eq:modequation} 
As we are no longer restricted to rotations around the z-axis only, diffusion 
is also not restricted to occur in the azimuthal angle. 

A different way to interpret the newly induced diffusion is artificial scattering. 
Rotating the quadrature set and interpolating the angular flux can be seen as 
particles scattering into new directions.
As it is the case for scattering, the rotation and interpolation procedure is 
conservative by construction. Furthermore, as described in Alg.~3, we can write 
the interpolation procedure as $\tilde{\psi}_{i,j} = I \psi_{i,j}$. Here, 
$\psi_{i,j}$ is the angular flux in a given spatial cell $c_{i,j}$ before the 
rotation and interpolation step and $\tilde{\psi}_{i,j}$ is the angular flux at 
the new quadrature points after the rotation and interpolation step. The matrix 
$I$ contains the interpolation weights and has only three non zero elements per 
row. It does not depend on the spatial cell index, but is different for each 
time step since the rotation axis differs from time step to time step.
Thus, the angular flux at $\Omega_q$ "scatters" into the directions for which 
it is used to interpolate function values at the new quadrature set.
Since more scattering implies fewer ray-effects, the \rSN method 
mitigates these undesired effects.

\section{Implementation}
\label{sec:Implementation}
The method can easily be implemented into an existing code for the discrete 
ordinates method as it is minimally inversive. Only the construction of the 
quadrature set has to be modified and a function for the rotation and 
interpolation has to be implemented. In Alg.~\ref{alg:sn} we see the main 
components of a discrete ordinates method. We assume available routines that 
return the quadrature points and quadrature weights given an order $N$, 
indicated by \texttt{getOctahedronQuadraturePoints(N)} and 
\texttt{getOctahedronQuadratureWeights(N)}.
These implementations are in accordance with the explanation in Section 
\ref{sec:quadpoints} and Section \ref{sec:quadweights}. 
For the 
standard \SN method, the number of quadrature points in the angular variable 
scales as $N_q=N^2$. Inside the while loop, the angular flux $\psi$ is updated 
via a finite 
volume scheme as described in Section \ref{sec:FVM}.
The modifications that have to be implemented to obtain the \rSN 
method are then highlighted in blue in Alg.~\ref{alg:rsn}.

\begin{minipage}[t]{7cm}
	\vspace{0pt}
	\centering 
\begin{algorithm}[H]	
	\underline{function \SN}$(\Delta t,t_{end},N,N_x,N_y,\psi_0) $	\;
	%\SetKwInOut{Input}{Input}
	%\SetKwInOut{Output}{Output}
	%\Input{Two nonnegative integers $a$ and $b$}
	%\Output{$\gcd(a,b)$}
	$ $
	
	$ t=0 $ 
	
	$N_q = N^2$
	
	$Q = \text{getQuadraturePoints}(N) \in \mathbb{R}^{3 \times N_q }$
	
	$W = \text{getQuadratureWeights}(N) \in \mathbb{R}^{N_q }$

	$\psi=\psi_0%\text{getInitialCondition}(N_q,N_x,N_y) 
	\in \mathbb{R}^{N_q 
	\times 
N_x \times N_y}$

	\While{$t<t_{end}$}{
	 $F = \text{computeFluxes}(\psi,Q,W) $
	 
	 $\psi = \psi+\Delta t\cdot F$
	 
	 $ $
	 
	 $t = t + \Delta t$
	}	
	
	$ $ 
	
	return $\psi$
    \caption{The \SN method}
	\label{alg:sn}
\end{algorithm}
\end{minipage}
\hspace{1cm}
\begin{minipage}[t]{9cm}
	\centering
	\vspace{0pt}
	\begin{algorithm}[H]
			\underline{function 
			\rSN}  $(\Delta t,t_{end},N,N_x,N_y,\psi_0,{\color{blue}\delta}) $	
	\;
	%\SetKwInOut{Input}{Input}
	%\SetKwInOut{Output}{Output}
	%\Input{Two nonnegative integers $a$ and $b$}
	%\Output{$\gcd(a,b)$}
	$ $
	
	$ t=0 $ 
	
	{\color{blue}$N_q = 4N^2-8N+6$}
	
	{\color{blue}$Q = \text{getOctahedronQuadraturePoints}(N) \in 
	\mathbb{R}^{3\times N_q }$}
	
{\color{blue}$W = \text{getOctahedronQuadratureWeights}(N) \in \mathbb{R}^{N_q 
}$}

	$\psi=\psi_0%\text{getInitialCondition}(N_q,N_x,N_y) 
\in \mathbb{R}^{N_q 
	\times 
	N_x \times N_y}$

	\While{$t<t_{end}$}{
		$F = \text{computeFluxes}(\psi,Q,W) $
		
		$\psi = \psi+\Delta t\cdot F$
		
		{\color{blue} $\psi,Q = \text{rotateAndInterpolate}(\psi,Q,\delta\cdot 
		\Delta t / N_q)$}
		
		$t = t + \Delta t$
	}	
	$ $ 
	
	return  $\psi$
	\caption{The \rSN method}
	\label{alg:rsn}
\end{algorithm}
\end{minipage}

%\end{minipage}
The algorithm that performs the rotation and interpolation in line 11 is 
presented in Alg.~\ref{alg:rotinterp}.
Here we assume that a method to determine the triangle into 
which a rotated quadrature point falls into is available. Such a method can easily be implemented in 
an efficient way as we can assume that any point will fall into one of its six 
neighboring triangles.
After having found the corresponding triangle, the interpolation weights are 
being computed and stored in a matrix $W$ as explained in Alg.~\ref{alg:rotinterp} in line 11 to 14. 
After having computed the 
interpolation weights, the interpolation is applied in each spatial cell. We 
only need to compute the interpolation weights once since the quadrature set is 
the same in each spatial cell. Alg.~\ref{alg:rotinterp} is furthermore well 
suited for parallel implementations.

In general, we cannot expect the interpolation procedure to be conservative.
It does however preserve positivity. One way to force the method to be 
conservative, is by rescaling the scalar flux at each quadrature 
point, such that the mass before the interpolation step is conserved in each 
cell separately.

In our comparisons, \rSN with $\delta=0$ is different from traditional \SN 
(relying on a tensorized angular grid as described in Section 
\ref{sec:RotationSimple}), because we use different quadrature sets.
We have included both in the comparisons to distinguish the effect of the new 
quadrature from the rotation procedure. Furthermore, due to the different 
relation between $N$ and 
$N_q$, we will later on compare the methods according to their total number of quadrature 
points $N_q$ and not their order $N$.

The rotation magnitude $\delta$ is scaled by $\Delta t/N_q$, shown in line 11 
of 
Alg.~\ref{alg:rsn}. This makes the observable effect of the rotation not depend 
on different time step sizes.

Since the matrix $W$  in Alg.~\ref{alg:rotinterp} contains only three 
nonzero entries per row, the procedure can also be implemented in a sparse 
manner. For readability, we chose not to do so in this example.

In all numerical experiments, we observed an increase of the runtime by $5\%$ 
to $10\%$ when adding the rotation and interpolation procedure to the standard 
\SN implementation. Within the context of all performed simulations we use the 
implementation of the \rSN method as described in Alg.~\ref{alg:rsn}. To 
compare simulations with different rotations strengths $\delta$ we denote by 
r${}_\delta$S$_N$ the \rSN method with rotation strength $\delta$.

	\begin{algorithm}[H]
	\underline{function rotateAndInterpolate}$(\psi,Q,\delta) $	
	\;
	%\SetKwInOut{Input}{Input}
	%\SetKwInOut{Output}{Output}
	%\Input{Two nonnegative integers $a$ and $b$}
	%\Output{$\gcd(a,b)$}
	$ $
	
	$N_q,N_x,N_y= \text{size}(\psi)$
	
	$n = \text{getRandomAxis}() \in \mathbb{S}^2$ \hfill The axis we are going 
	to rotate around.
	
	$R = \text{getRotationMatrix}(n,\delta) \in \mathbb{R}^3$  \hfill  See Eq. 
	\eqref{eq:rotmatrix}.
	
	$\hat{Q} = R \cdot Q \in \mathbb{R}^{3 \times N_q}$   \hfill  The rotated 
	quadrature 
	points.
	
	$W = \text{zeros}(N_q,N_q) \in \mathbb{R}^{N_q \times N_q}$ \hfill	
	The matrix 
	$W$ will store interpolation weights.
	
		$\hat{\psi} = \text{zeros}(N_q,N_x,N_y) \in \mathbb{R}^{N_q \times N_x 
		\times N_y}$ \hfill  The tensor $\hat{\psi}$ will store the 
		interpolated angular flux.
	
	\For{$q=1,\cdots,N_q$}{
		$\hat{q} = \hat{q}[:,q]$  \hfill  Store in 
		$\hat{q}$ a quadrature 
		point from the rotated quadrature.
		
		$i,j,k = \text{interpolateFrom}(Q,\hat{q})$ \hfill  Compute the three 
		vertices 
		of the triangle that
		
		$ $ \hfill   the 
		rotated quadrature point $\hat{q}$ falls into. 
		
		$w_i,w_j,w_k = \text{interpolationWeights}(Q,\hat{q})$ \hfill  
		The interpolation weights as 
		described 
		in Section \ref{sec:quadinterp}.
		
		$W[q,i] = w_i$,\, $W[q,j] = w_j$,\, $W[q,k] =w_k$ \hfill Store the 
		interpolation weights in $W$.
	}

\For{$i=1,\cdots,N_x$}{
	\For{$j=1,\cdots,N_y$}{
		
		$\hat{\psi}[:,i,j] = W \cdot \psi[:,i,j]$ \hfill Apply interpolation in 
		each spatial cell.
		
	}
}

$ $ 
	
	return  $\hat{\psi}$, $\hat{Q}$
	\caption{The rotation and interpolation routine}
	\label{alg:rotinterp}
\end{algorithm}

\subsection{Alternative implementations}
Implementing the rotation and interpolation step into the standard \SN method 
is possible in several ways. We will now describe two of these modified 
implementations. All numerical simulations were performed with the method 
described above. This is due to the fact that all tested implementations show 
similar qualitative behavior. This indicates the importance of the presence of 
the interpolation and rotation step. For example the order in which these steps 
are being executed as well as other details play a smaller role.

\paragraph*{Rotating forth and back}
One implementation performed a rotation around a random axis with strength 
$\delta$ in one time step, followed by a rotation with strength $-\delta$ 
around the same axis in the next time step. This procedure is similar to the 
idea described in Sec.~\ref{sec:modeq} and allowed to perform a more rigorous 
theoretical analysis, but showed no differences in the computed solution.

\paragraph*{Two opposite rotations within one time step}
The proposed implementation requires to update the quadrature set outside of 
the rotation and interpolation step. There might be implementations of the \SN 
method were this is not possible. In those cases, two opposite rotations can be 
performed within one time update. This keeps the quadrature set outside the 
rotation and interpolation procedure fixed. To observe the same qualitative 
behavior as for the proposed \rSN implementation, two rotations with 
$\delta/2$ and $-\delta/2$ have to be performed instead of one rotation with 
 $\delta$.

\section{Numerical results}
\label{sec:Results}
In the following, we show results obtained with the standard \SN method, as 
well as with the \rSN method. Our results can be reproduced with the freely avaliable code \cite{JuliaWN}. To demonstrate the properties of 
ray effect mitigation, we choose test cases that are prone to this behavior.

\subsection{The line-source problem}
\label{sec:LineSourceResults}
The first problem we look at is the \emph{line-source problem}. In this test 
case, the solution is projected onto the $xy$-plane, i.e. the spatial domain is 
$[a,b]\times[a,b]$. There is only one medium with cross sections 
$\sigma_a,\sigma_b$. The initial density distribution is a dirac in the center 
of the spatial domain $x=0$, $y=0$, i.e. we have
\begin{align*}
\psi(0,\bx,\bOmega) = \frac{1}{4\pi}\delta_0(x)\delta_0(y)\;.
\end{align*}
Numerically, the initial condition is approximated by
\begin{align}\label{eq:IClineSource}
\psi(0,\bx,\bOmega) = \frac{1}{4\pi\sigma_{\text{IC}}^2} 
\cdot\exp\left(-\frac{x^2+y^2}{4\sigma_{\text{IC}}^2}\right)
\end{align}
with $\sigma_{\text{IC}}$ close to zero (see Tab.~\ref{tab:linesource}). We vary the number of quadrature 
points as well as the rotation strength $\delta$ to study effects on the 
solution. We use a second-order scheme with a minmod slope limiter in space as well as 
Heun's method in time. The code used in this work is a re-implementation of \cite{garrett2013comparison}.

\begin{figure}
	
	\begin{subfigure}{0.24\linewidth}
		\centering
		\includegraphics[scale=0.17]{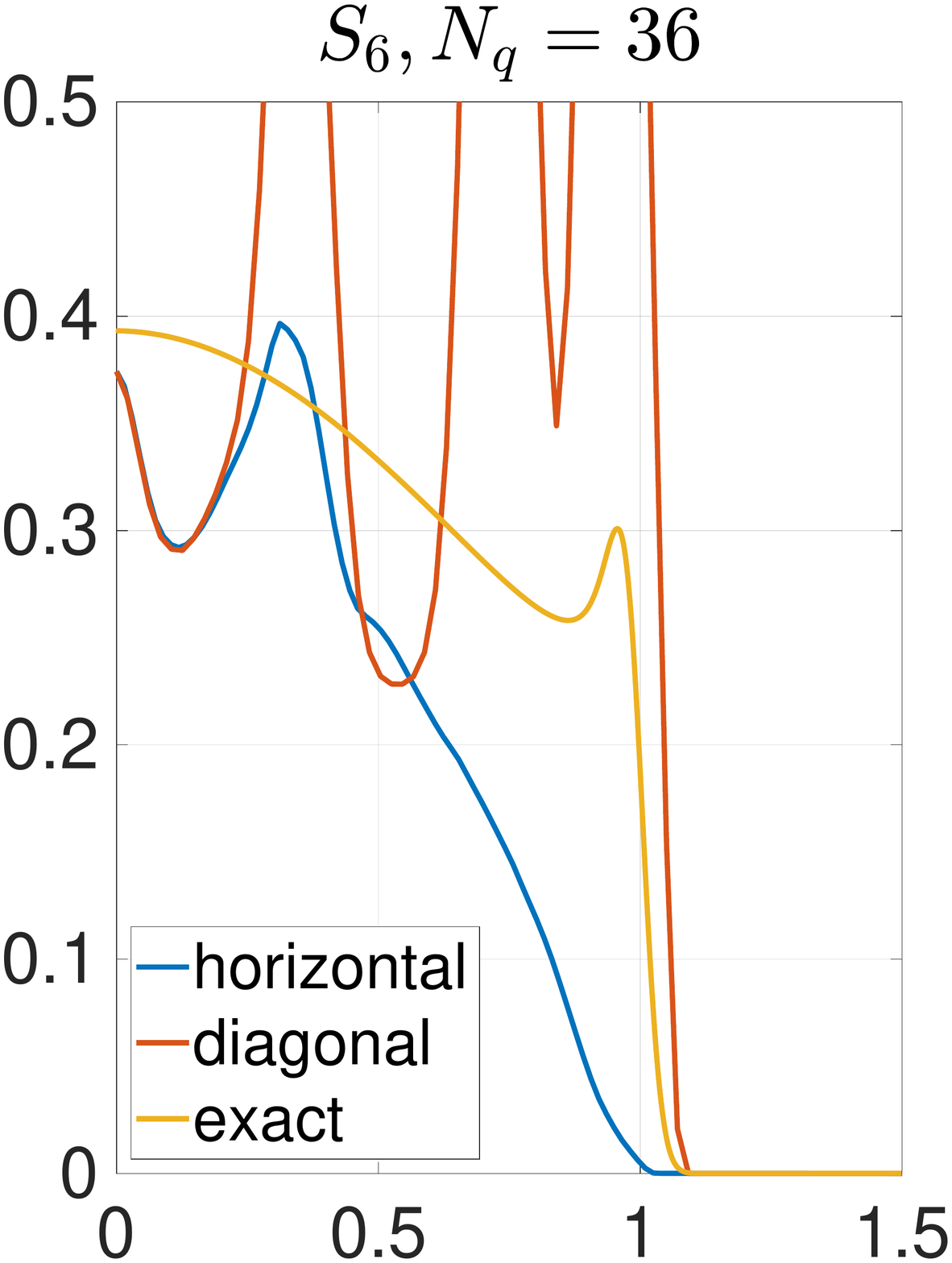}
		
		\label{fig:sub1}
	\end{subfigure}%
	\begin{subfigure}{0.24\linewidth}
		\centering
		\includegraphics[scale=0.17]{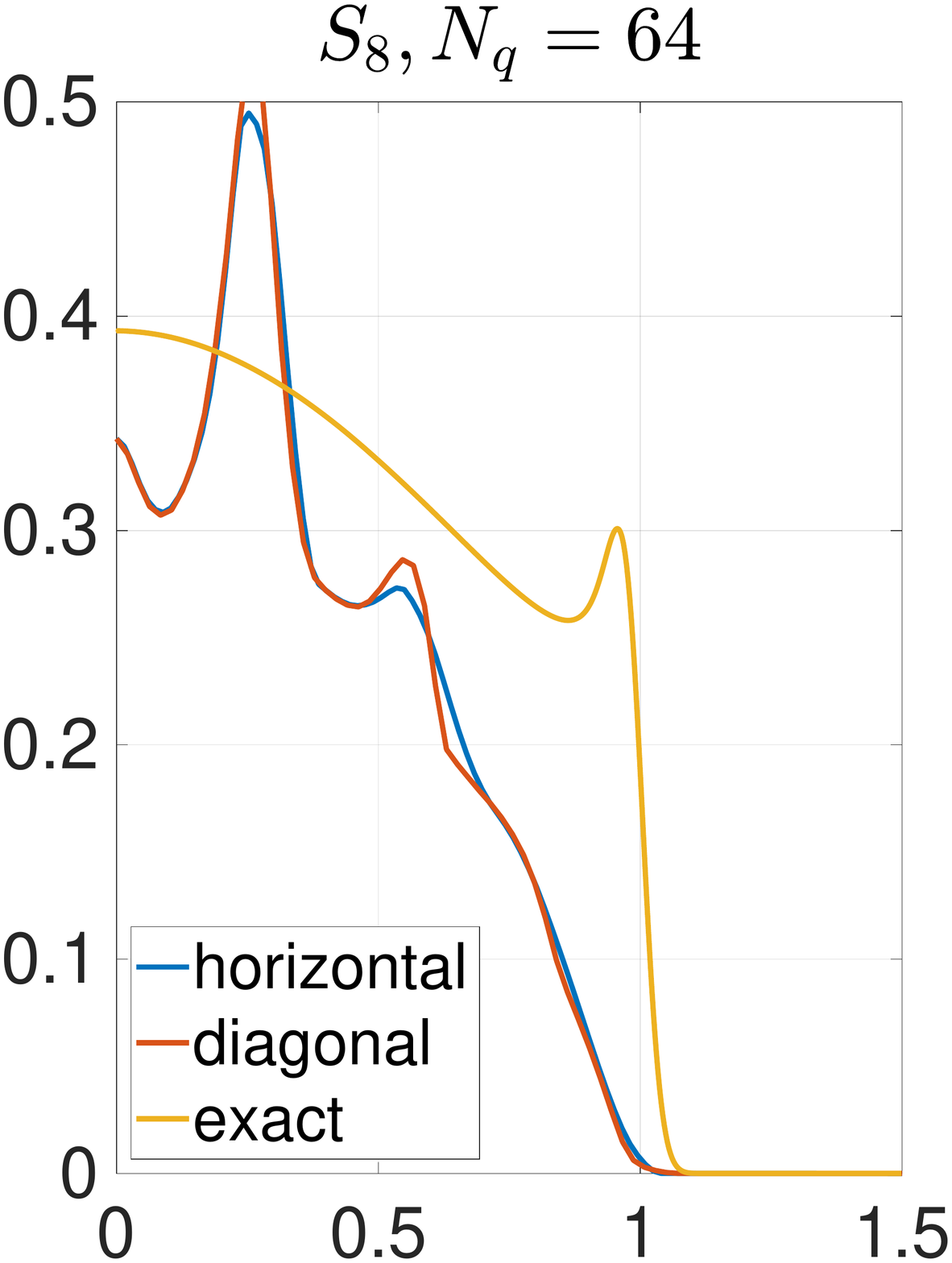}
		
		\label{fig:sub2}
	\end{subfigure}
	\begin{subfigure}{0.24\linewidth}
		\centering
		\includegraphics[scale=0.17]{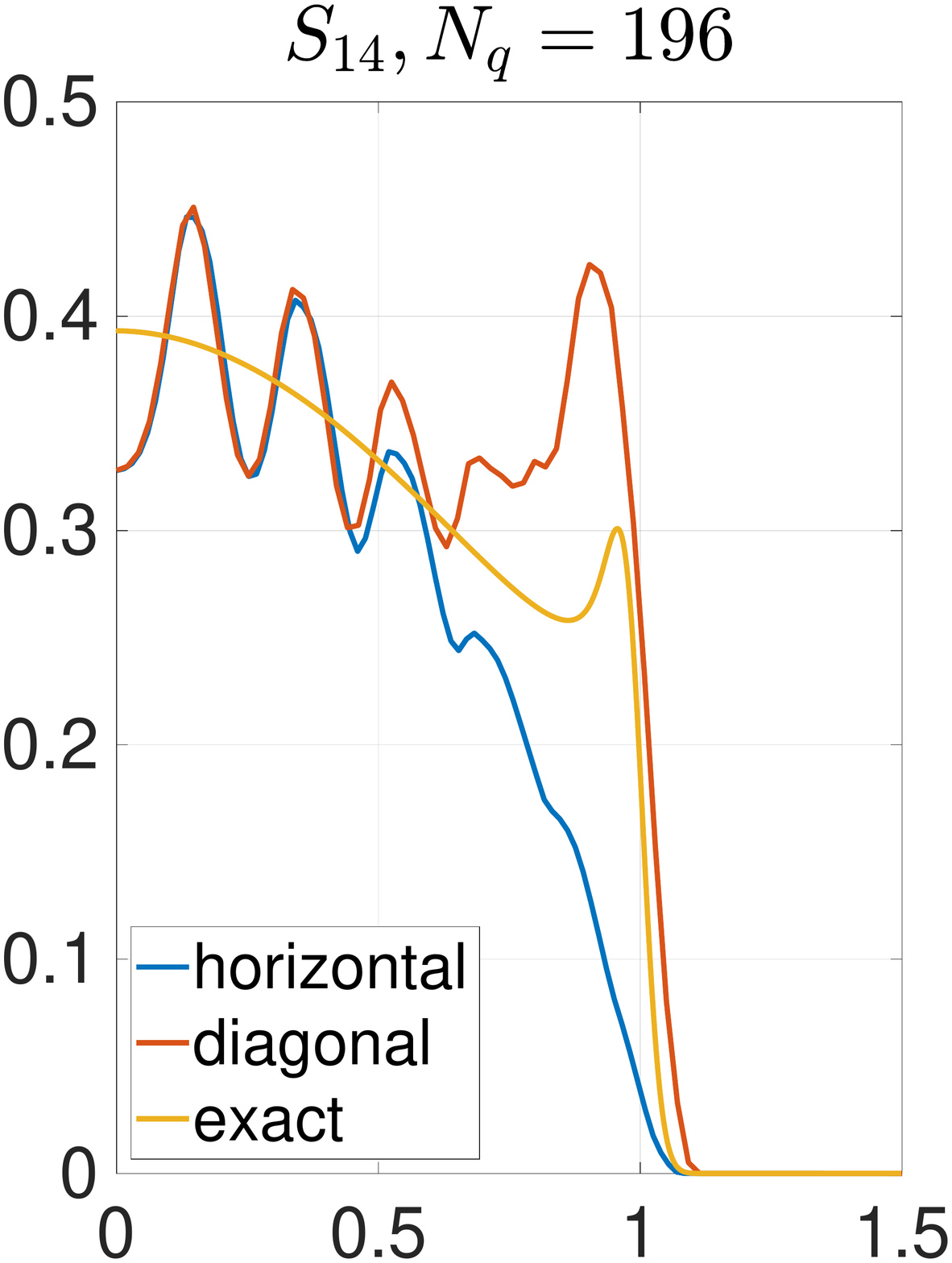}
		
		\label{fig:sub3}
	\end{subfigure}
	\begin{subfigure}{0.24\linewidth}
		\centering
		\includegraphics[scale=0.17]{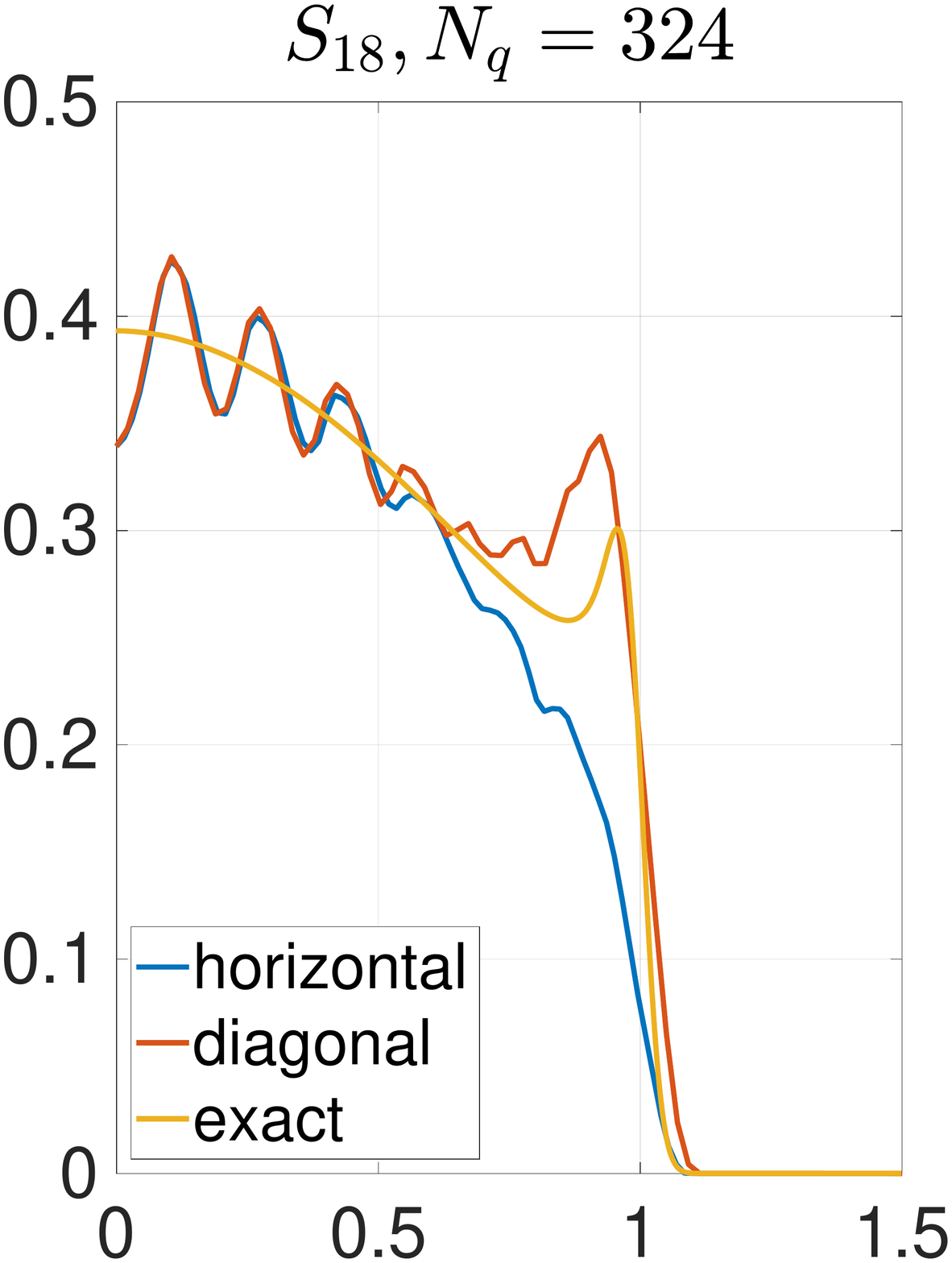}
		
		\label{fig:sub3}
	\end{subfigure}\\[-1ex]
	\begin{subfigure}{0.24\linewidth}
		\centering
		\includegraphics[scale=0.17]{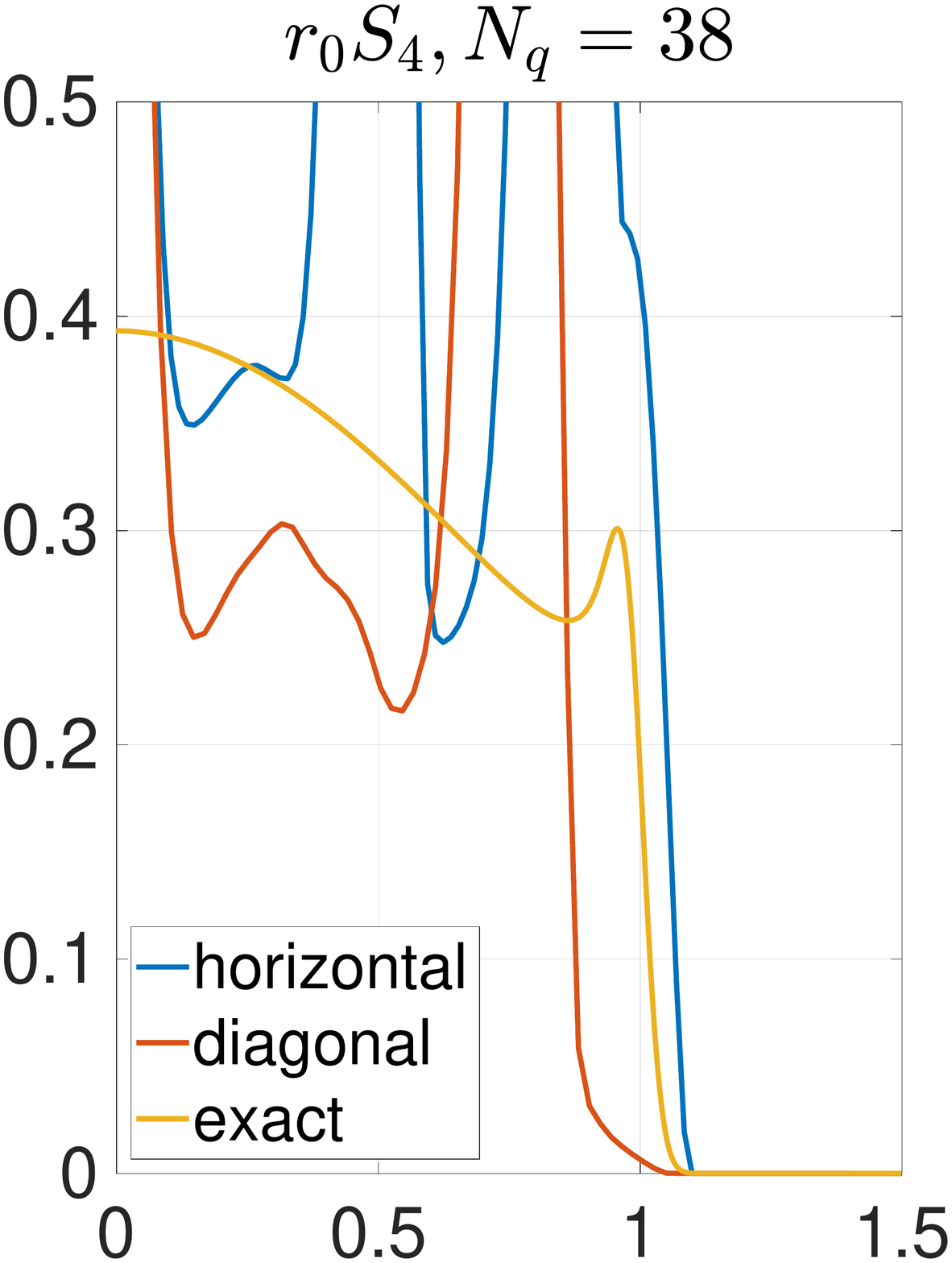}
		
		\label{fig:sub1}
	\end{subfigure}%
	\begin{subfigure}{0.24\linewidth}
		\centering
		\includegraphics[scale=0.17]{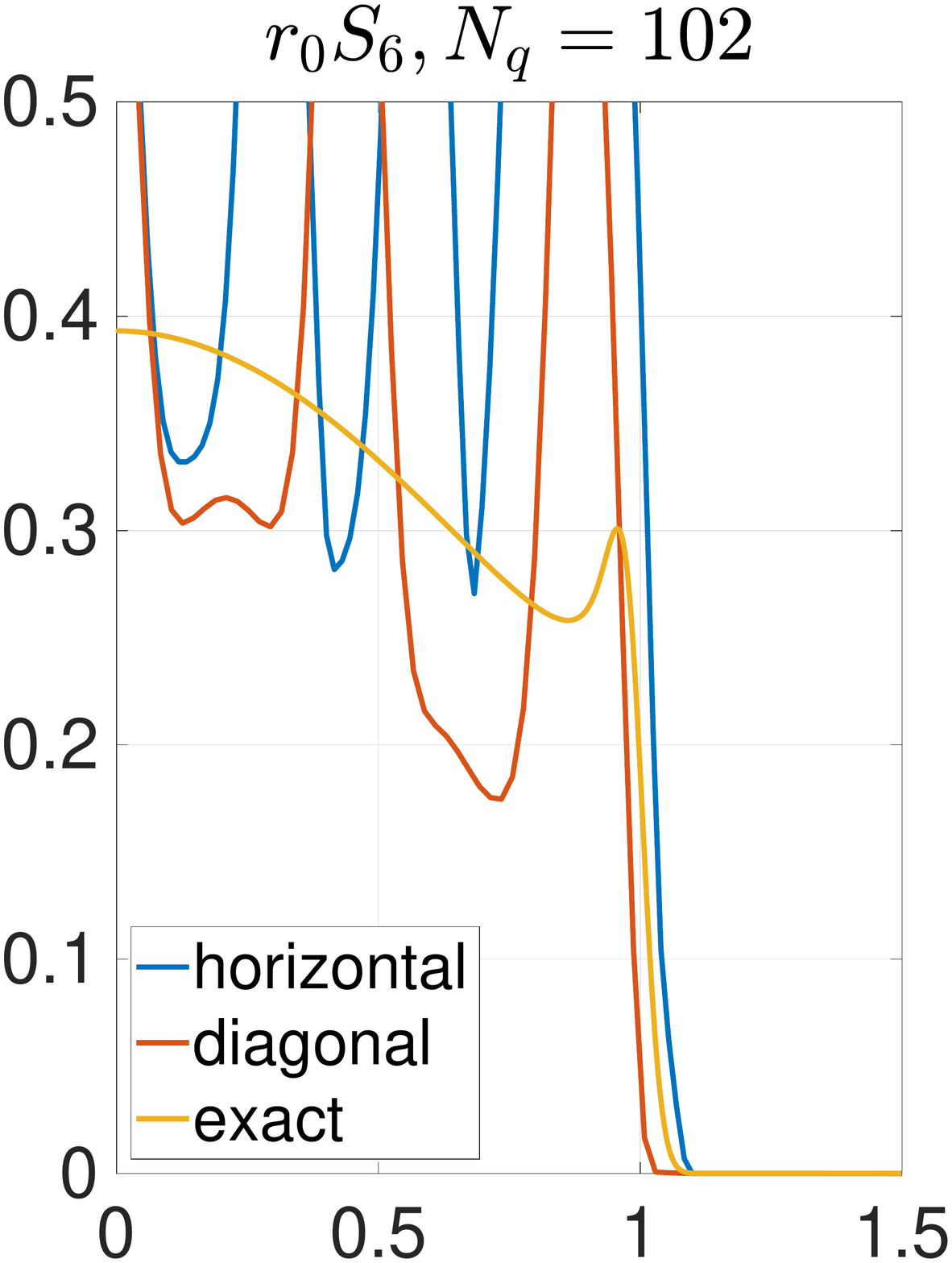}
		
		\label{fig:sub2}
	\end{subfigure}
	\begin{subfigure}{0.24\linewidth}
		\centering
		\includegraphics[scale=0.17]{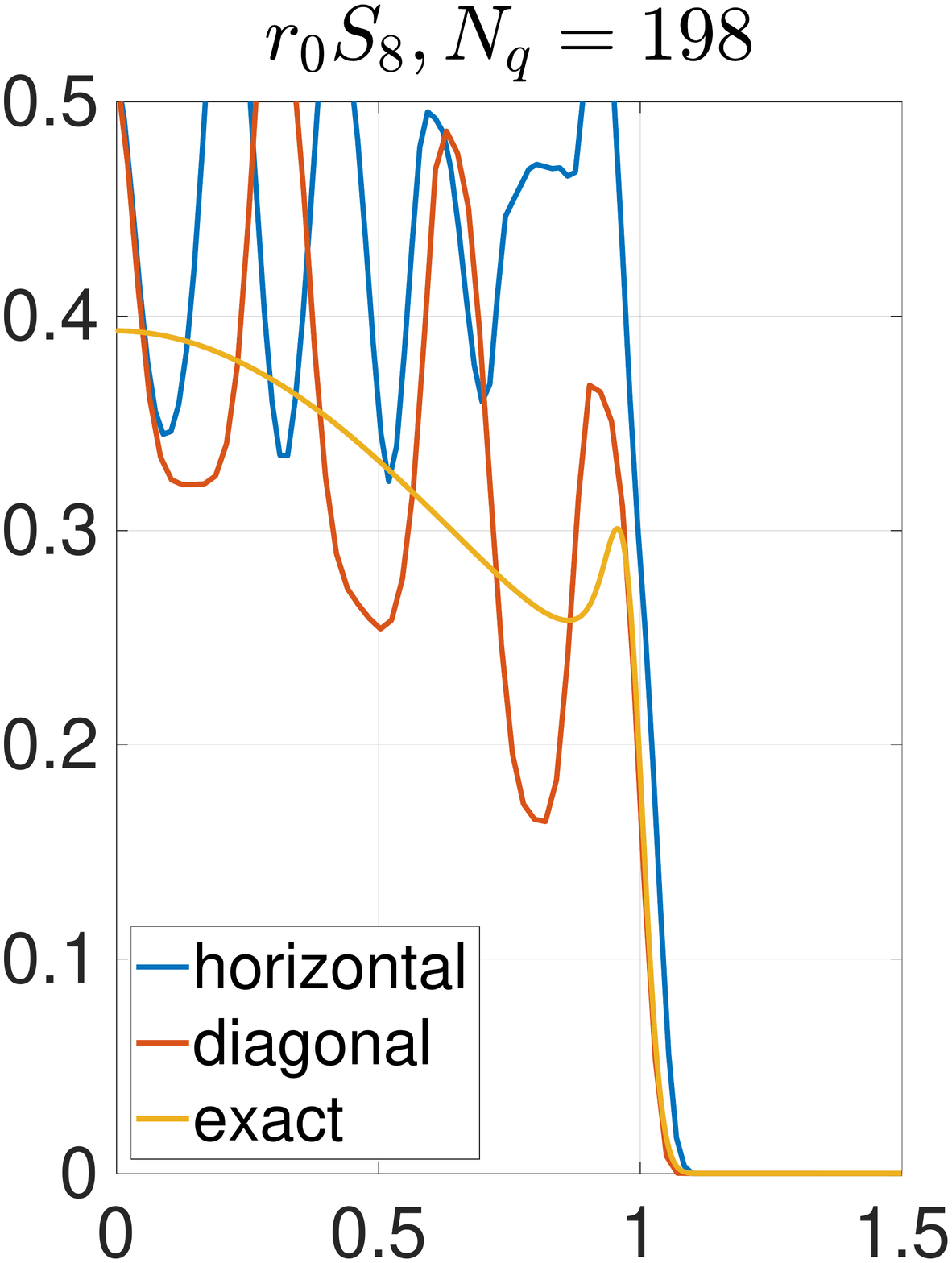}
		
		\label{fig:sub3}
	\end{subfigure}
	\begin{subfigure}{0.24\linewidth}
		\centering
		\includegraphics[scale=0.17]{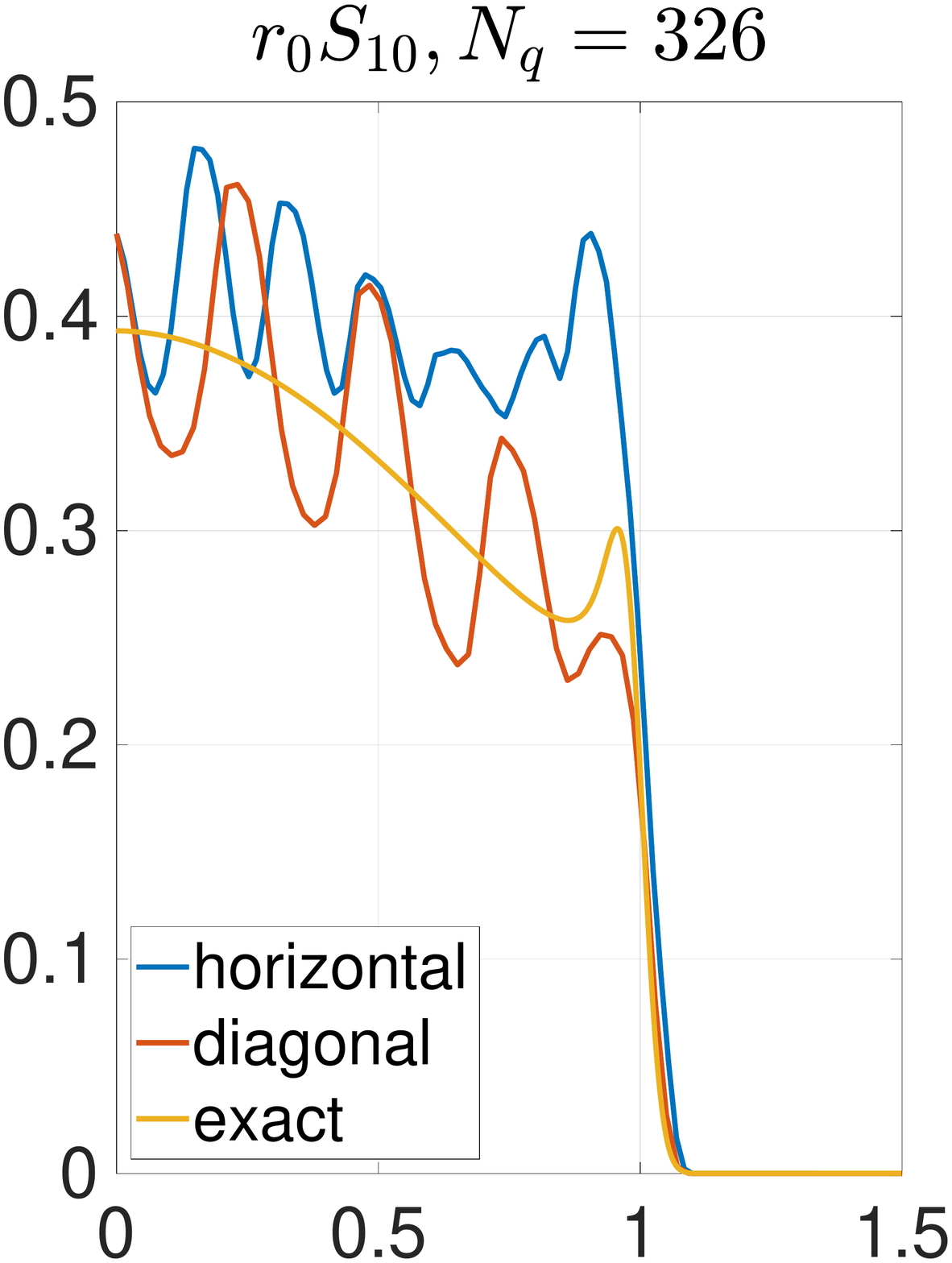}
		
		\label{fig:sub3}
	\end{subfigure}
	\\[-1ex]
	\begin{subfigure}{0.24\linewidth}
		\centering
		\includegraphics[scale=0.17]{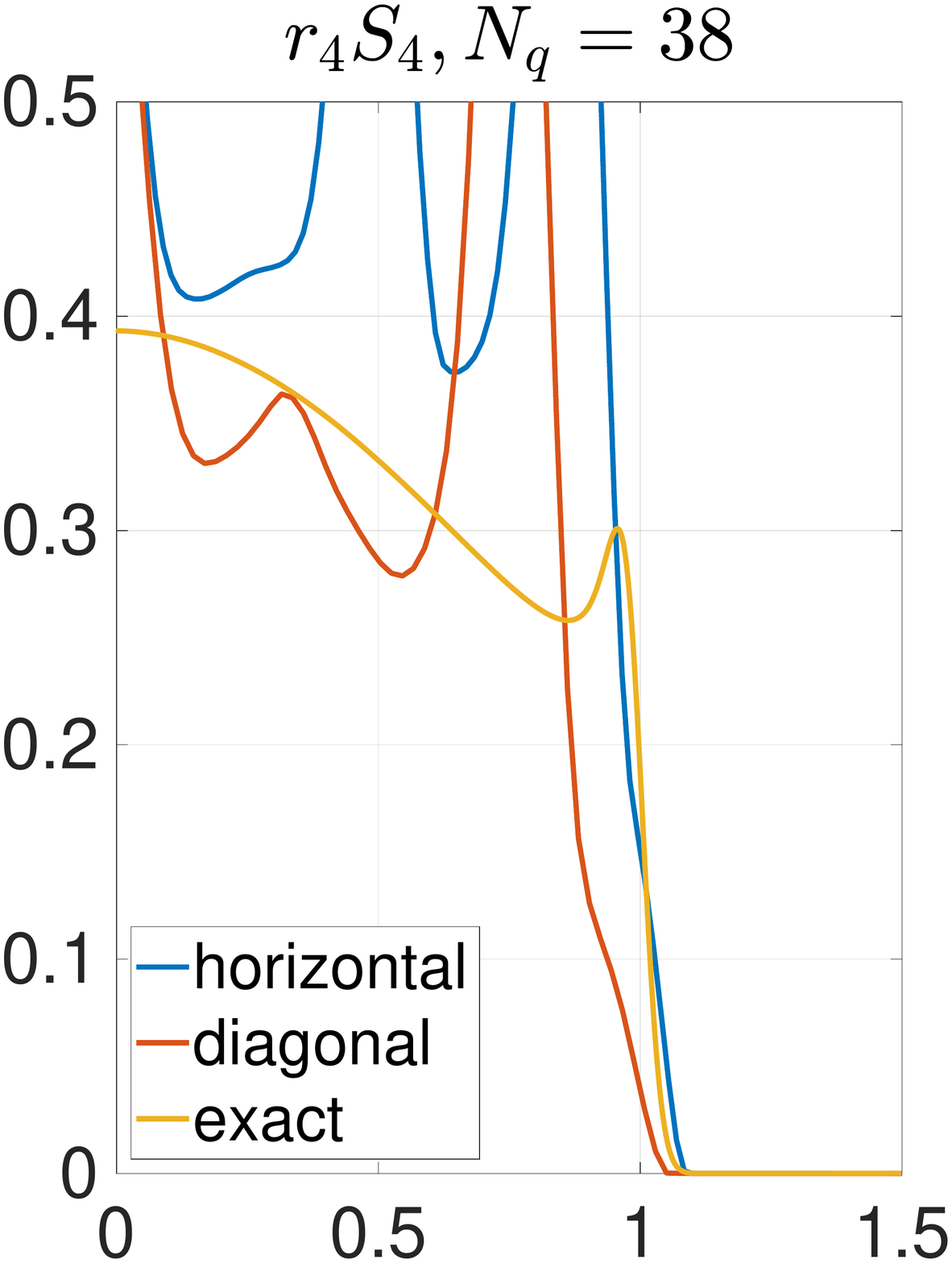}
		
		\label{fig:sub1}
	\end{subfigure}%
	\begin{subfigure}{0.24\linewidth}
		\centering
		\includegraphics[scale=0.17]{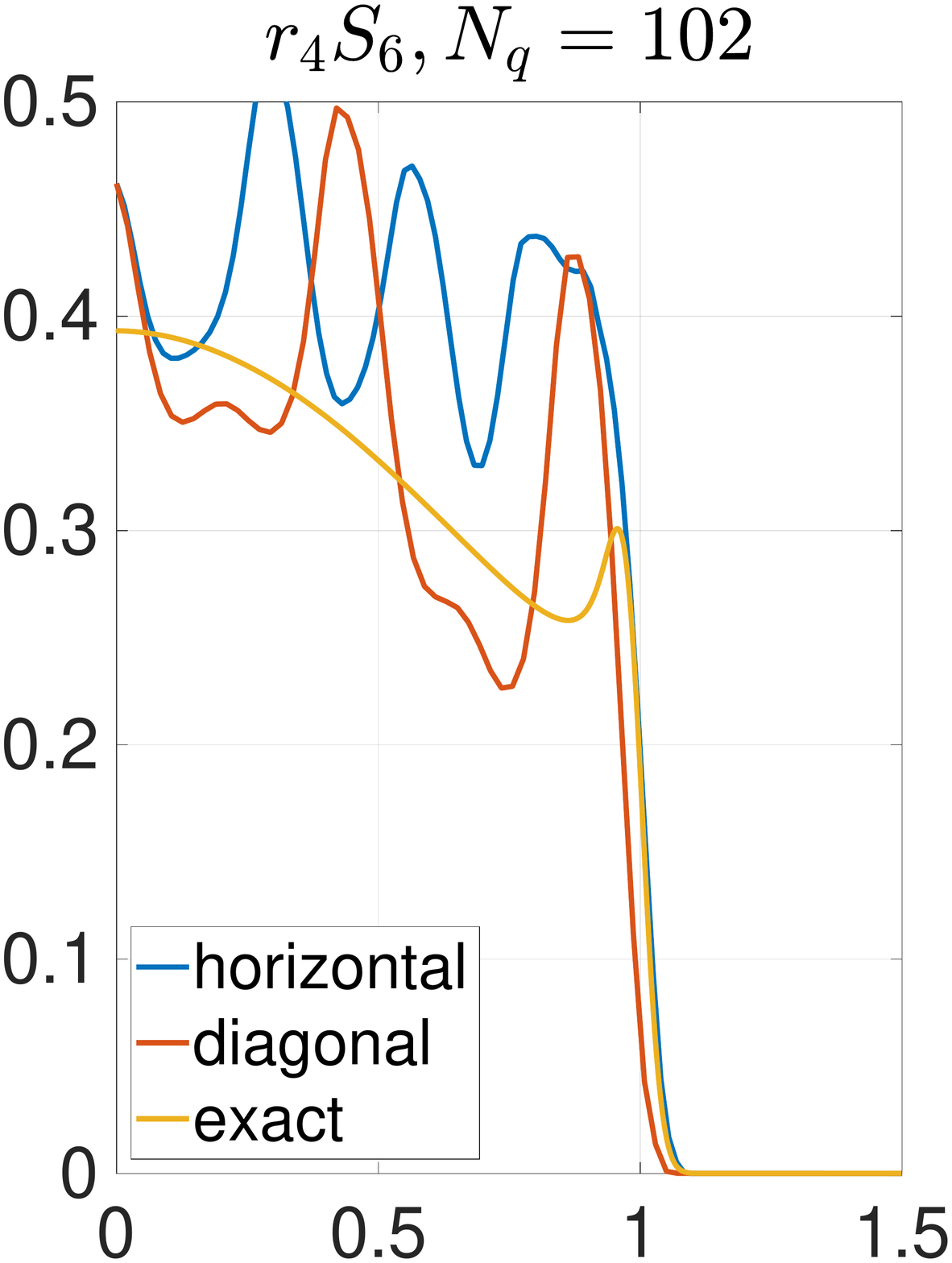}
		
		\label{fig:sub2}
	\end{subfigure}
	\begin{subfigure}{0.24\linewidth}
		\centering
		\includegraphics[scale=0.17]{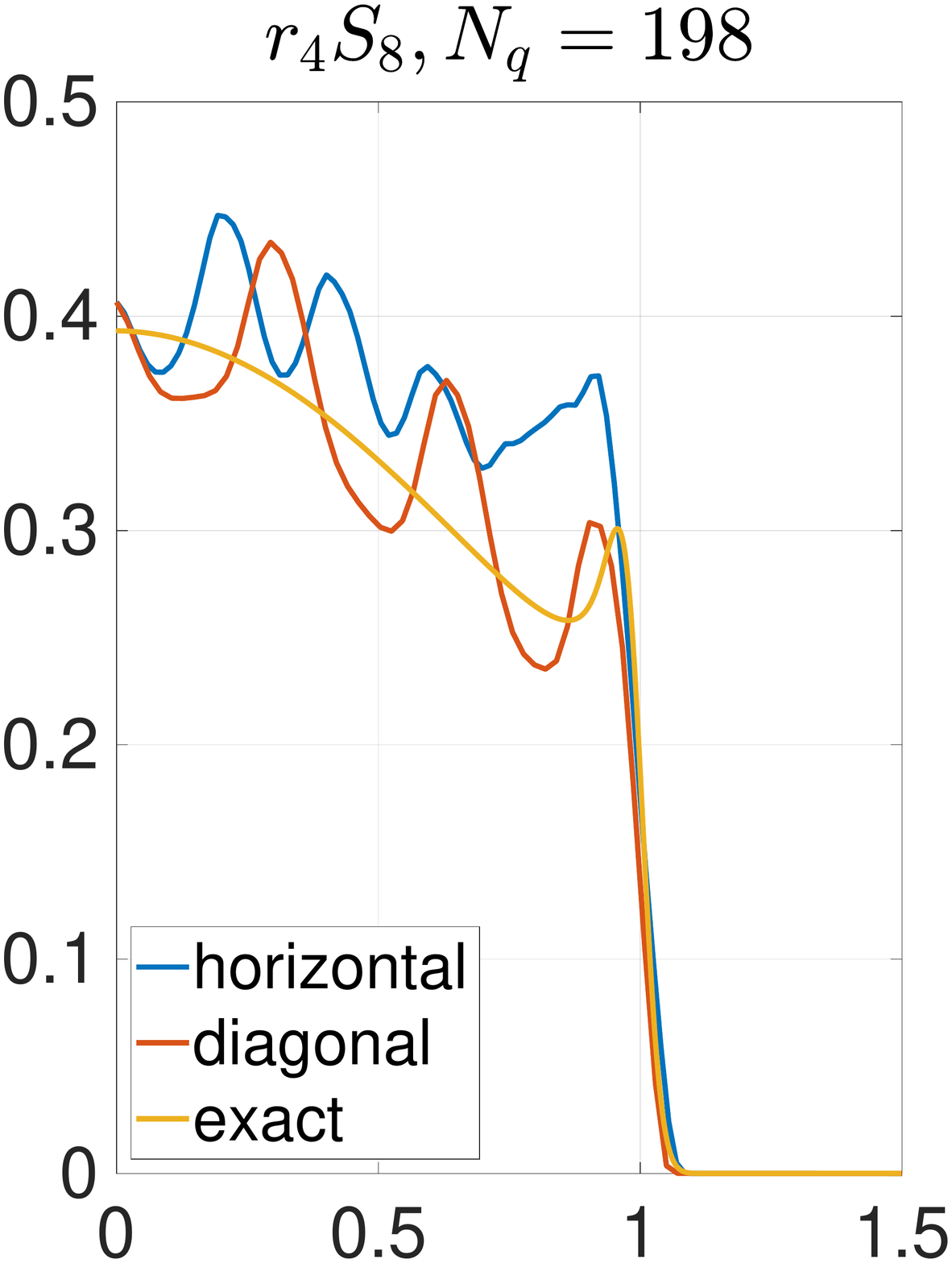}
		
		\label{fig:sub3}
	\end{subfigure}
	\begin{subfigure}{0.24\linewidth}
		\centering
		\includegraphics[scale=0.17]{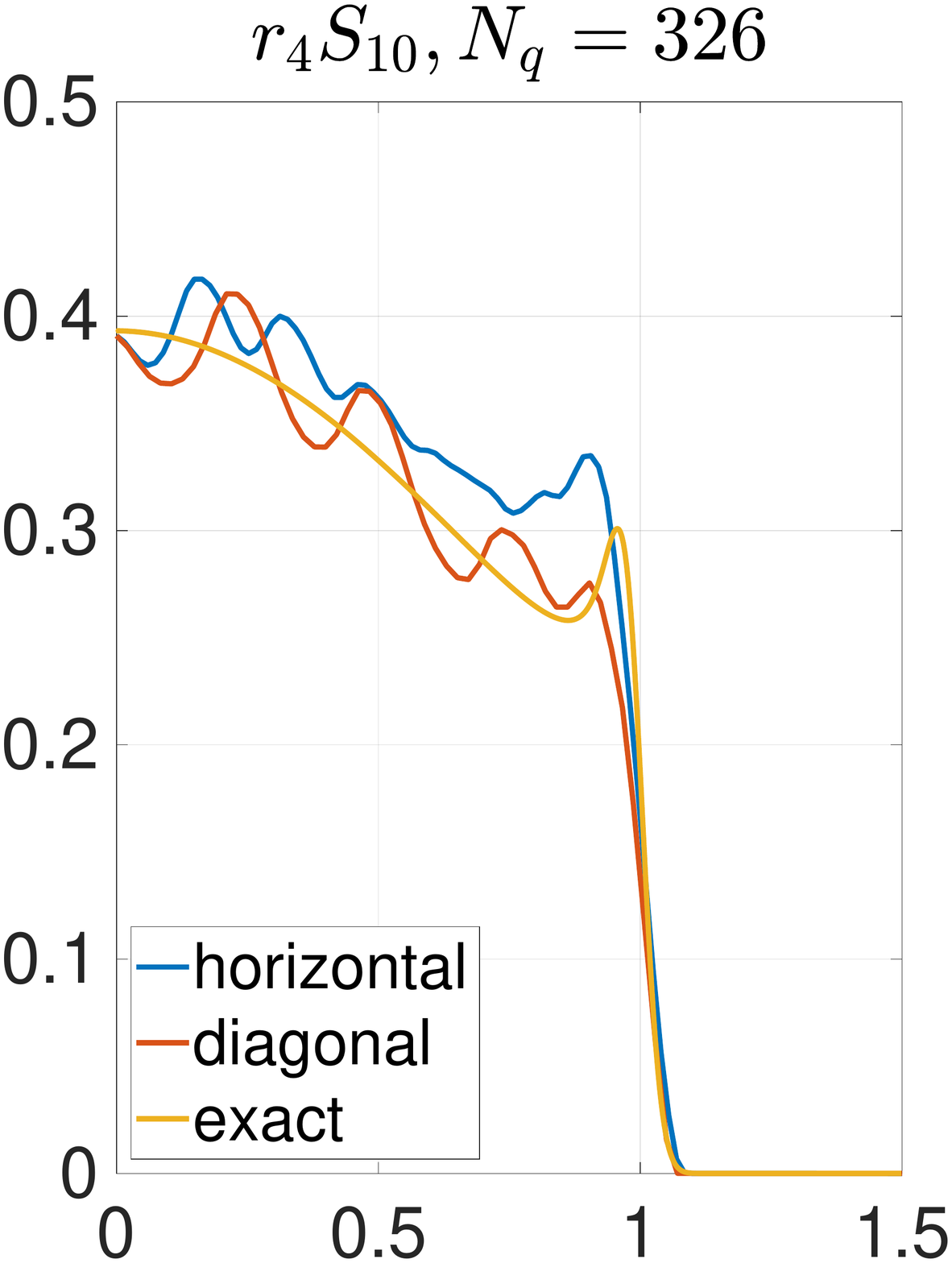}
		
		\label{fig:sub3}
	\end{subfigure}\\[-1ex]
	\begin{subfigure}{0.24\linewidth}
		\centering
		\includegraphics[scale=0.17]{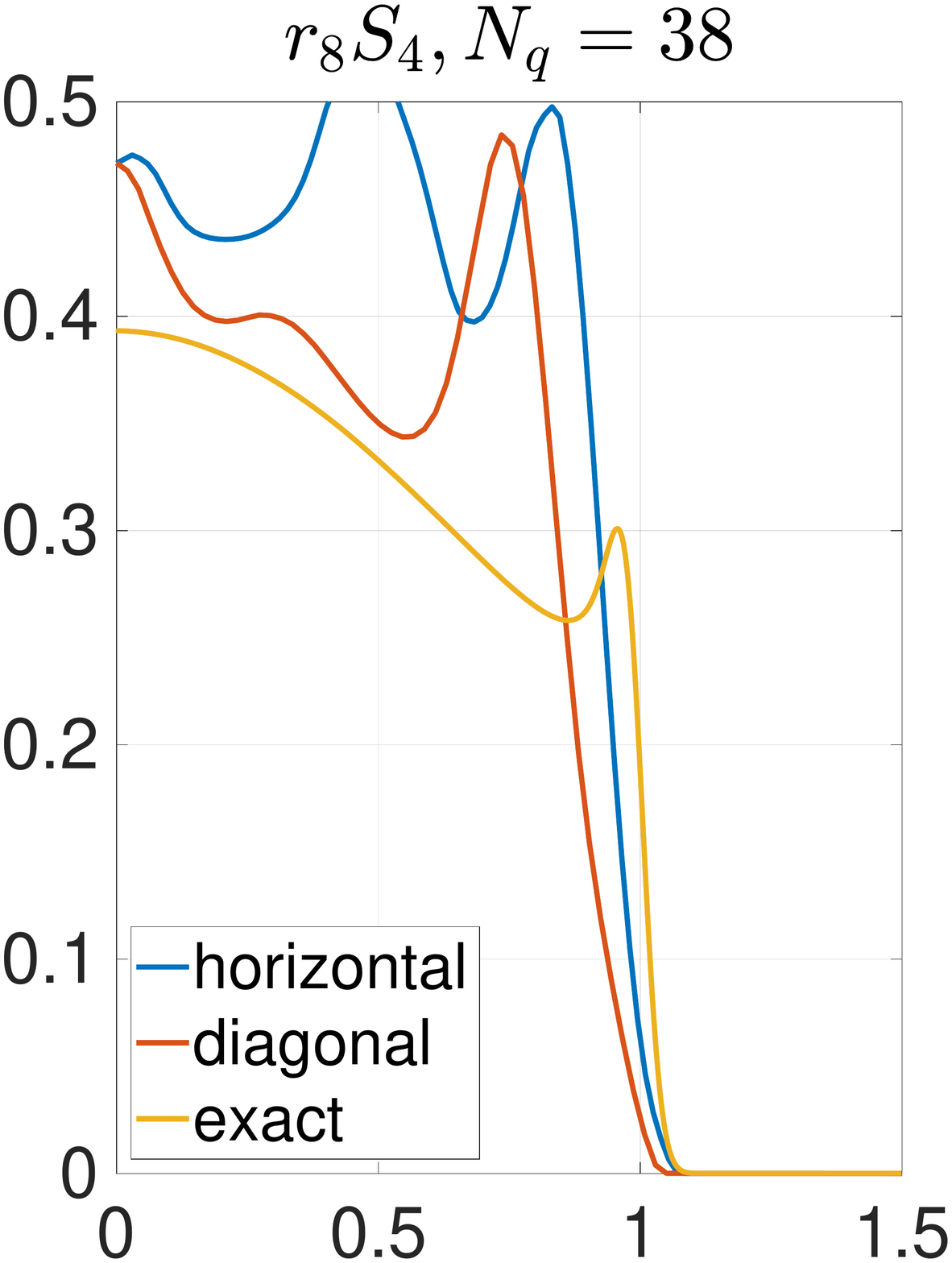}
		
		\label{fig:sub1}
	\end{subfigure}%
	\begin{subfigure}{0.24\linewidth}
		\centering
		\includegraphics[scale=0.17]{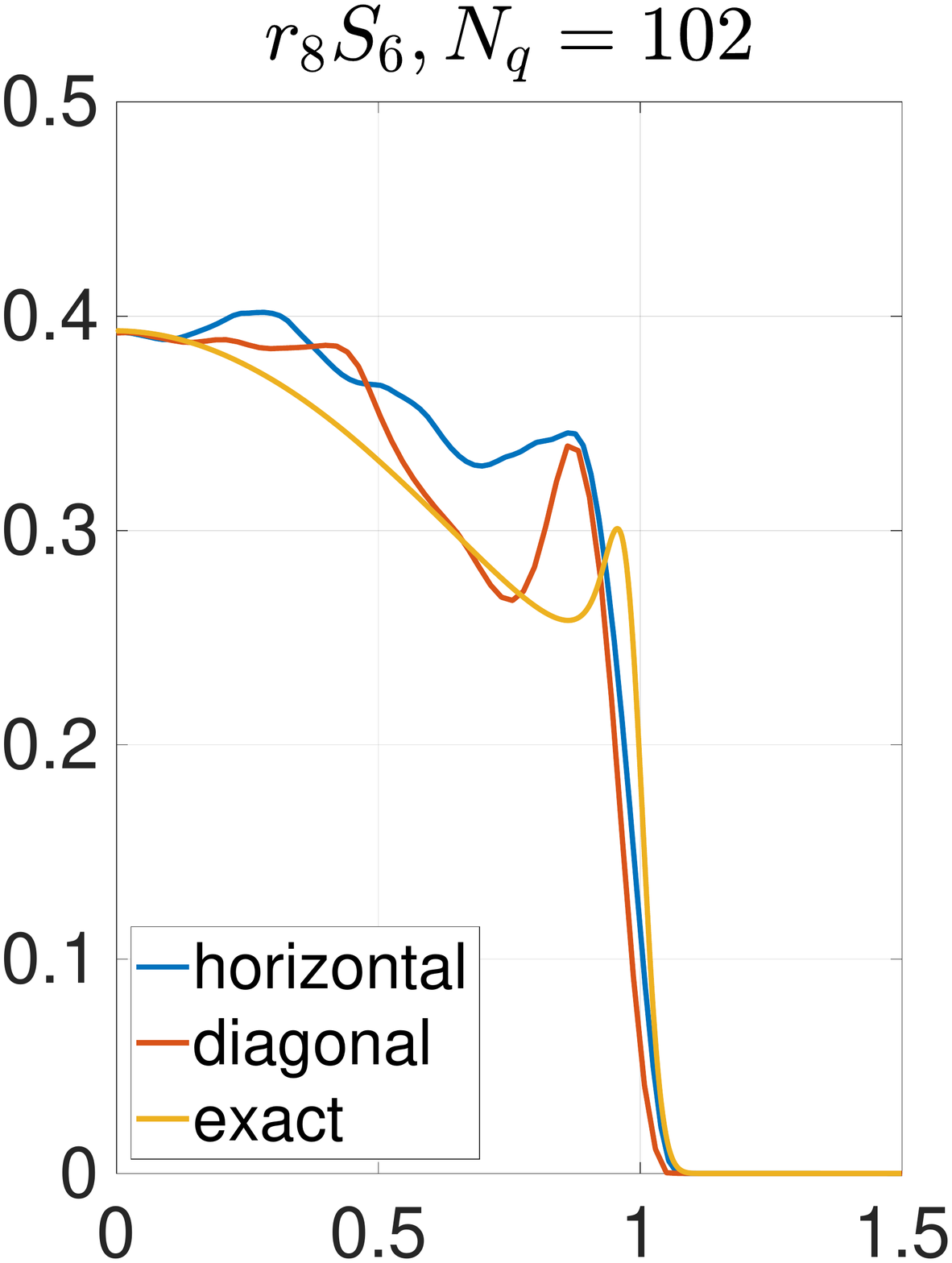}
		
		\label{fig:sub2}
	\end{subfigure}
	\begin{subfigure}{0.24\linewidth}
		\centering
		\includegraphics[scale=0.17]{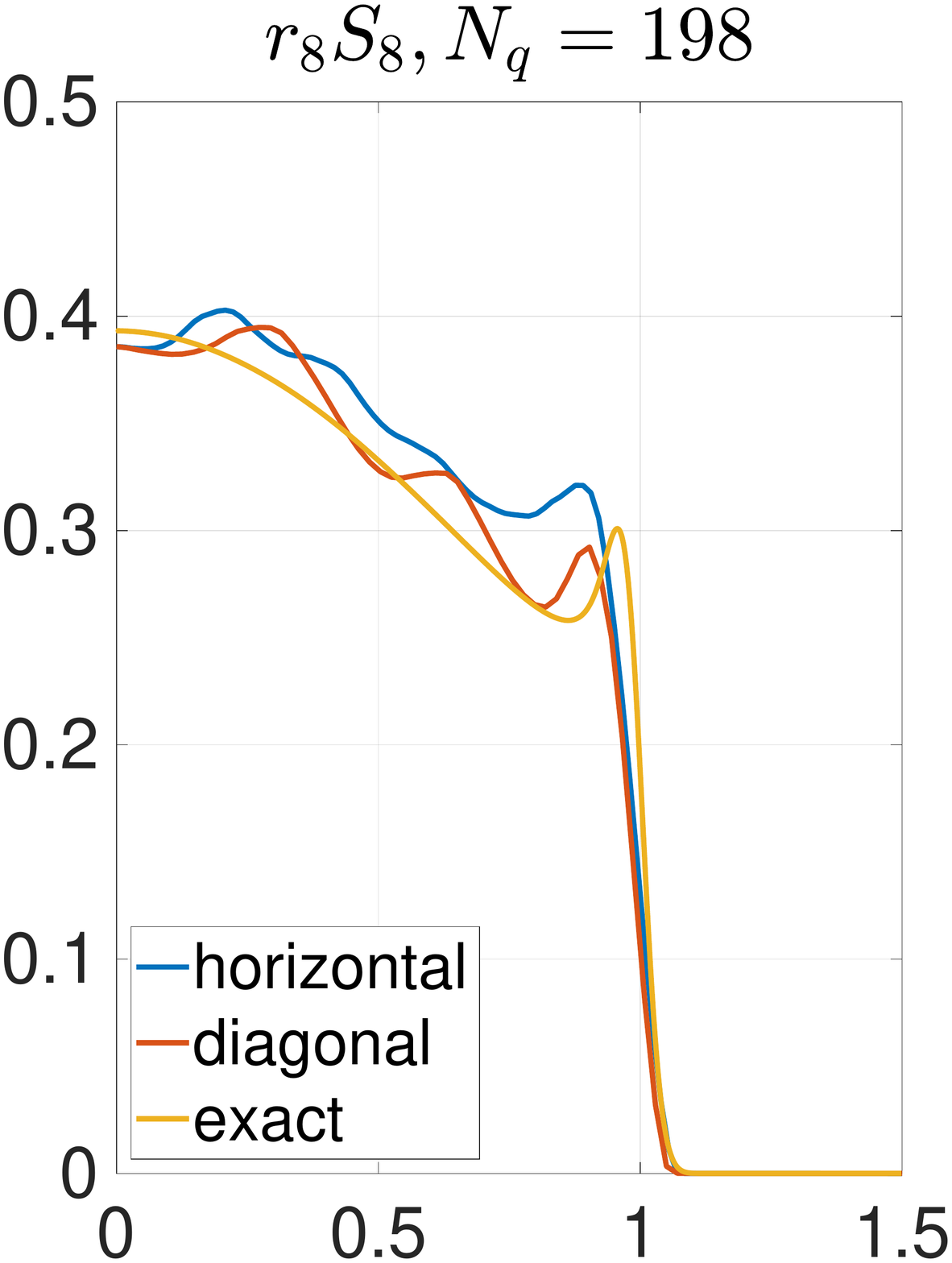}
		
		\label{fig:sub3}
	\end{subfigure}
	\begin{subfigure}{0.24\linewidth}
		\centering
		\includegraphics[scale=0.17]{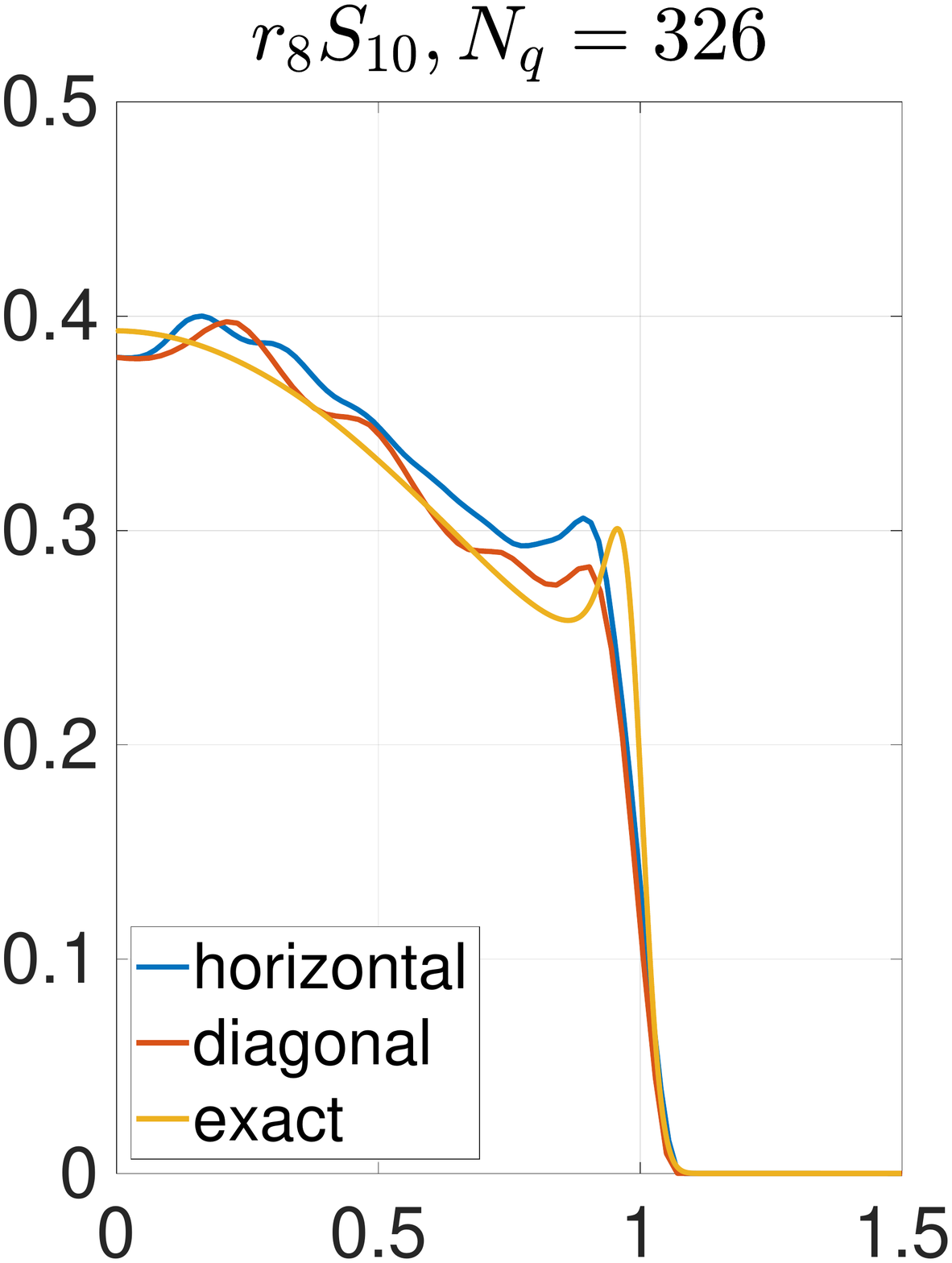}
		\label{fig:sub3}
	\end{subfigure}
	\caption{Cross sections for the line-source problem.}%
	\label{fig:LineSourcecut}
\end{figure}

The results of the densities for the line-source 
problem can be found in Fig.~\ref{fig:LineSourcecut}.
Fixed parameters used in the computation can be found in Table 
\ref{tab:linesource}.
\begin{table}
	\centering
    \begin{tabular}{ | l | p{6cm} |}
    \hline
    $a = -1.5$, $b=-1.5$ & spatial domain \\
    $N_x=200$, $N_y = 200$ & number of spatial cells \\
    $t_{end}=1.0$ & end time \\
    $\sigma_{\text{IC}}=0.03$ & parameters of initial condition \eqref{eq:IClineSource}\\
    $\sigma_a = 0$, $\sigma_s = 1.0$ & cross sections of material \\
    \hline
    \end{tabular}
    \caption{Parameters for the line-source problem.}
    \label{tab:linesource}
\end{table}

\begin{figure}
	\begin{subfigure}{0.24\linewidth}
		\centering
		\includegraphics[scale=0.17]{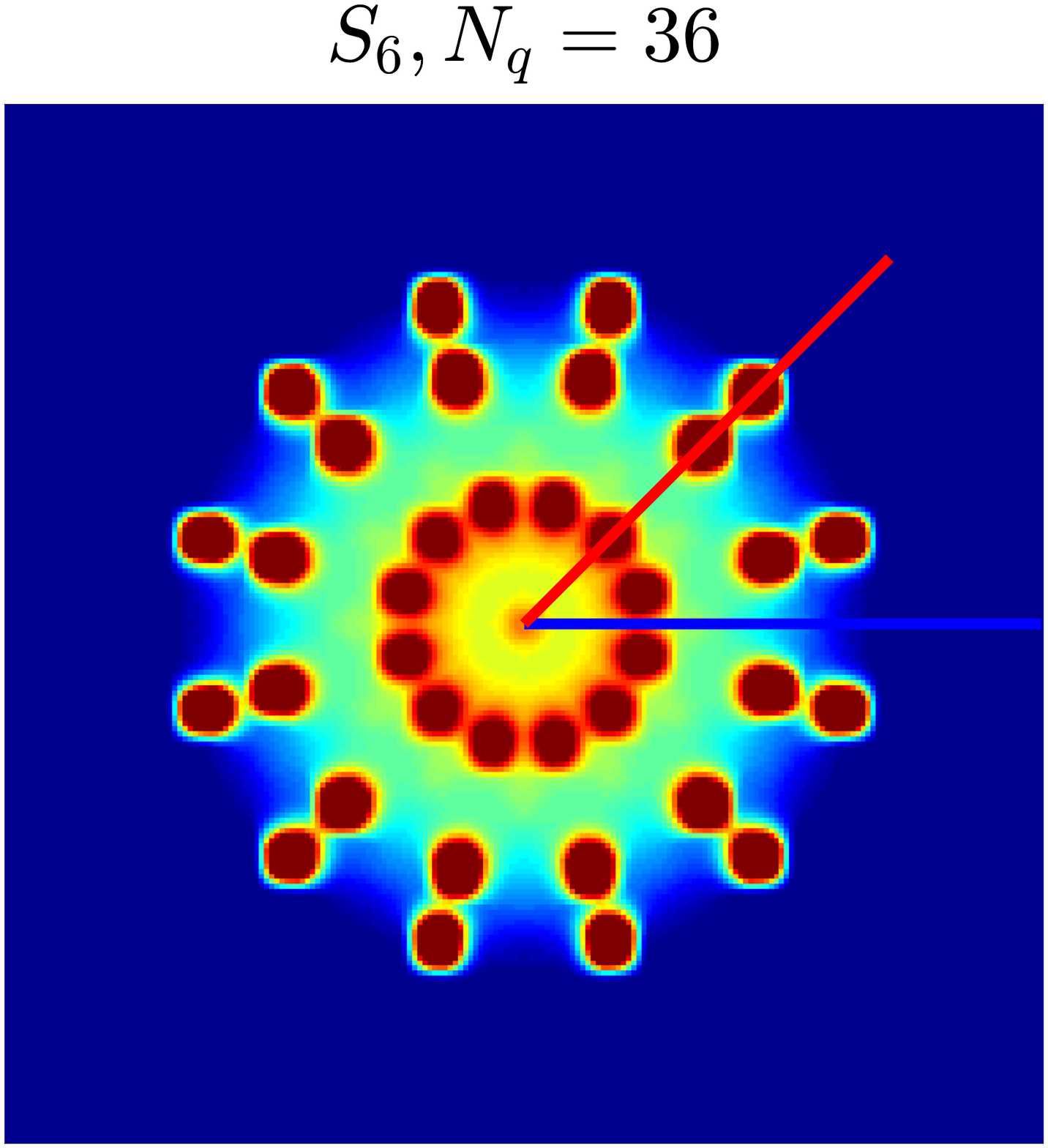}
		
		\label{fig:sub1}
	\end{subfigure}%
	\begin{subfigure}{0.24\linewidth}
		\centering
		\includegraphics[scale=0.17]{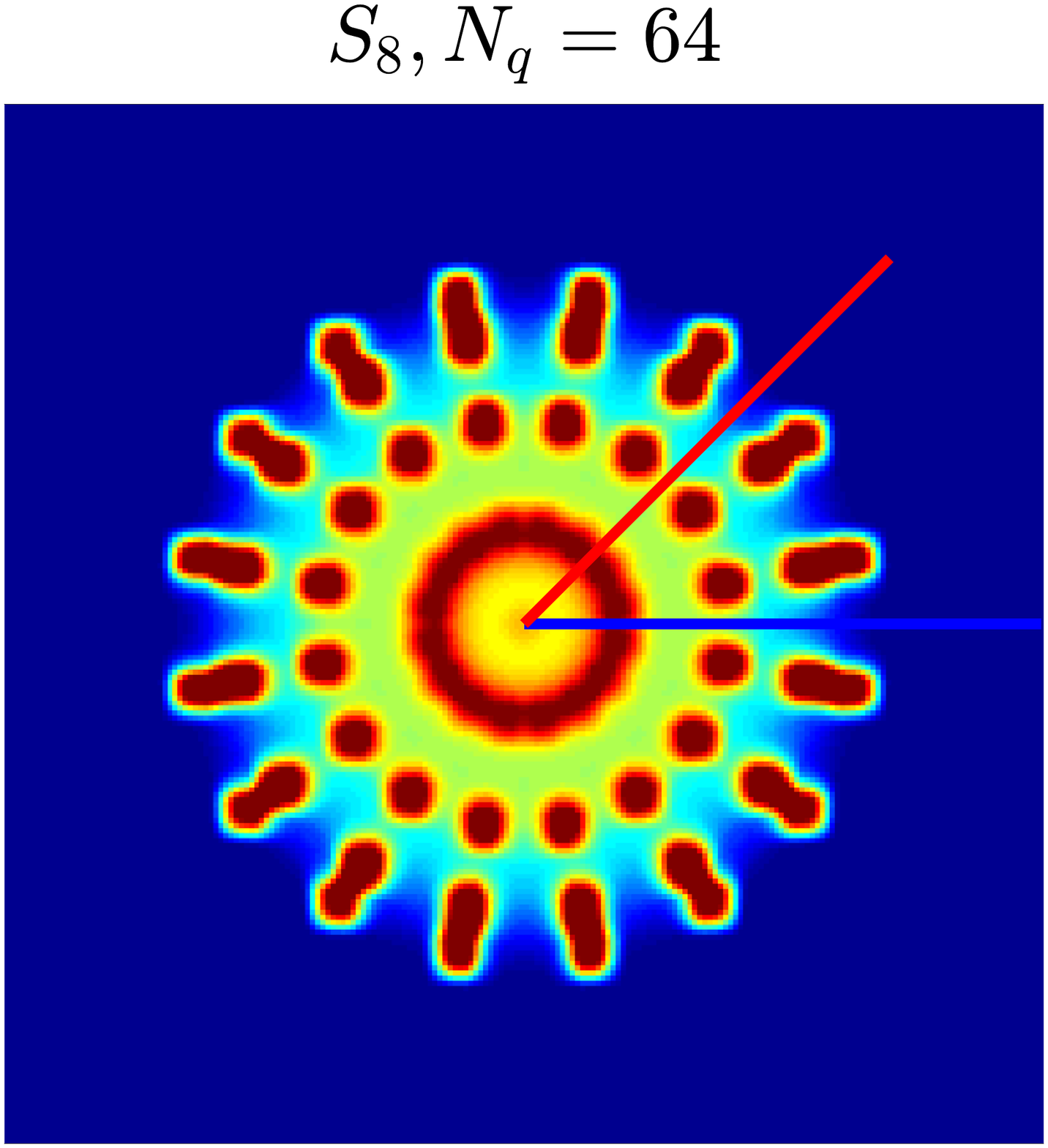}
		
		\label{fig:sub2}
	\end{subfigure}
	\begin{subfigure}{0.24\linewidth}
		\centering
		\includegraphics[scale=0.17]{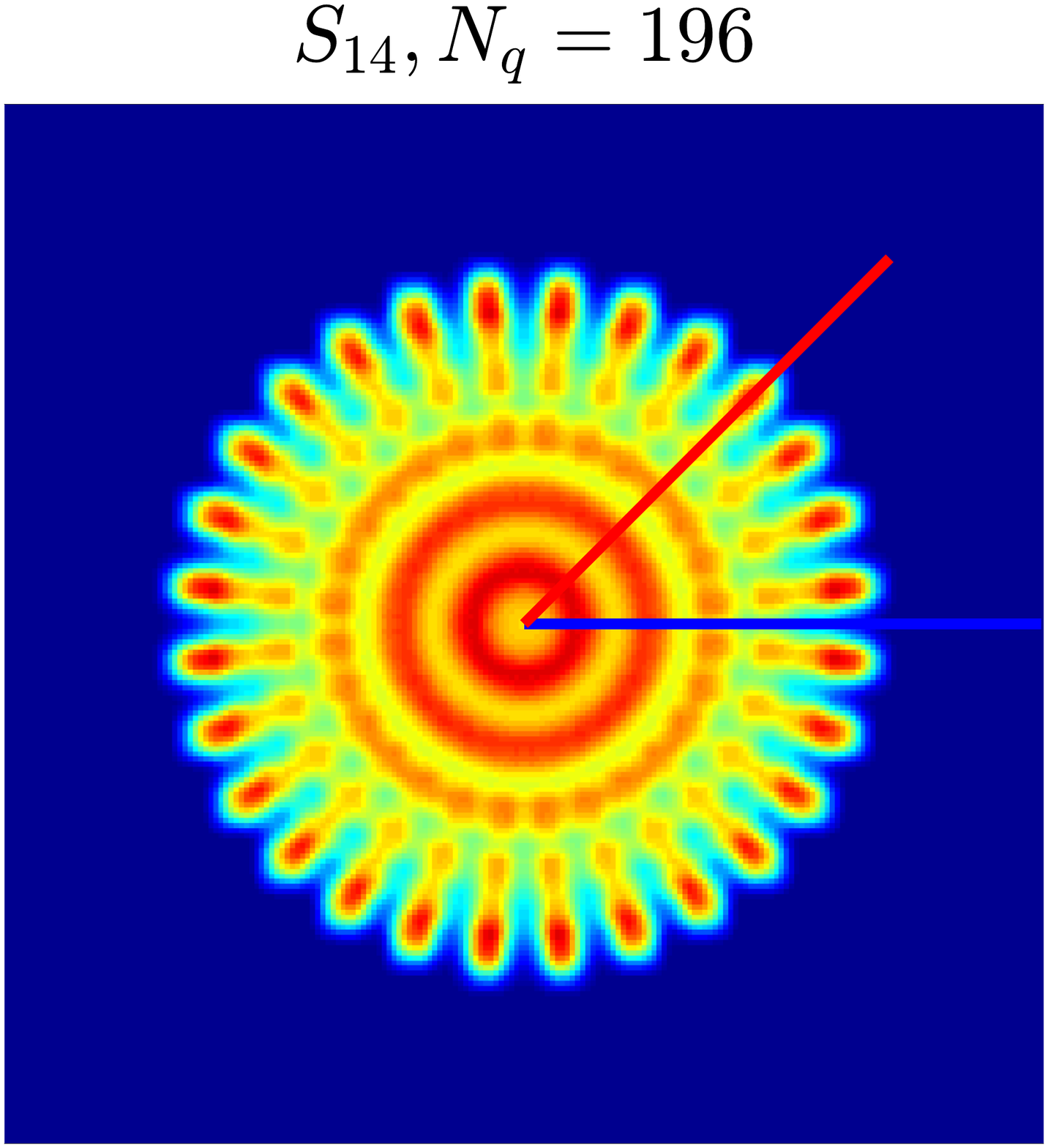}
		
		\label{fig:sub3}
	\end{subfigure}
	\begin{subfigure}{0.24\linewidth}
		\centering
		\includegraphics[scale=0.17]{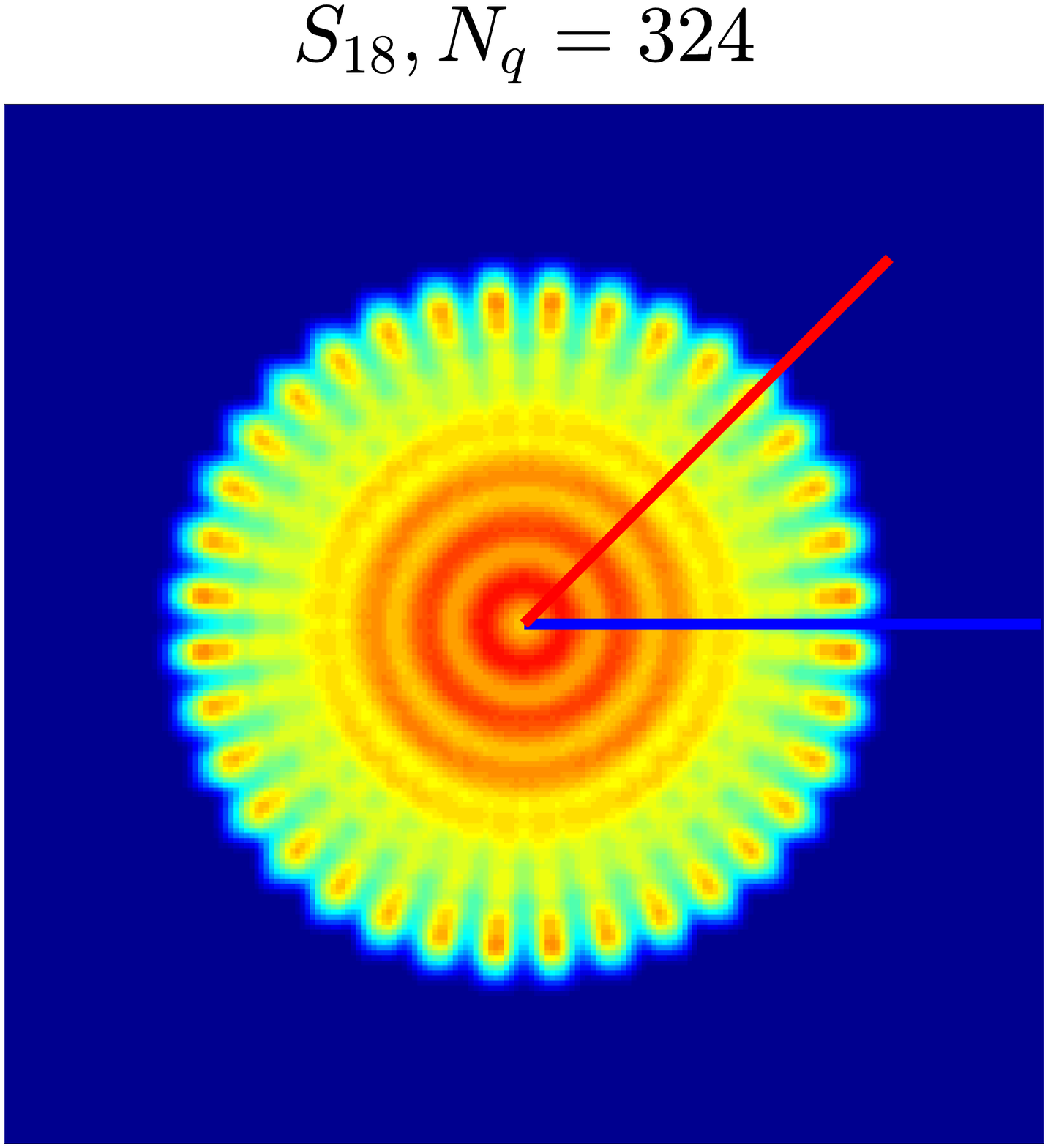}
		
		\label{fig:sub3}
	\end{subfigure}\\[-2ex]
	\begin{subfigure}{0.24\linewidth}
		\centering
		\includegraphics[scale=0.17]{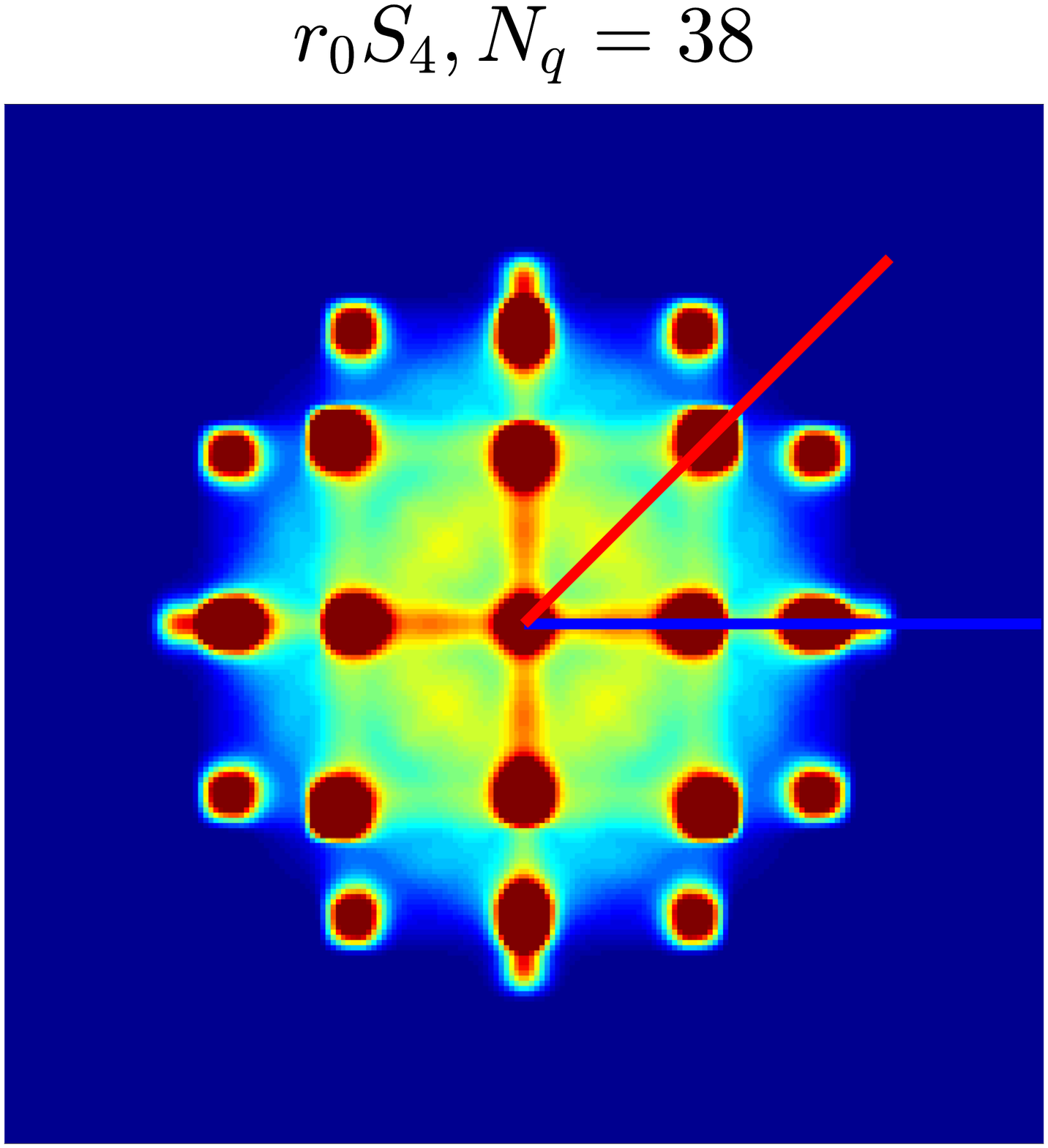}
		
		\label{fig:sub1}
	\end{subfigure}%
	\begin{subfigure}{0.24\linewidth}
		\centering
		\includegraphics[scale=0.17]{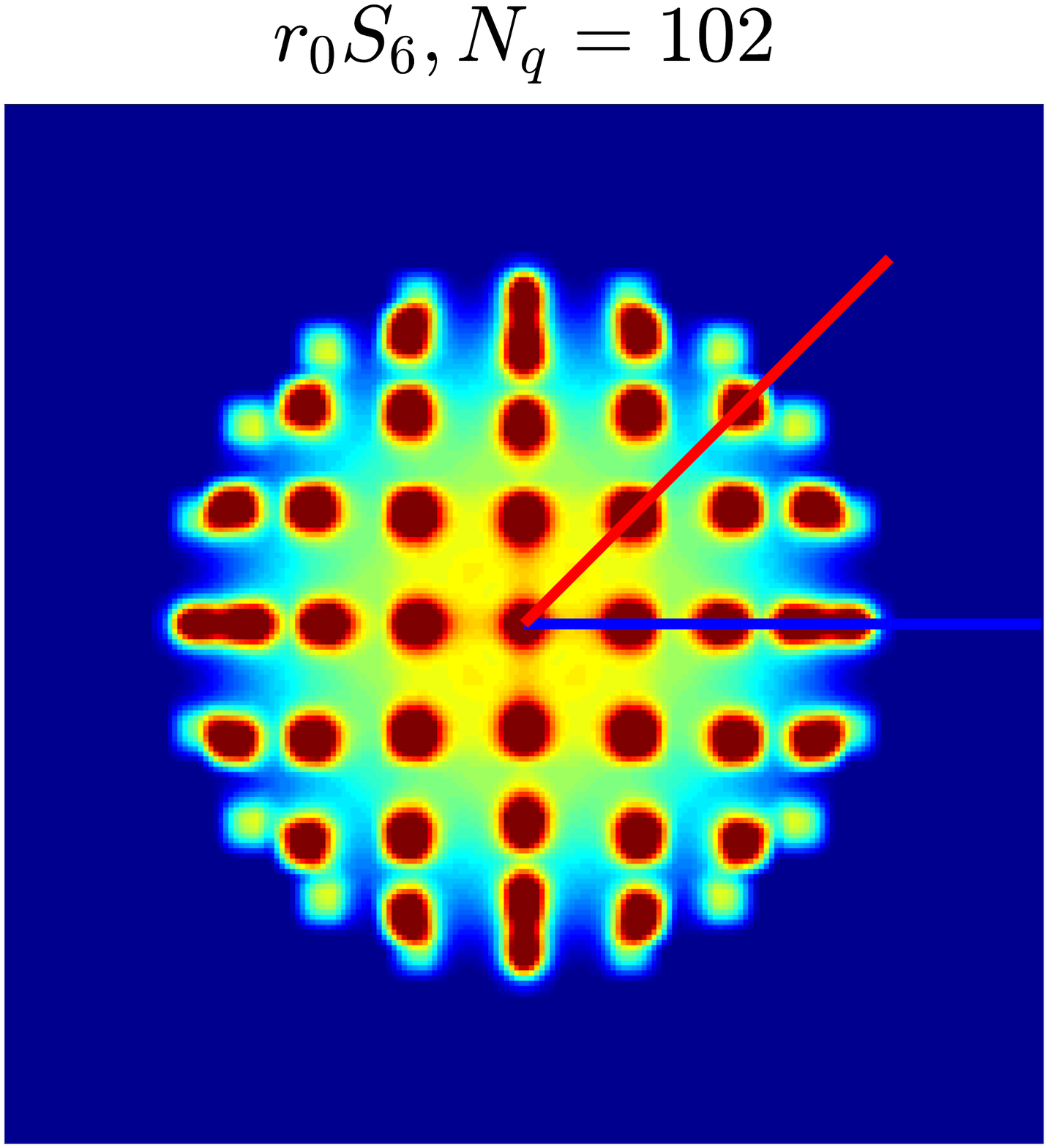}
		
		\label{fig:sub2}
	\end{subfigure}
	\begin{subfigure}{0.24\linewidth}
		\centering
		\includegraphics[scale=0.17]{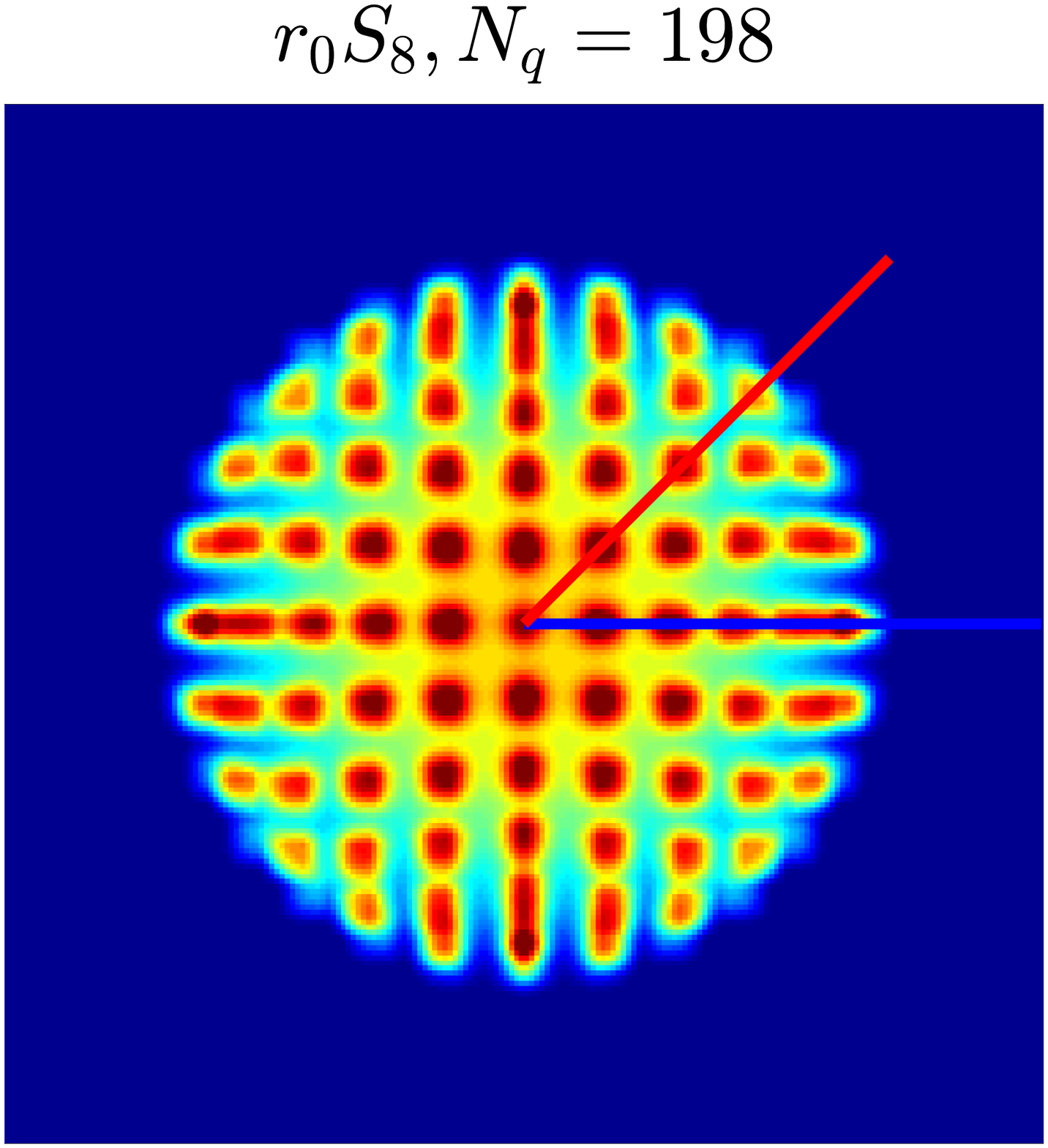}
		
		\label{fig:sub3}
	\end{subfigure}
	\begin{subfigure}{0.24\linewidth}
		\centering
		\includegraphics[scale=0.17]{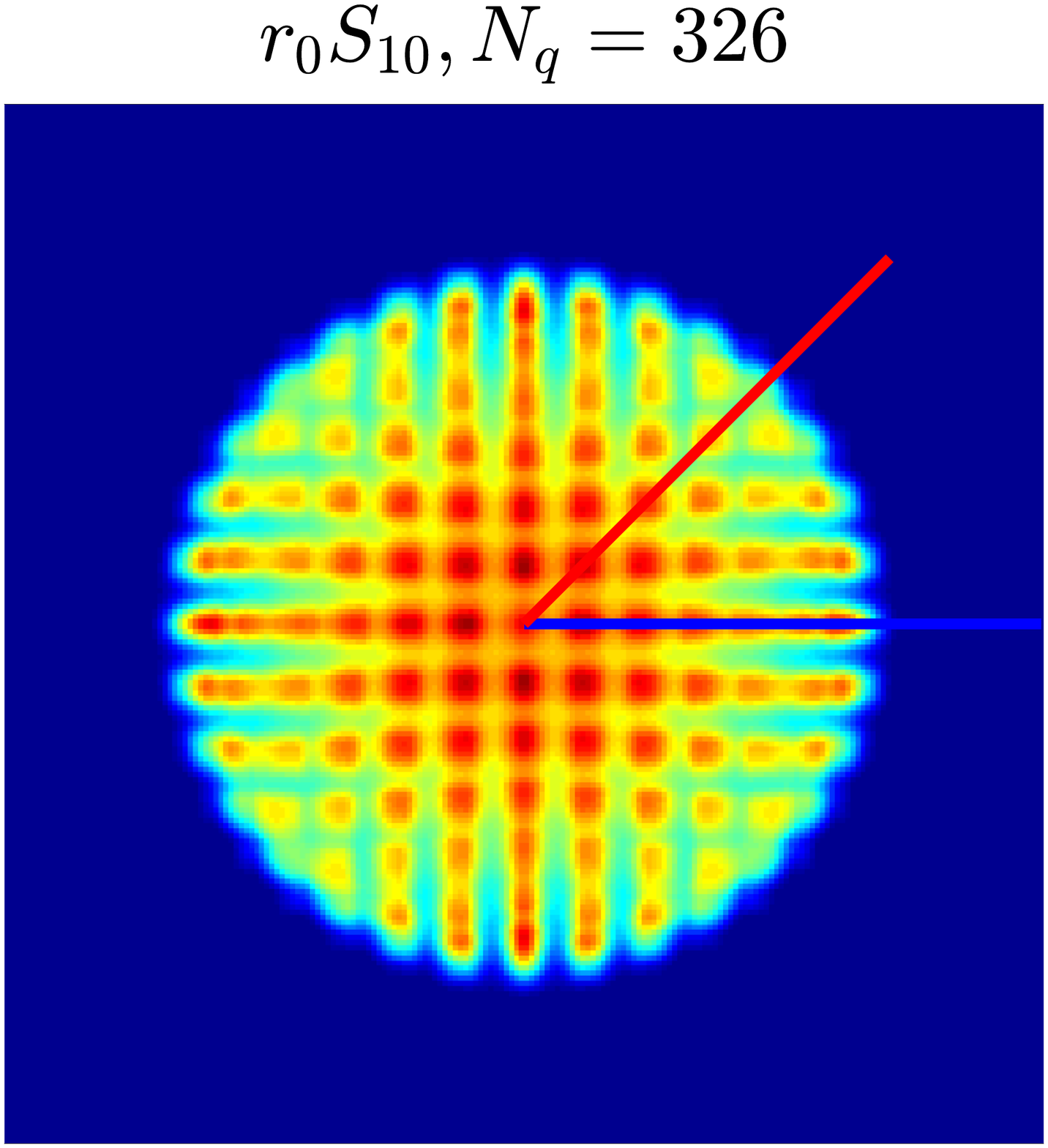}
		
		\label{fig:sub3}
	\end{subfigure}
	\\[-2ex]
	\begin{subfigure}{0.24\linewidth}
		\centering
		\includegraphics[scale=0.17]{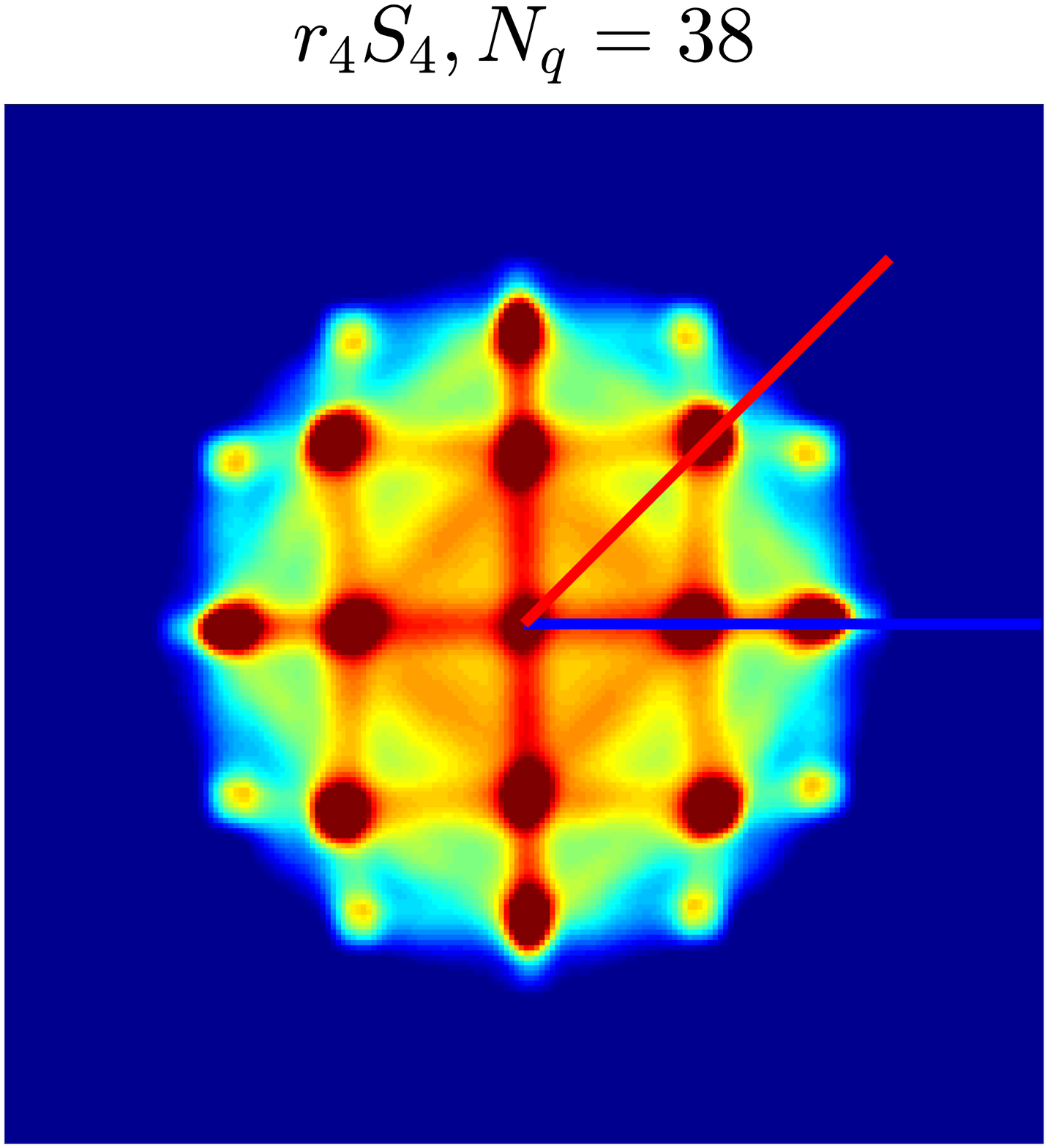}
		
		\label{fig:sub1}
	\end{subfigure}%
	\begin{subfigure}{0.24\linewidth}
		\centering
		\includegraphics[scale=0.17]{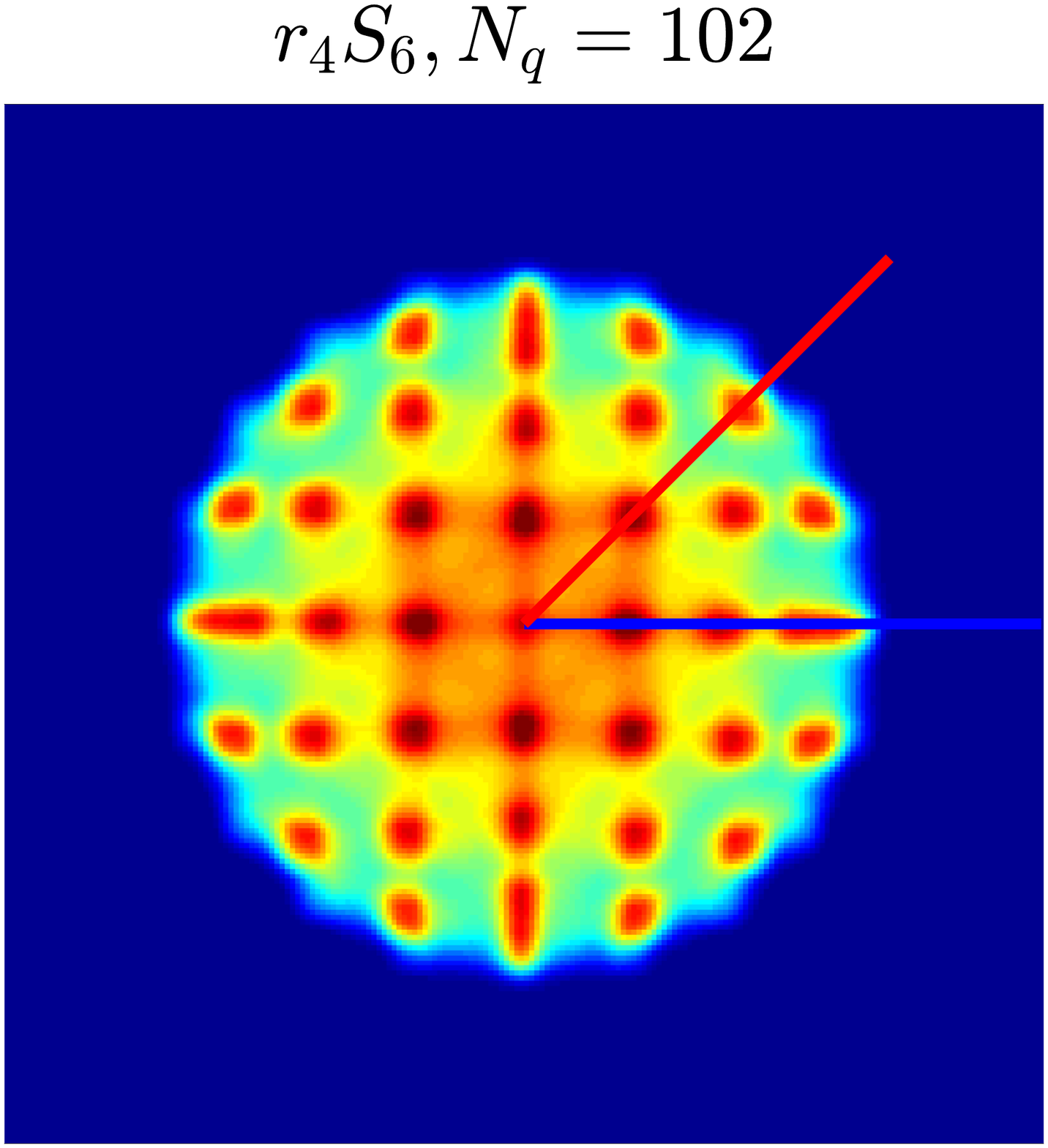}
		
		\label{fig:sub2}
	\end{subfigure}
	\begin{subfigure}{0.24\linewidth}
		\centering
		\includegraphics[scale=0.17]{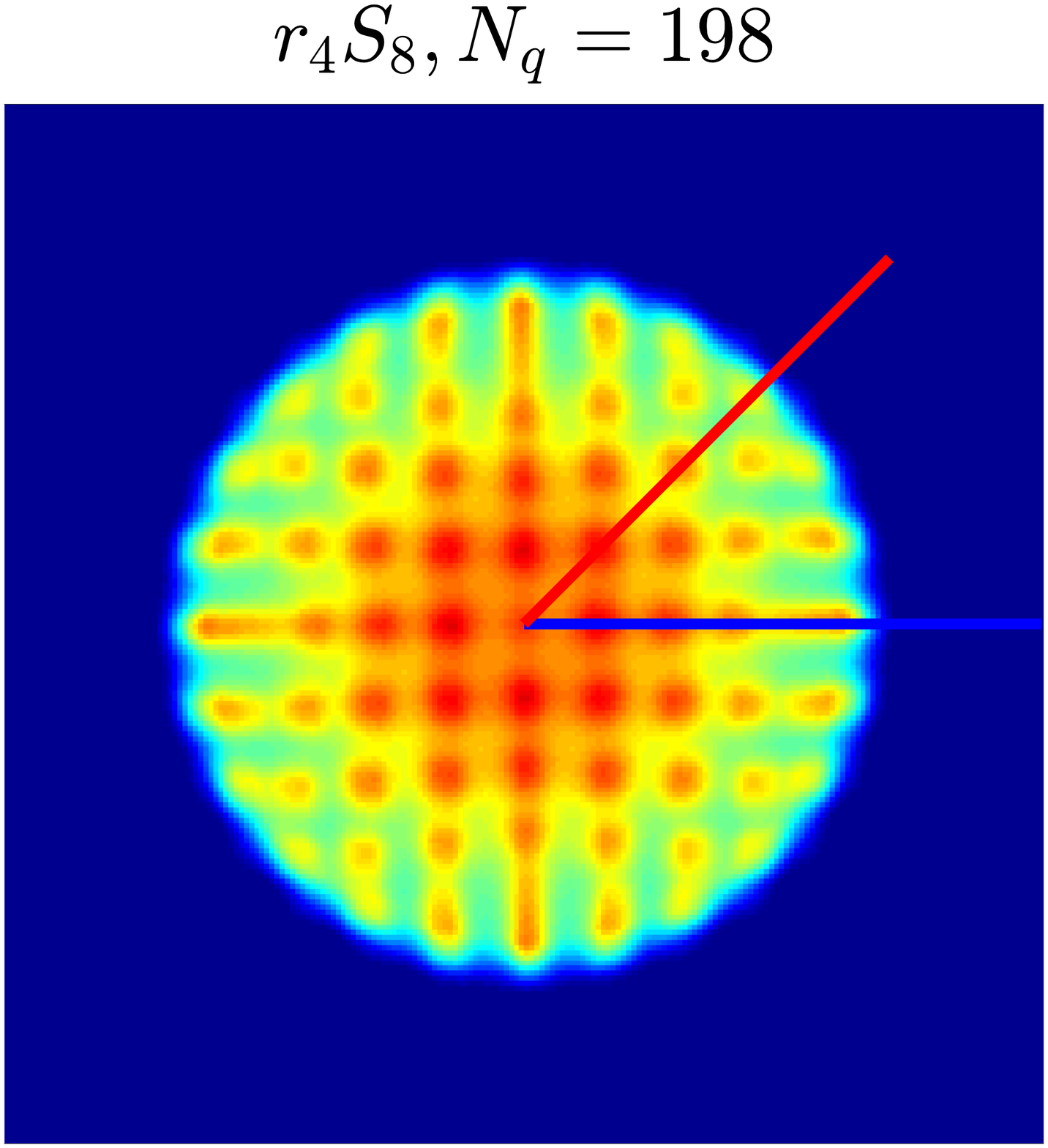}
		
		\label{fig:sub3}
	\end{subfigure}
	\begin{subfigure}{0.24\linewidth}
		\centering
		\includegraphics[scale=0.17]{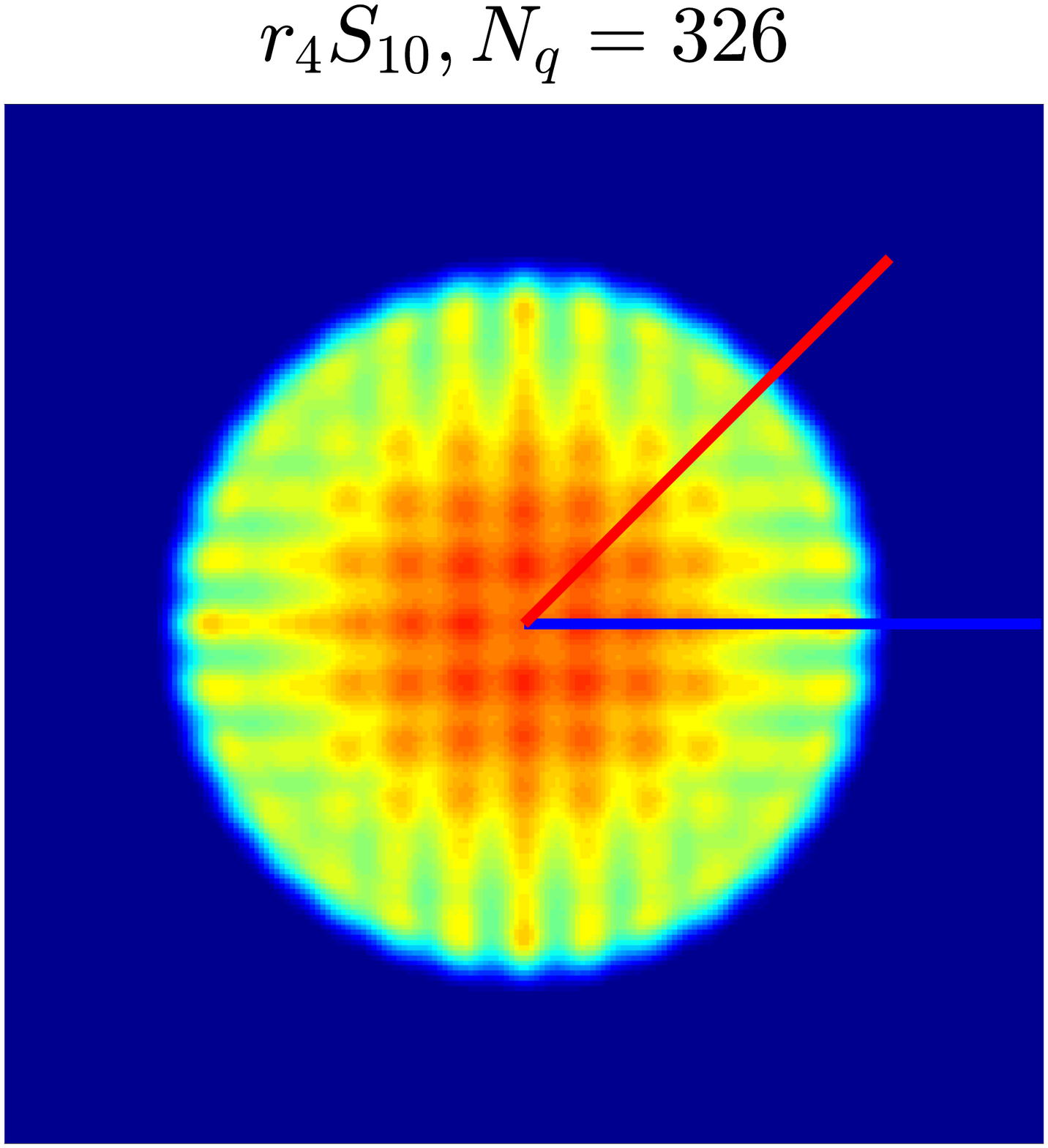}
		
		\label{fig:sub3}
	\end{subfigure}\\[-2ex]
	\begin{subfigure}{0.24\linewidth}
		\centering
		\includegraphics[scale=0.17]{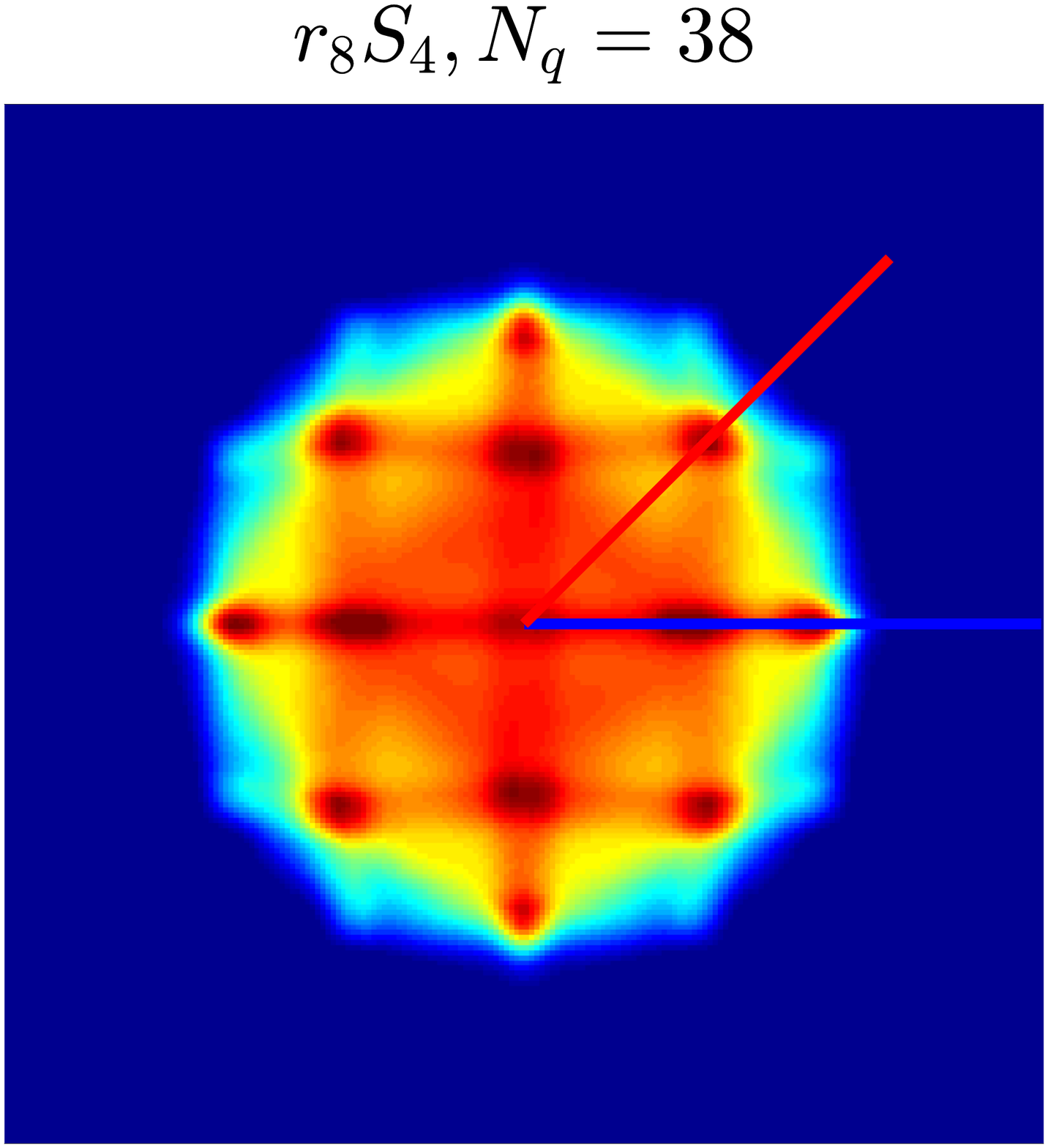}
		
		\label{fig:sub1}
	\end{subfigure}%
	\begin{subfigure}{0.24\linewidth}
		\centering
		\includegraphics[scale=0.17]{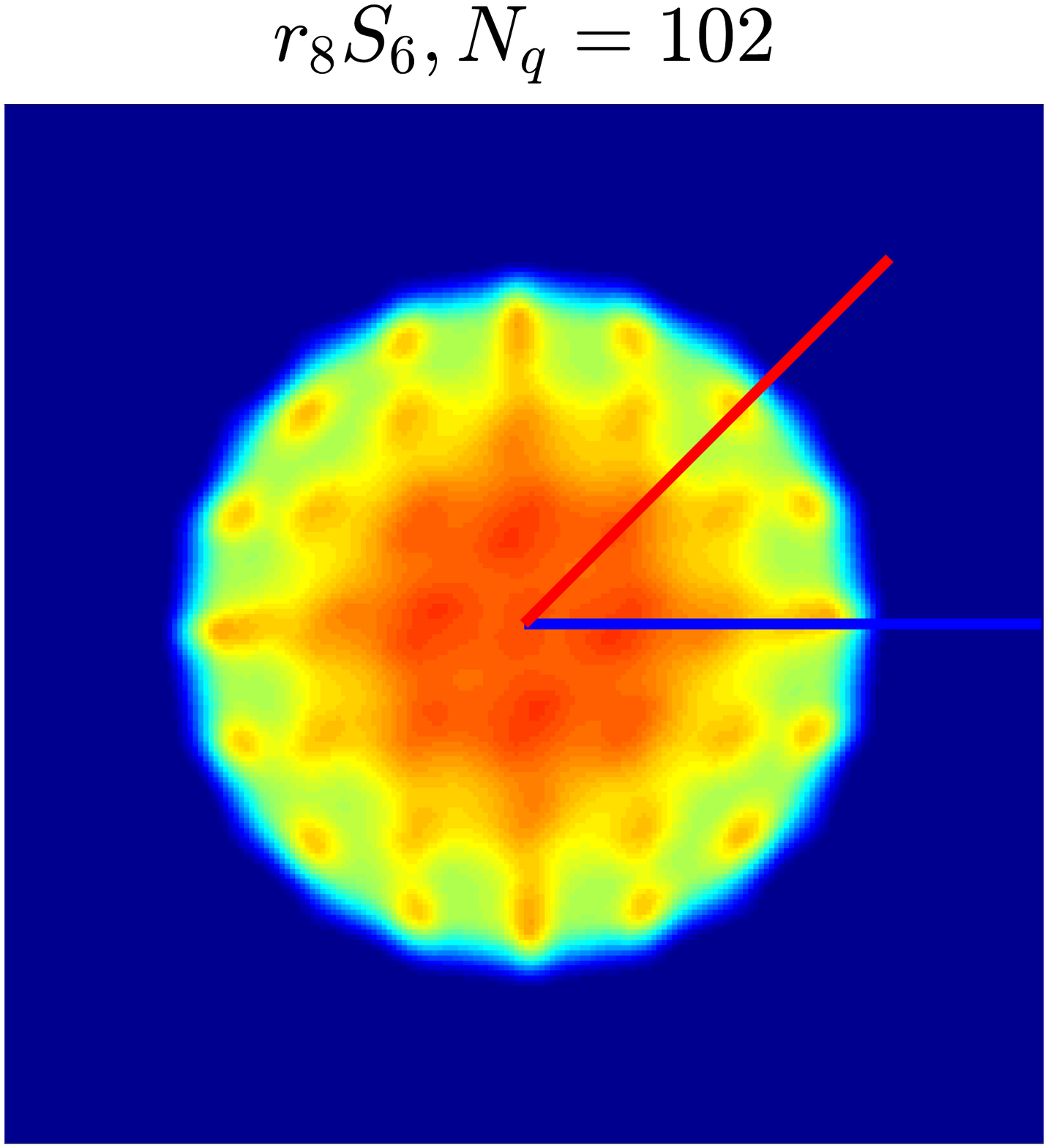}
		
		\label{fig:sub2}
	\end{subfigure}
	\begin{subfigure}{0.24\linewidth}
		\centering
		\includegraphics[scale=0.17]{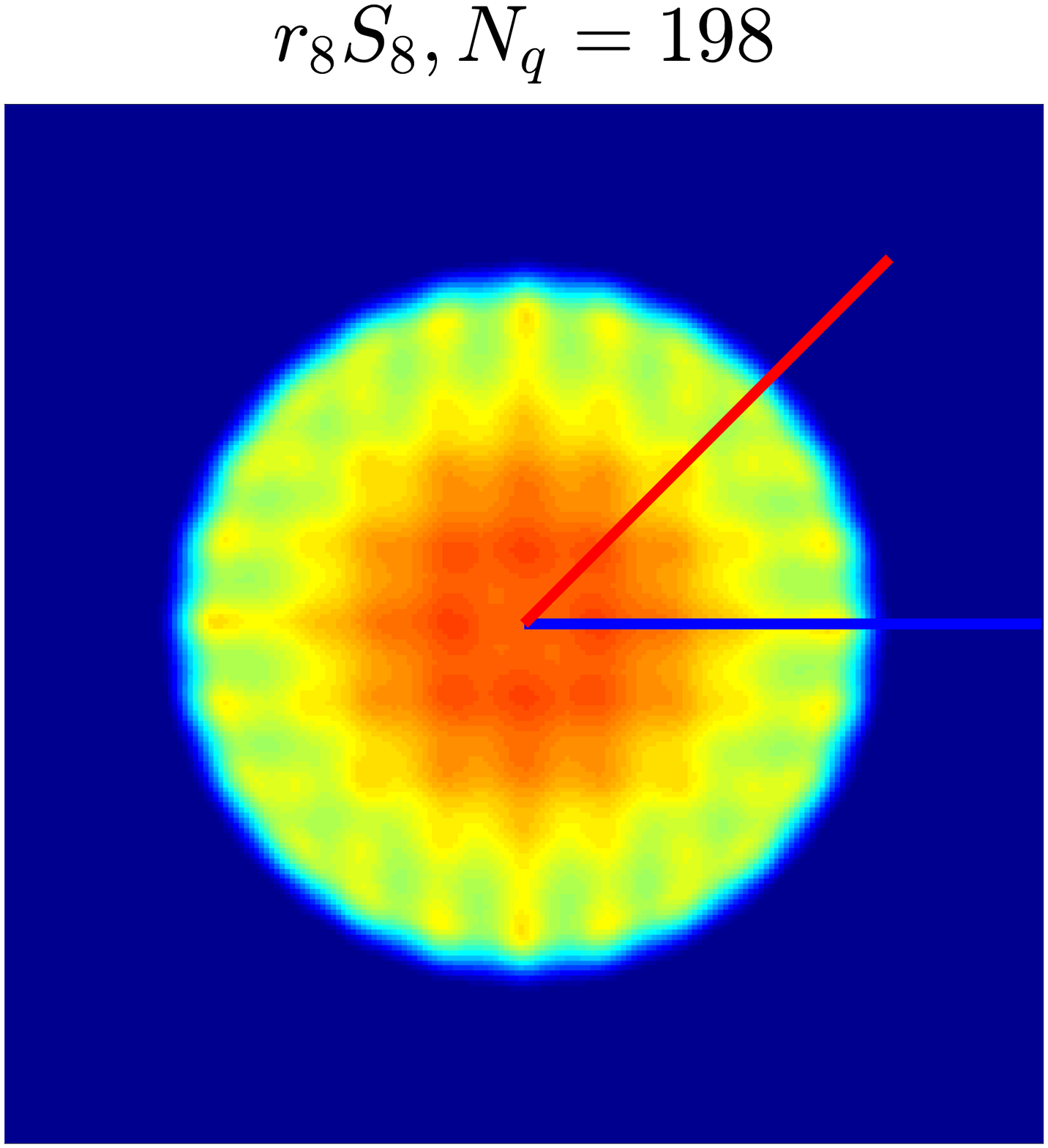}
		\label{fig:sub3}
	\end{subfigure}
	\begin{subfigure}{0.24\linewidth}
		\centering
		\includegraphics[scale=0.17]{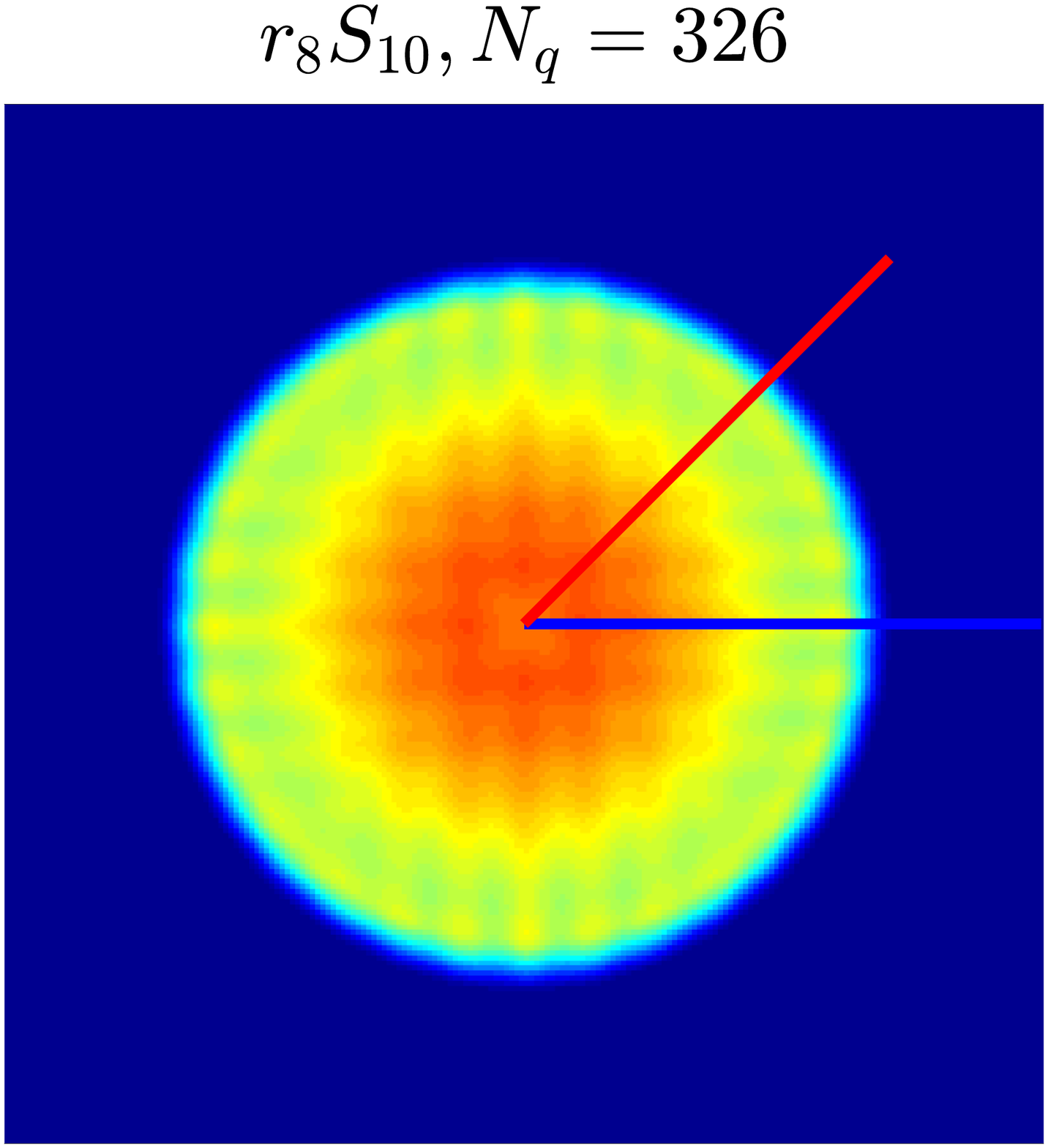}
		\label{fig:sub3}
	\end{subfigure} 
	\caption{Density for the line-source problem.}
	\label{fig:testlinesource}
\end{figure}

The first row of Fig. \ref{fig:testlinesource} shows the density for the 
line-source problem, computed with the standard \SN method for different 
numbers of quadrature points. Ray effects are predominant especially for the 
computations with 36 and 64 quadrature points. 
The solution along the blue horizontal and red diagonal line 
is shown in Fig.~\ref{fig:LineSourcecut} together with the reference solution. 
While the density along the horizontal cut is mostly underestimated, the 
solution along the vertical cut is mostly overestimated. Even for 324 
quadrature points, the numerical solution has not yet converged against the 
reference solution.

The second to fourth row show the results for the \rSN method with different 
rotation magnitudes $\delta$ and number of quadrature points, i.e. 
$\displaystyle 
\delta\in \{0,4,8\}$ and $N_q\in \{38,102,198,326\}$. Since the number of 
quadrature points for the \SN method scales different than for the \rSN, it is 
not possible to match the number of quadrature points exactly. For the case of 
$\delta=0$, i.e. no rotation, two quadrature points fall onto another when 
being projected into the $xy$-plane. This effectively halves the number of 
quadrature points.
We still observe ray effects when using the new quadrature set. However, when 
rotating the quadrature set in every time step by $\delta=4$ or $\delta=8$, the 
presence of these ray effects reduces dramatically. This can be seen in 
Fig.~\ref{fig:testlinesource}, but more precisely in 
Fig~\ref{fig:LineSourcecut}. The solution along the horizontal and diagonal cut 
are both closer to another, as well as closer to the reference solution. 
Oscillations in radial direction can be reduced significantly, which is 
observable when comparing the first and last row of 
Fig~\ref{fig:LineSourcecut}, respectively.
When comparing the \rSN method for 
different values of $\delta$ but the same number of quadrature points, we 
observe a slight decrease in the propagation speed of the solution. The 
wavefront of the solution slightly moves back for higher values of $\delta$.
However, comparing this with the wavefront of the standard \SN method, the \rSN 
method still manages to capture the actual wavefront more accurately in all 
configurations. 

The superiority of the \rSN method can be observed when comparing \rSN with 
$\delta=8$ and $N_q=102$ against the standard \SN method with $N_q=324$. With 
less than a third of the number of quadrature points, the \rSN method has fewer 
oscillations and varies less when comparing the horizontal cut with the 
diagonal cut. 
Since the additional computational cost of the rotation and interpolation 
procedure never exceeds 10\%, the \rSN method can yield similar, if not better, 
results for the line-source problem with a third of the costs due to fewer 
quadrature points. Maybe more importantly, the memory footprint is reduced as well.

\subsection{The lattice problem}
\begin{table}[]%
	\centering
	\centering
	\begin{tabular}{ | l | p{6cm} |}
		\hline
		$a = 0$, $b=7$ & spatial domain \\
		$N_x=280$, $N_y = 280$ & number of spatial cells \\
		$t_{end}=3.2$ & end time \\
		$\sigma_a = 0$, $\sigma_s = 1.0$ & cross sections of 
		background \\
		$\displaystyle \sigma_a = 10.0$, $\sigma_s = 0$ & cross sections of 
		squares \\
		$Q = 1$ & strength of isotropic source \\
		\hline
	\end{tabular}
	\caption{Parameters for the lattice problem.}
	\label{fig:cbSetting}
\end{table}

%\qquad
\begin{figure}[]
	\centering
	\includegraphics[width=0.15\textwidth]{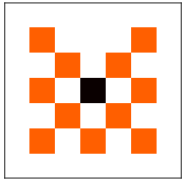}
	\caption{Layout of the lattice problem.}
	\label{fig:cblayout}
\end{figure}

\begin{figure}
	\begin{subfigure}{0.24\linewidth}
		\centering
		\includegraphics[scale=0.17]{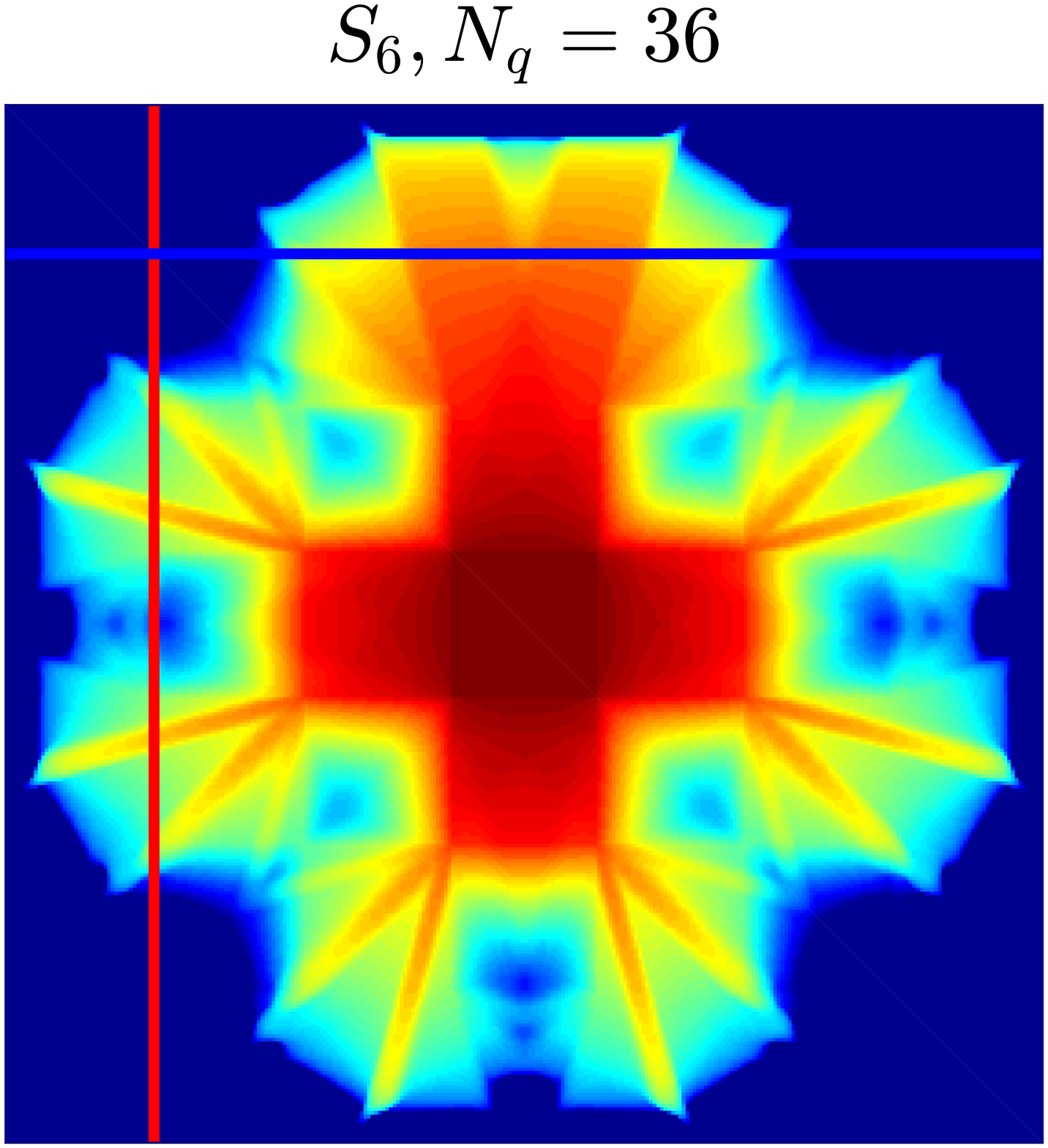}
		
		\label{fig:sub1}
	\end{subfigure}%
	\begin{subfigure}{0.24\linewidth}
		\centering
		\includegraphics[scale=0.17]{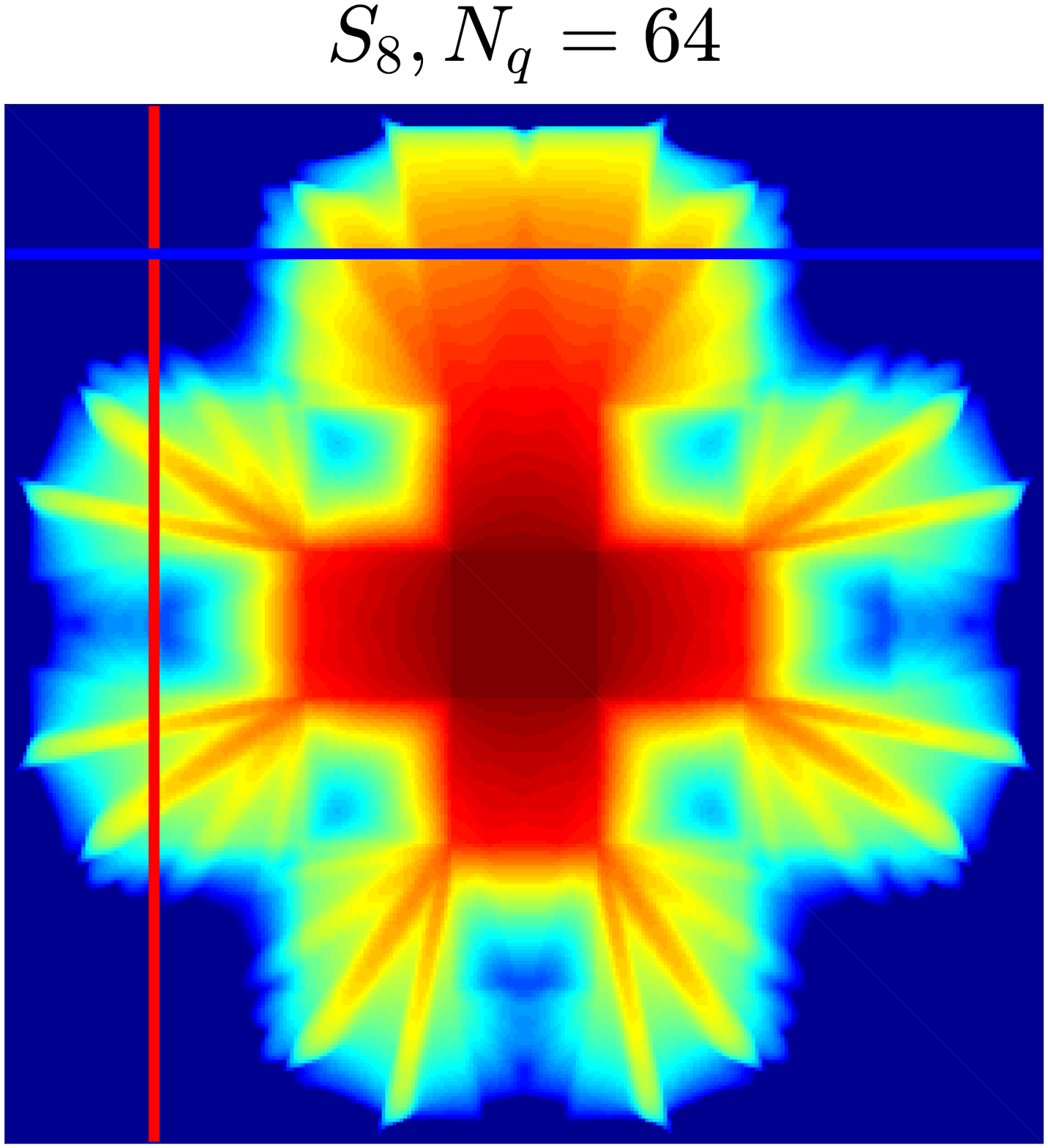}
		
		\label{fig:sub2}
	\end{subfigure}
	\begin{subfigure}{0.24\linewidth}
		\centering
		\includegraphics[scale=0.17]{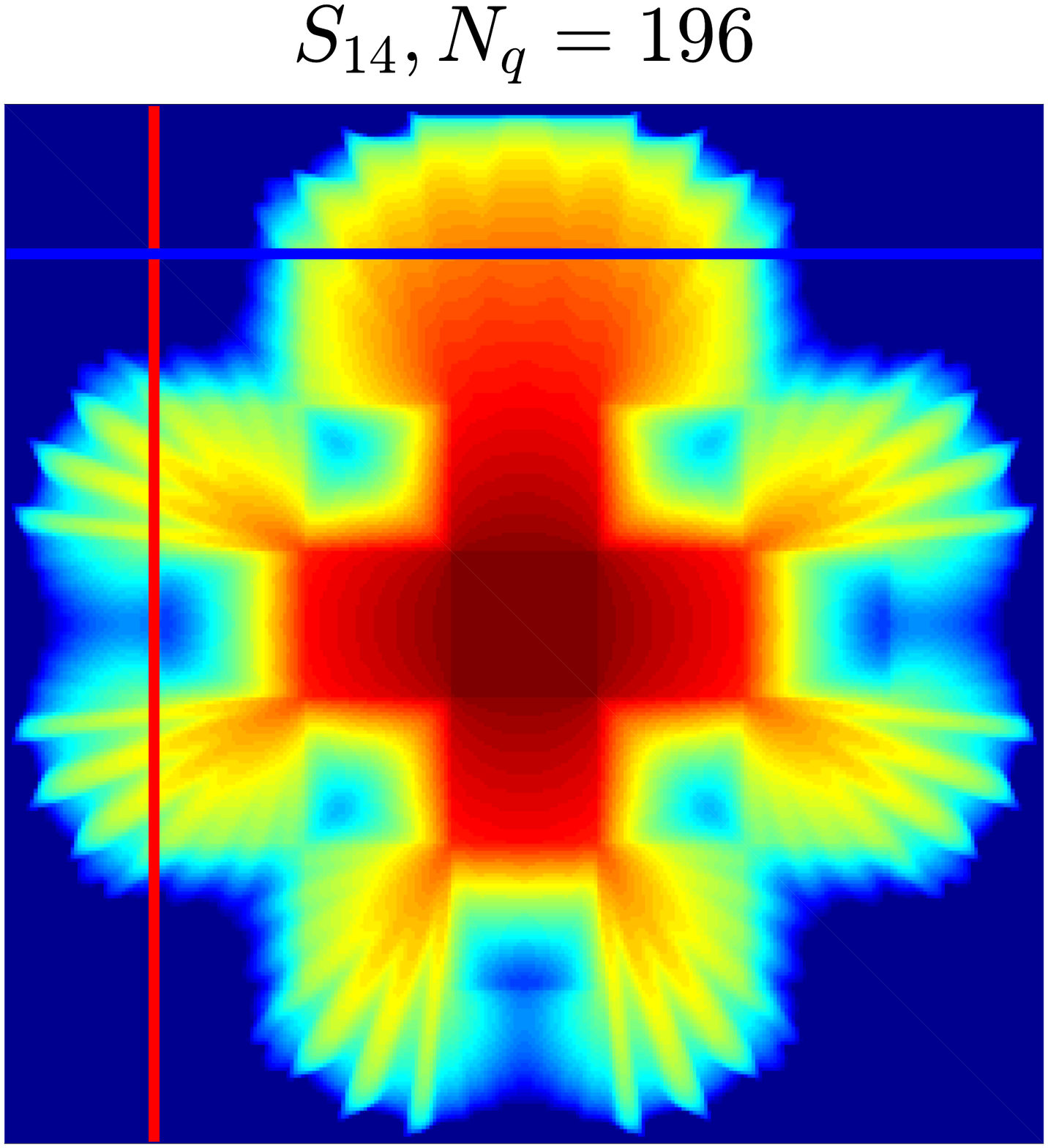}
		
		\label{fig:sub3}
	\end{subfigure}
	\begin{subfigure}{0.24\linewidth}
		\centering
		\includegraphics[scale=0.17]{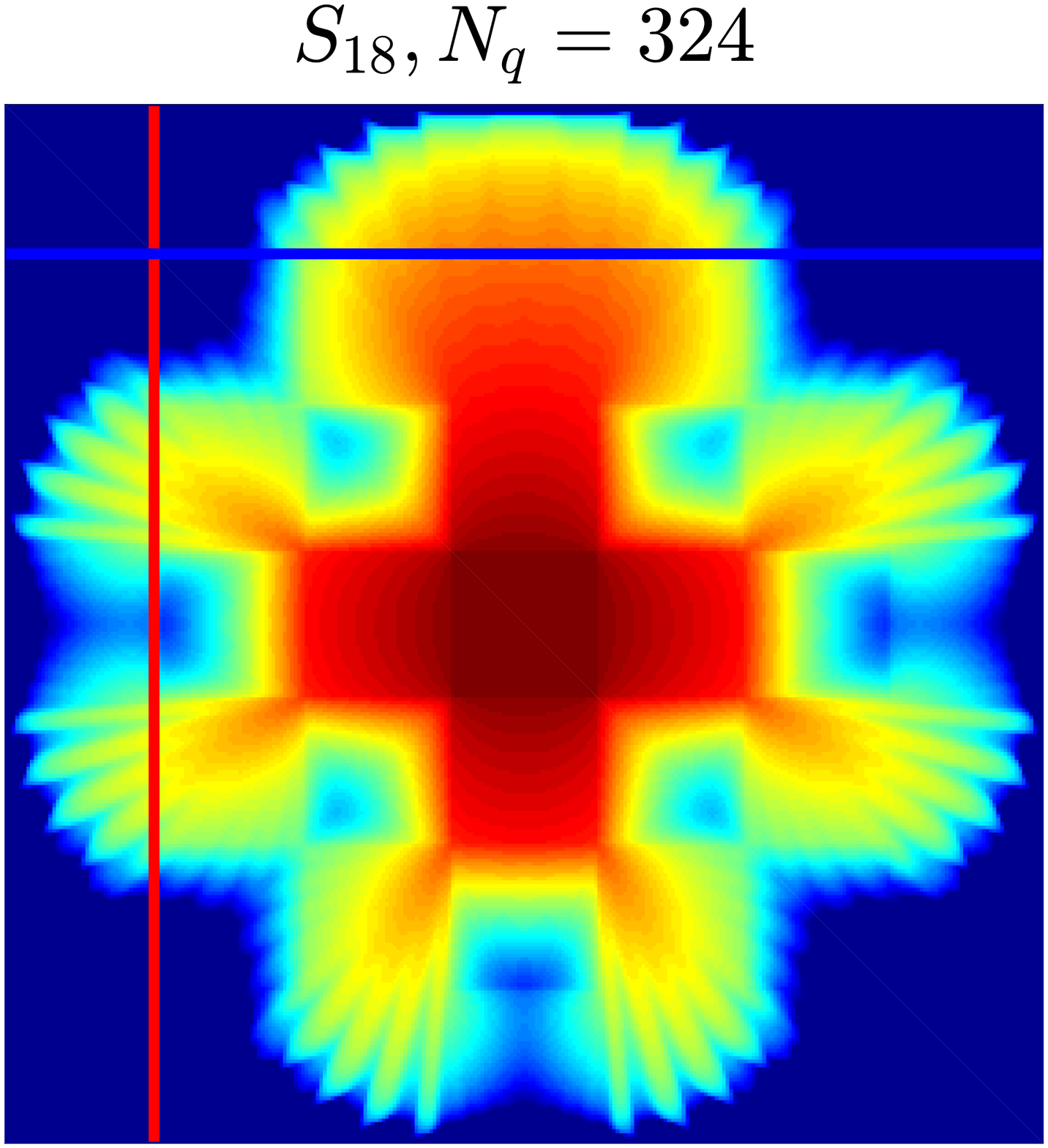}
		
		\label{fig:sub3}
	\end{subfigure}\\[-1ex]
	\begin{subfigure}{0.24\linewidth}
		\centering
		\includegraphics[scale=0.17]{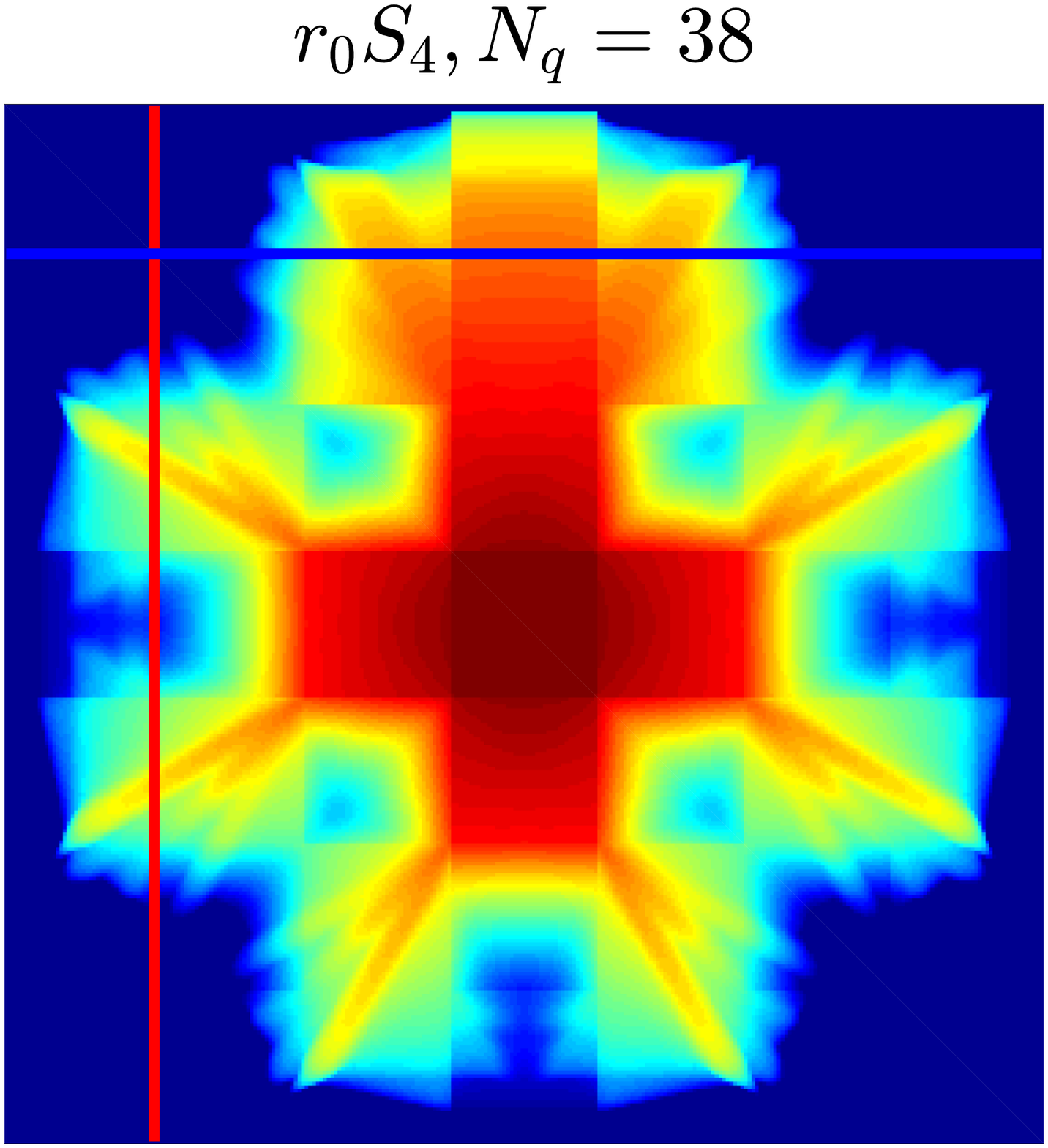}
		
		\label{fig:sub1}
	\end{subfigure}%
	\begin{subfigure}{0.24\linewidth}
		\centering
		\includegraphics[scale=0.17]{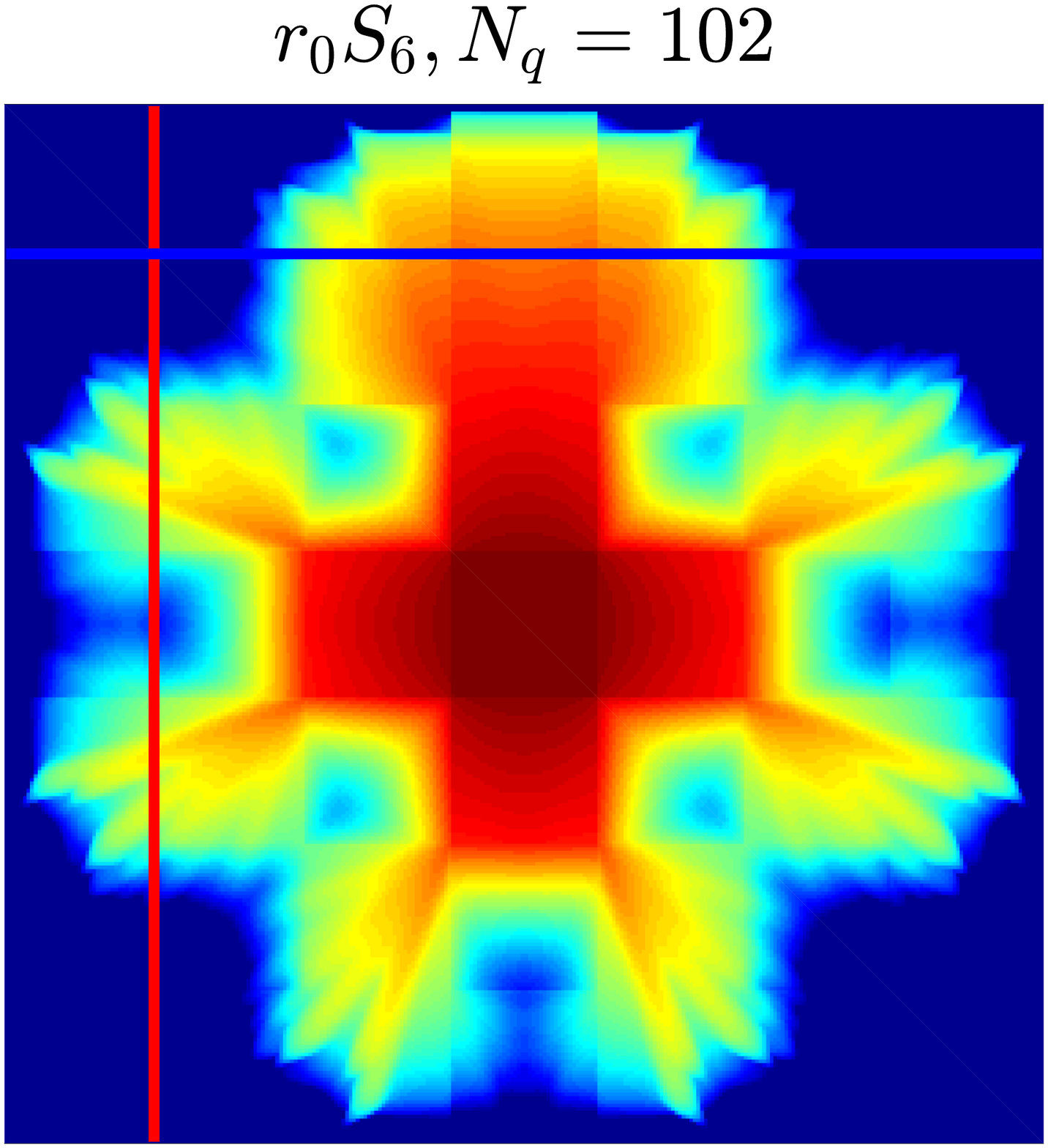}
		
		\label{fig:sub2}
	\end{subfigure}
	\begin{subfigure}{0.24\linewidth}
		\centering
		\includegraphics[scale=0.17]{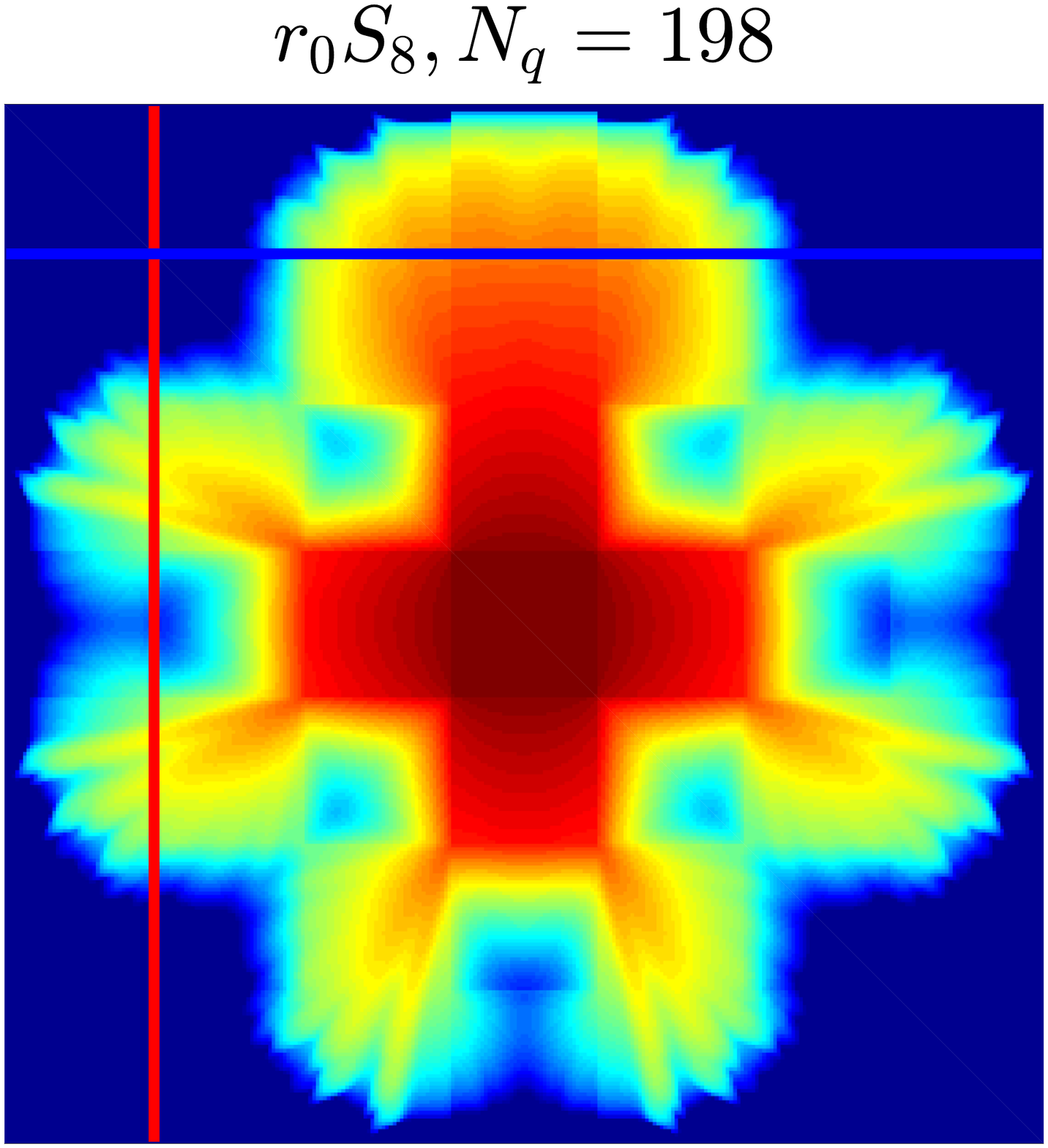}
		
		\label{fig:sub3}
	\end{subfigure}
	\begin{subfigure}{0.24\linewidth}
		\centering
		\includegraphics[scale=0.17]{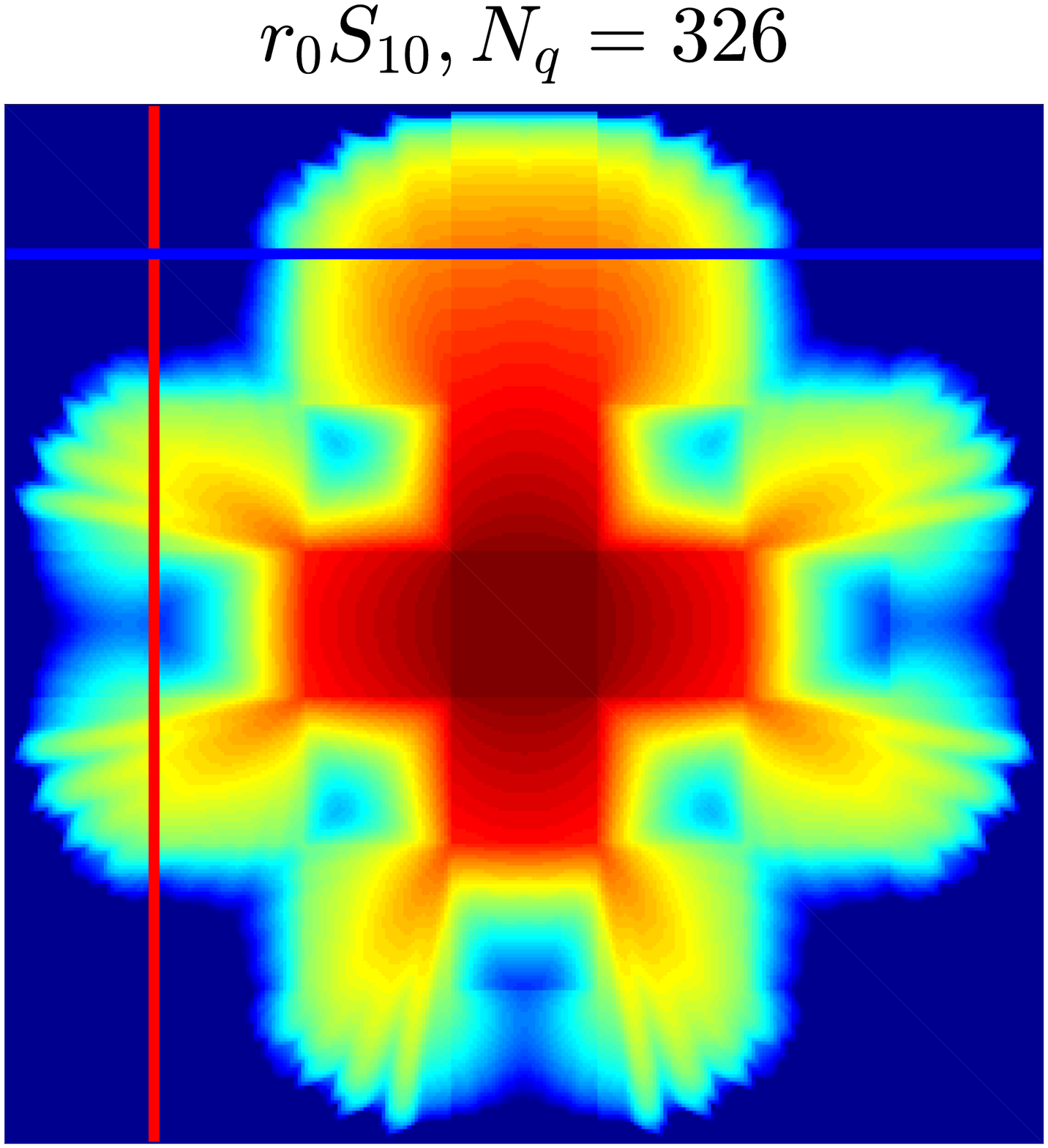}
		
		\label{fig:sub3}
	\end{subfigure}
	\\[-1ex]
	\begin{subfigure}{0.24\linewidth}
		\centering
		\includegraphics[scale=0.17]{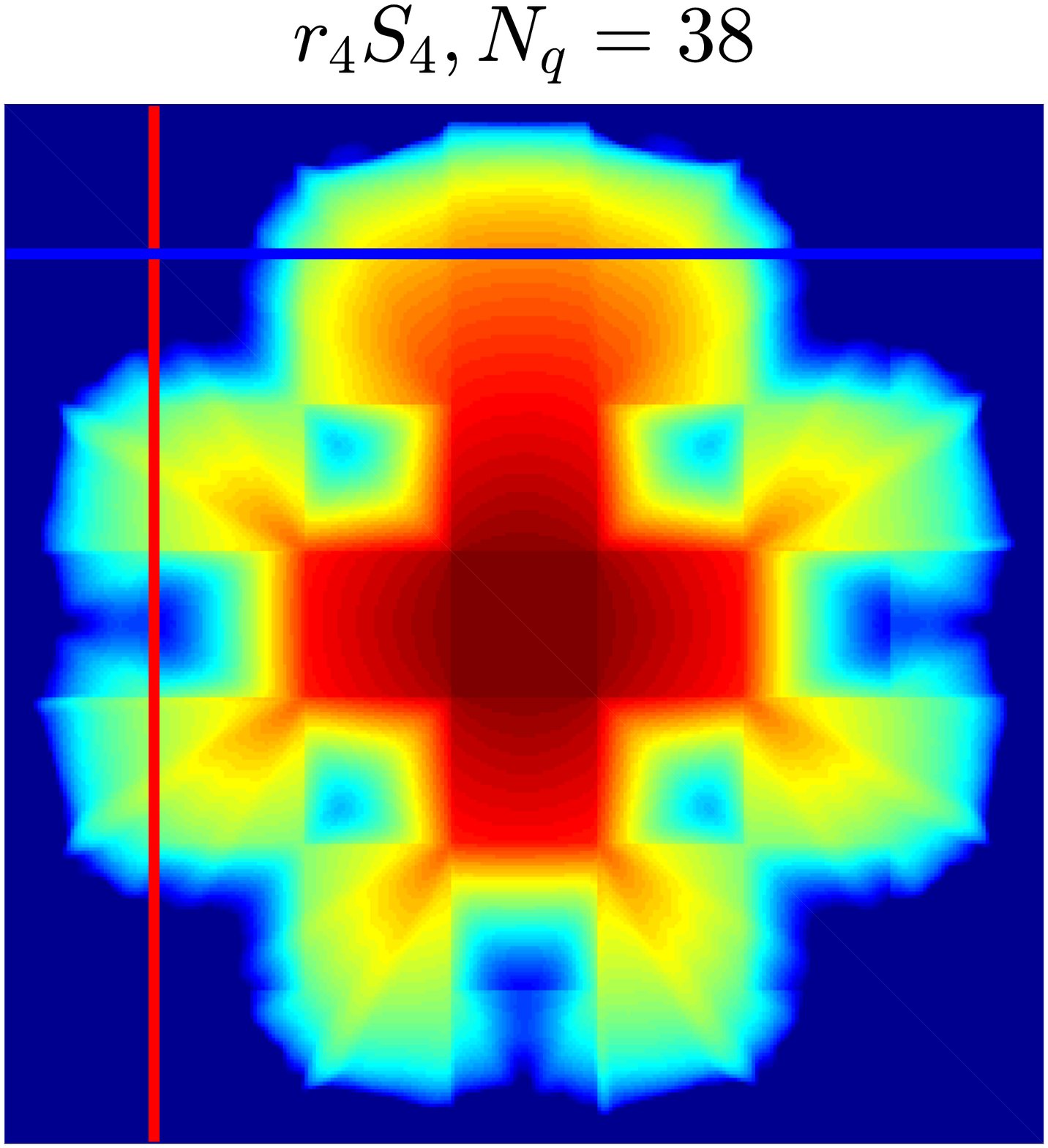}
		
		\label{fig:sub1}
	\end{subfigure}%
	\begin{subfigure}{0.24\linewidth}
		\centering
		\includegraphics[scale=0.17]{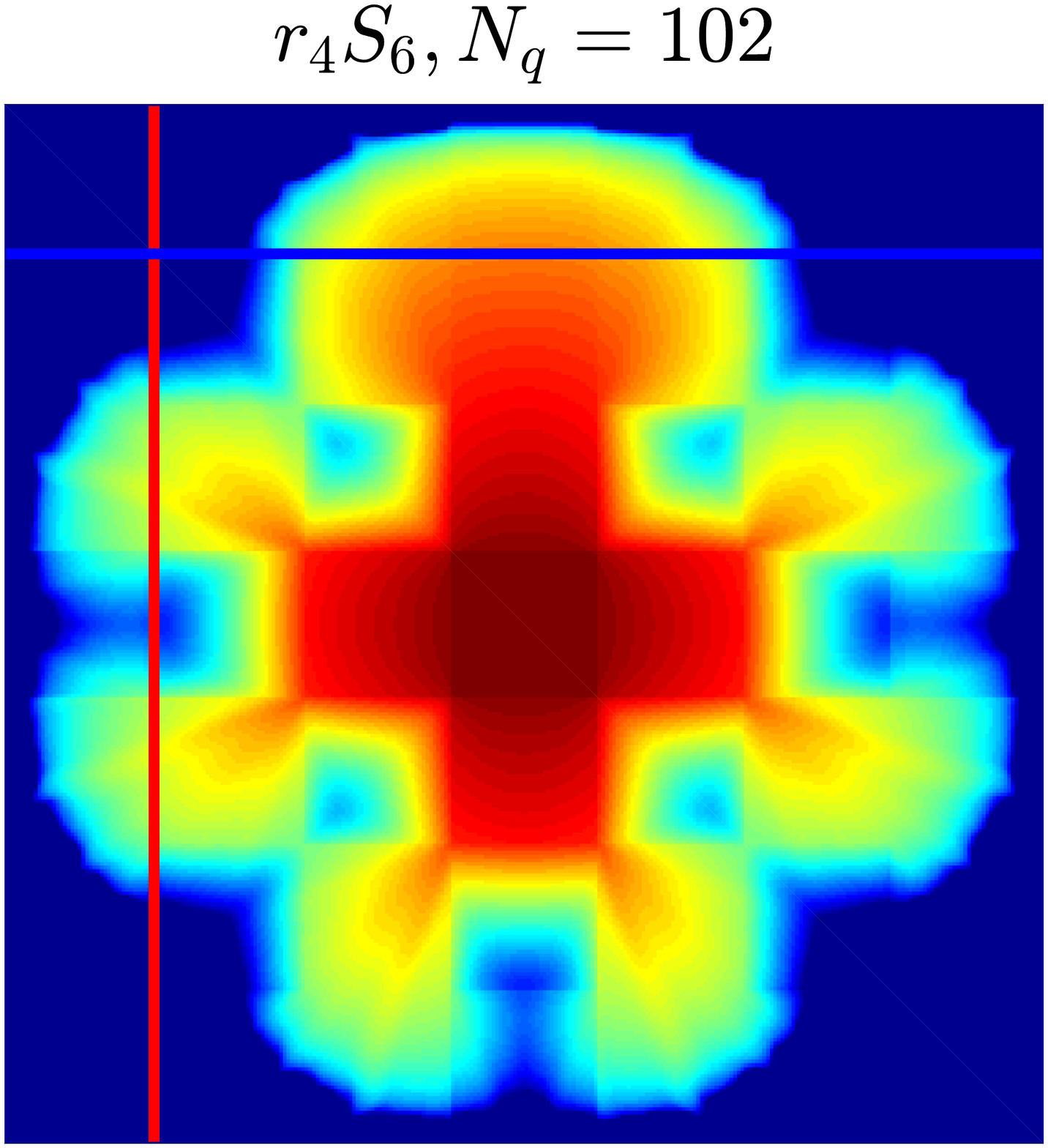}
		
		\label{fig:sub2}
	\end{subfigure}
	\begin{subfigure}{0.24\linewidth}
		\centering
		\includegraphics[scale=0.17]{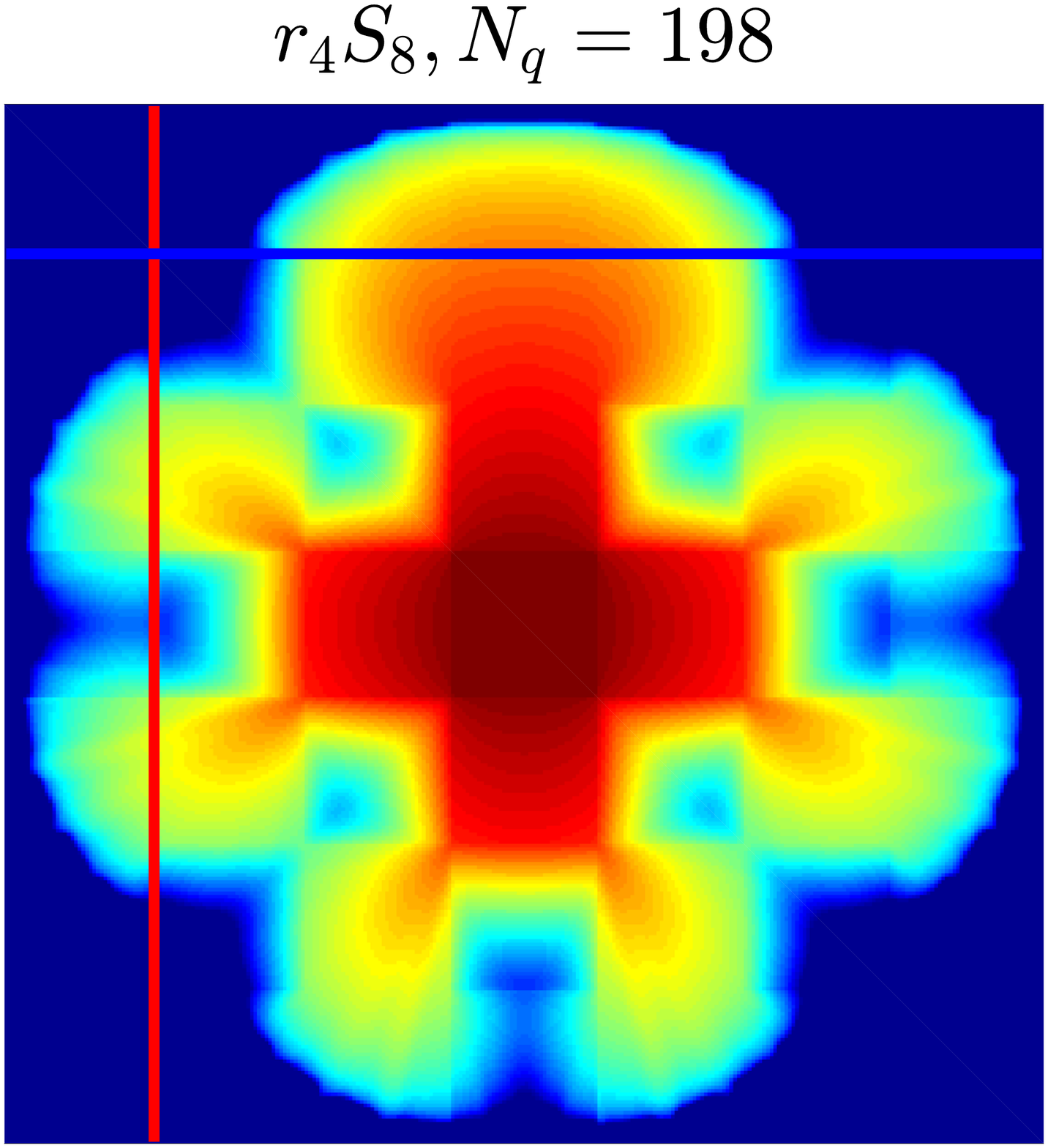}
		
		\label{fig:sub3}
	\end{subfigure}
	\begin{subfigure}{0.24\linewidth}
		\centering
		\includegraphics[scale=0.17]{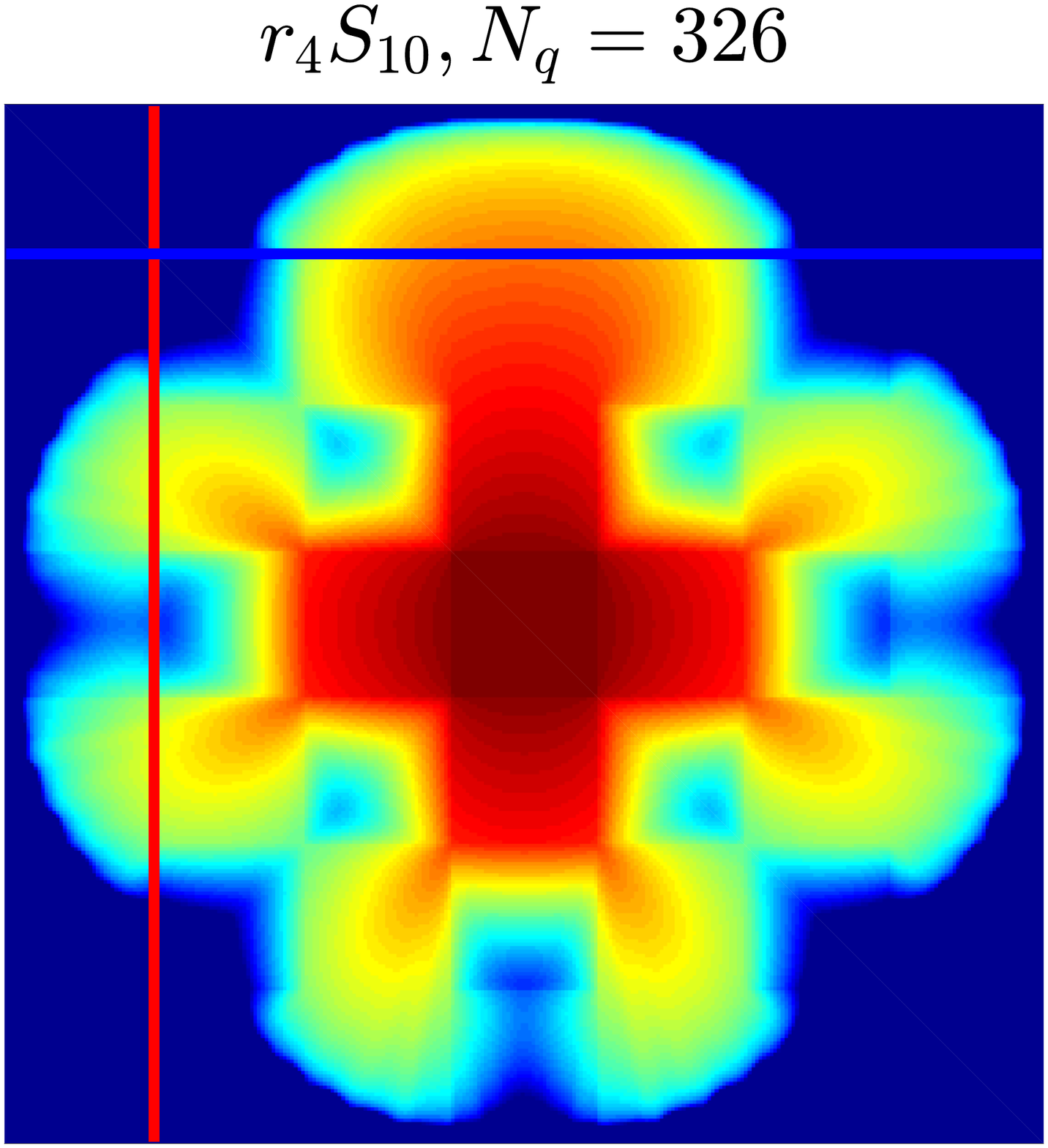}
		
		\label{fig:sub3}
	\end{subfigure}\\[-1ex]
	\begin{subfigure}{0.24\linewidth}
		\centering
		\includegraphics[scale=0.17]{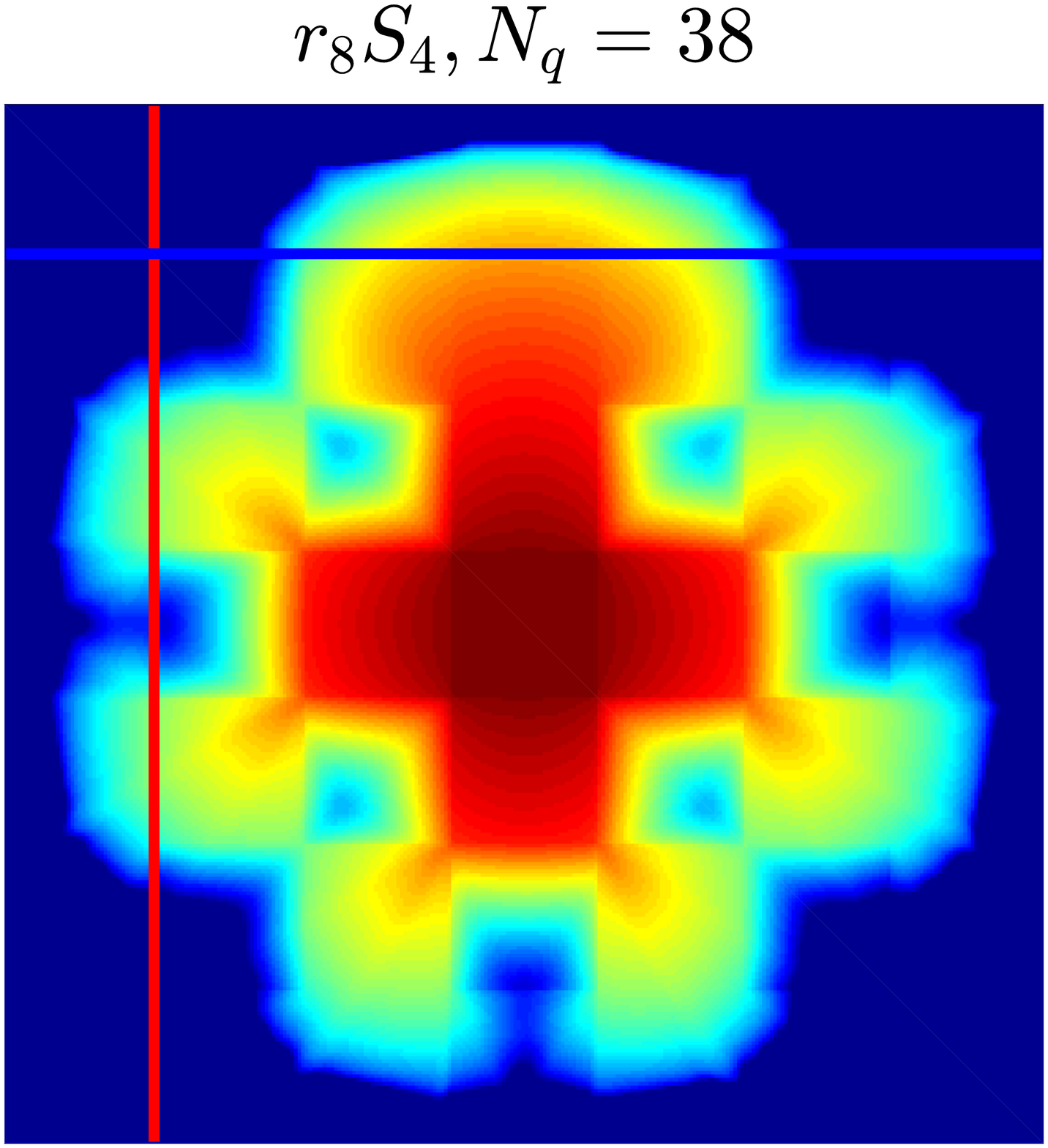}
		
		\label{fig:sub1}
	\end{subfigure}%
	\begin{subfigure}{0.24\linewidth}
		\centering
		\includegraphics[scale=0.17]{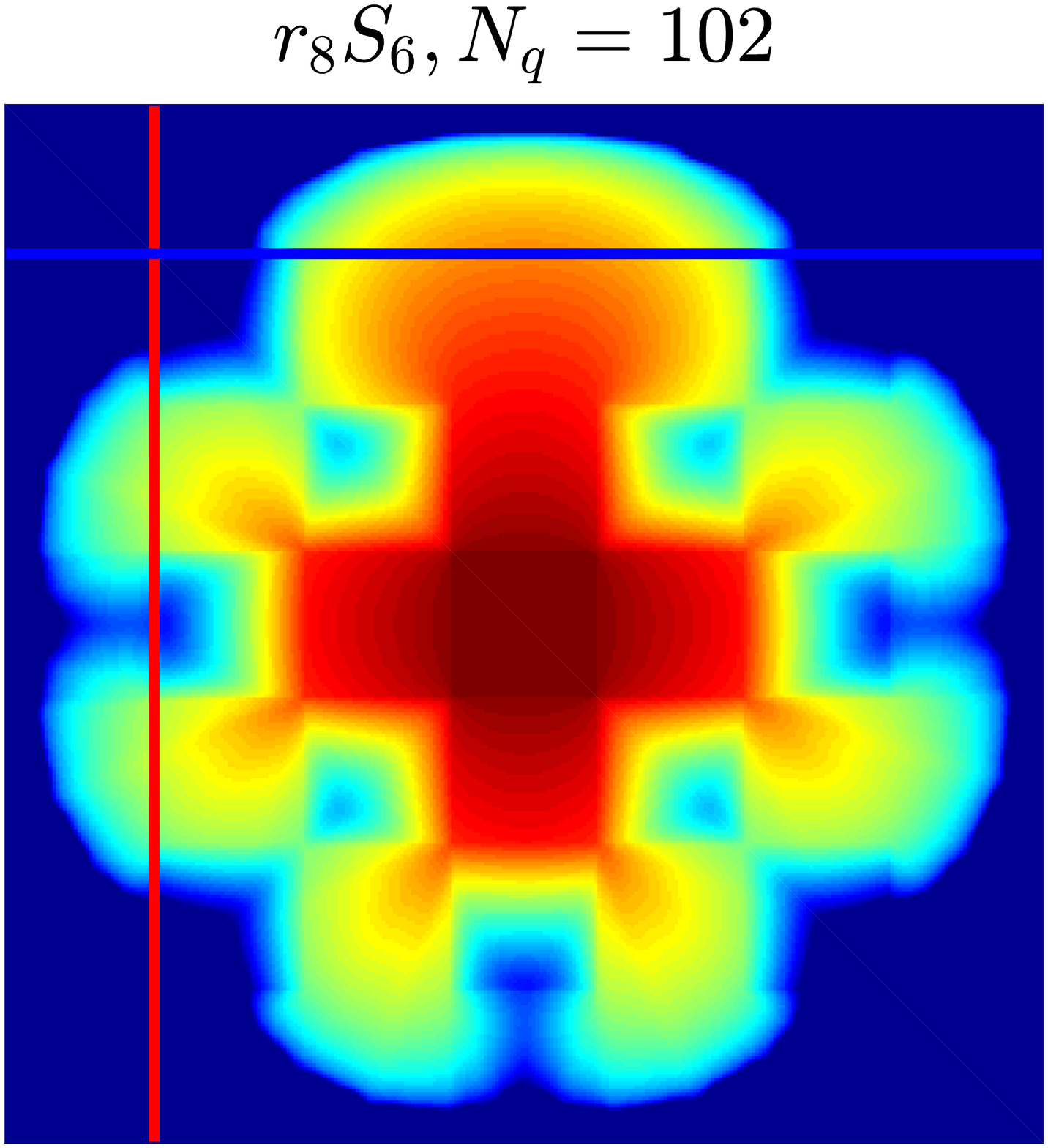}
		
		\label{fig:sub2}
	\end{subfigure}
	\begin{subfigure}{0.24\linewidth}
		\centering
		\includegraphics[scale=0.17]{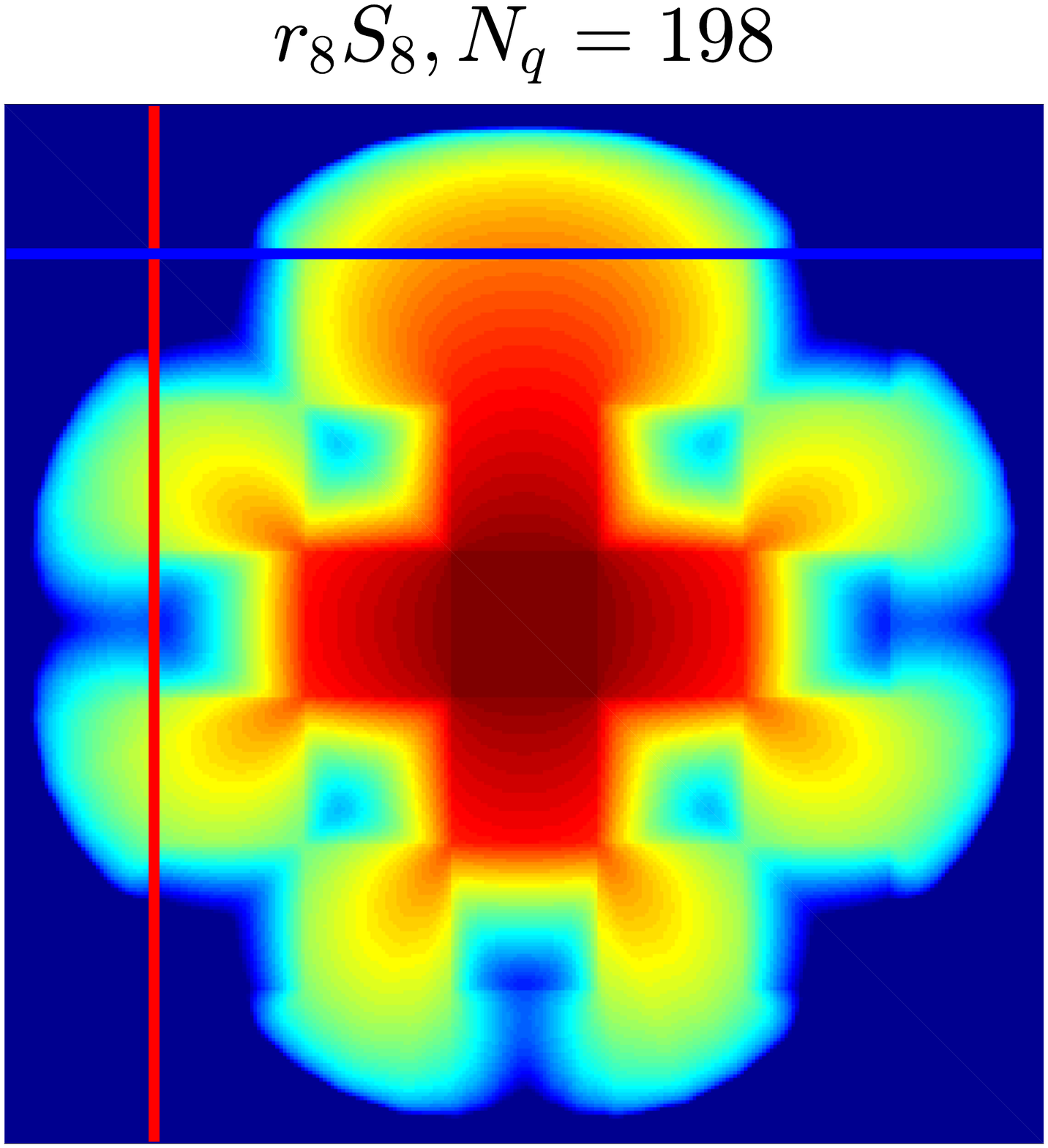}
		\label{fig:sub3}
	\end{subfigure}
	\begin{subfigure}{0.24\linewidth}
		\centering
		\includegraphics[scale=0.17]{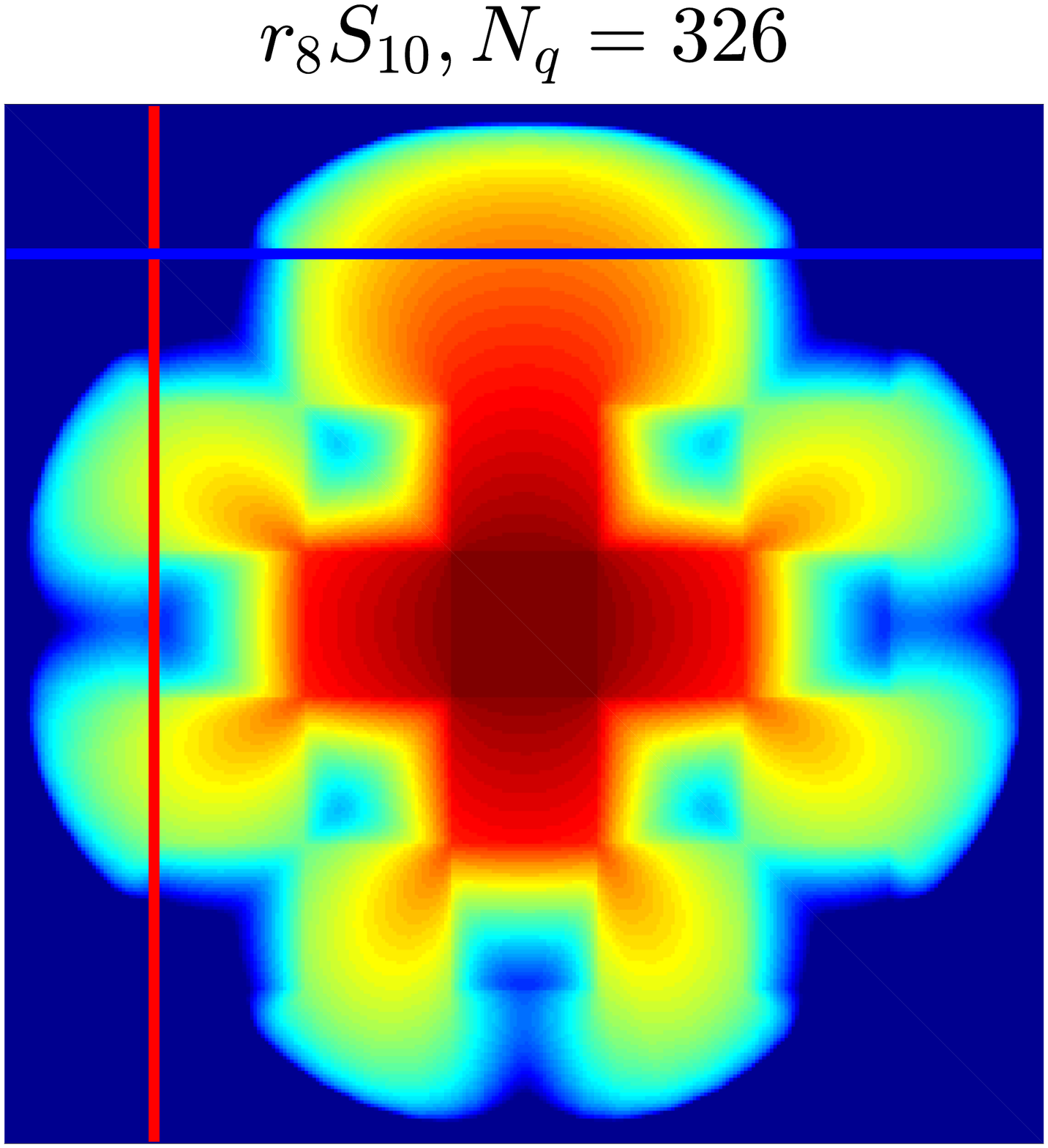}
		\label{fig:sub3}
	\end{subfigure} 
	\caption{Logarithmic density for the lattice problem.}
	\label{fig:checkerboard}
\end{figure}

\begin{figure}
	
	\begin{subfigure}{0.24\linewidth}
		\centering
		\includegraphics[scale=0.17]{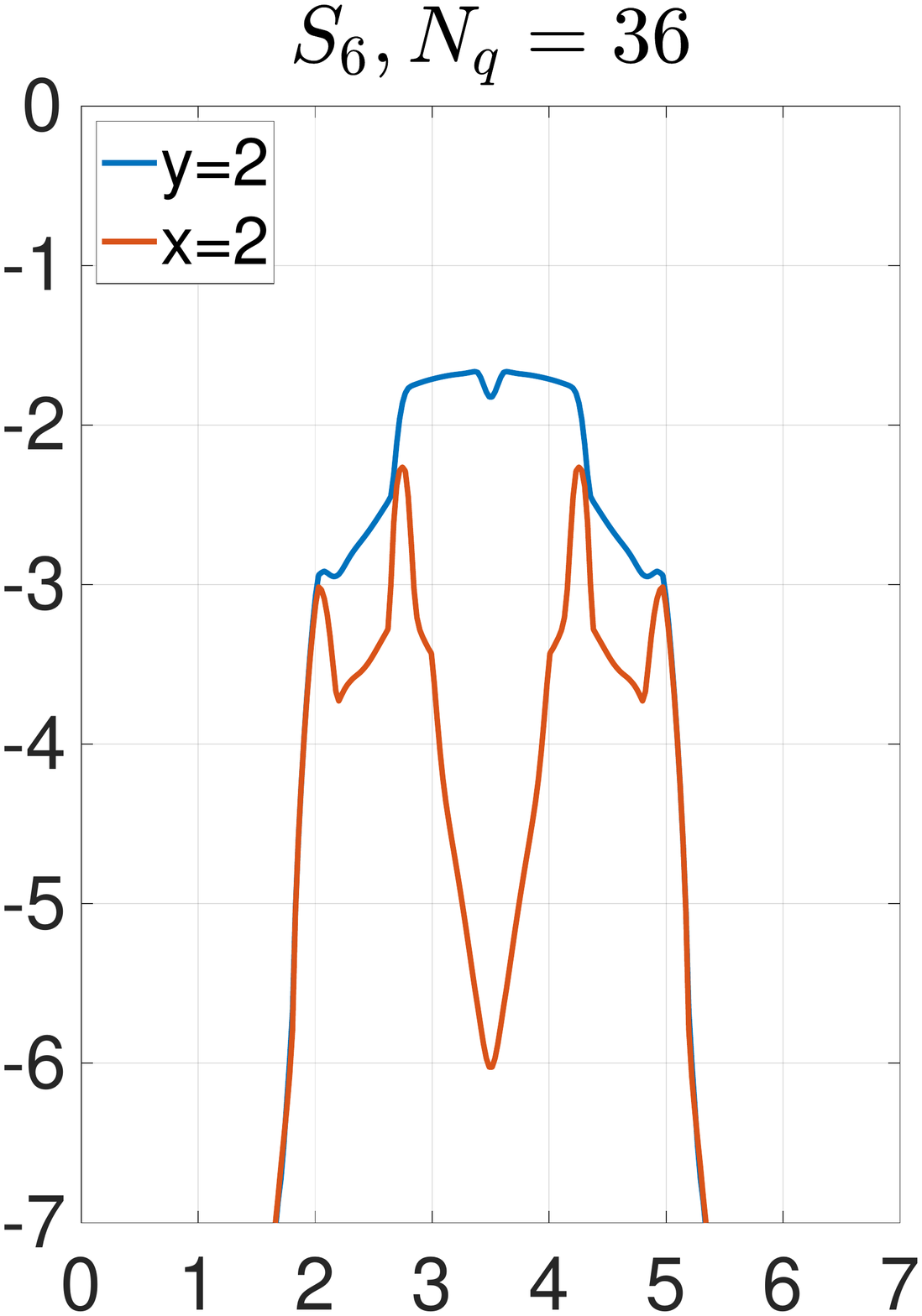}
		
		\label{fig:sub1}
	\end{subfigure}%
	\begin{subfigure}{0.24\linewidth}
		\centering
		\includegraphics[scale=0.17]{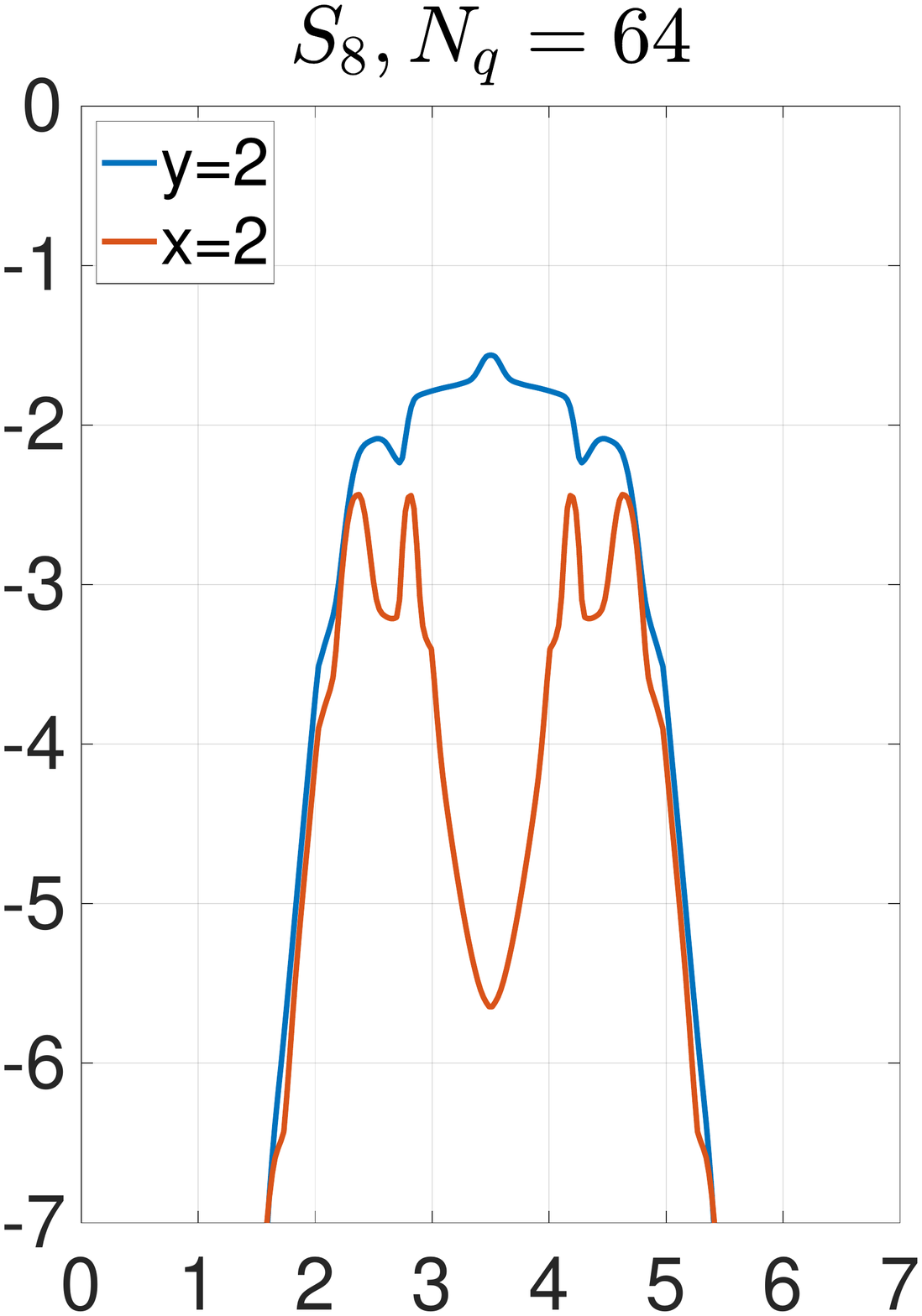}
		
		\label{fig:sub2}
	\end{subfigure}
	\begin{subfigure}{0.24\linewidth}
		\centering
		\includegraphics[scale=0.17]{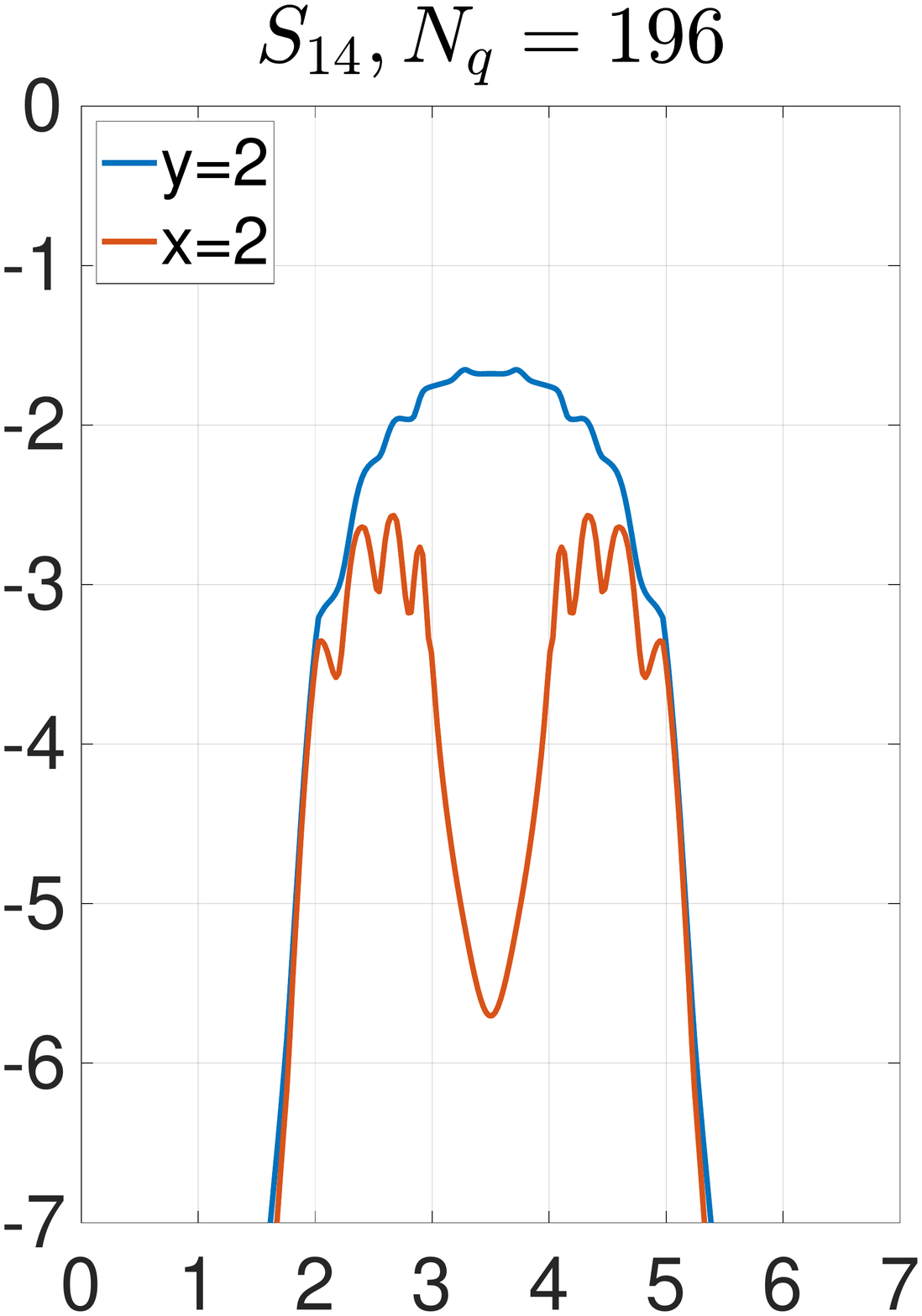}
		
		\label{fig:sub3}
	\end{subfigure}
	\begin{subfigure}{0.24\linewidth}
		\centering
		\includegraphics[scale=0.17]{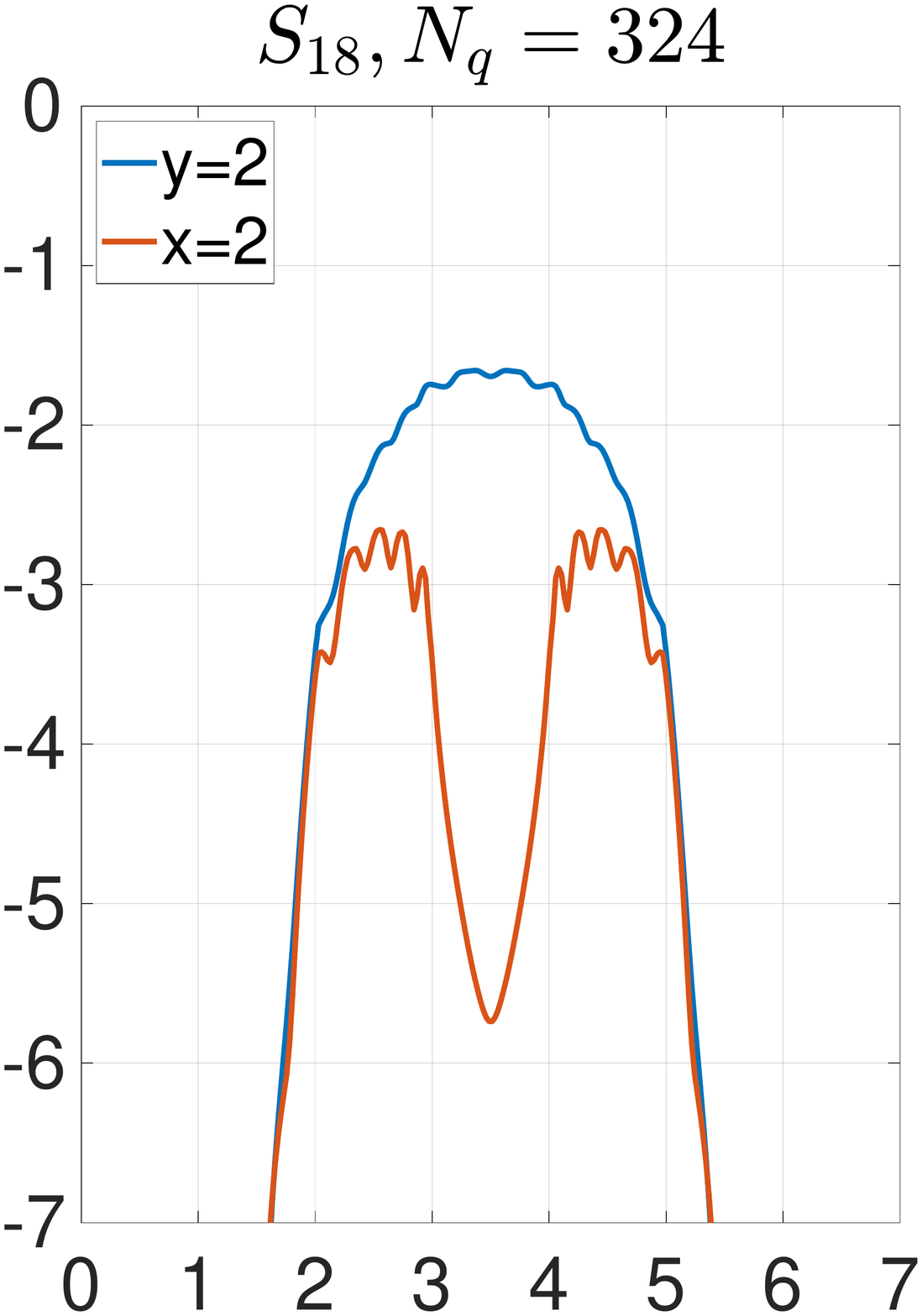}
		
		\label{fig:sub3}
	\end{subfigure}\\[-1ex]
	\begin{subfigure}{0.24\linewidth}
		\centering
		\includegraphics[scale=0.17]{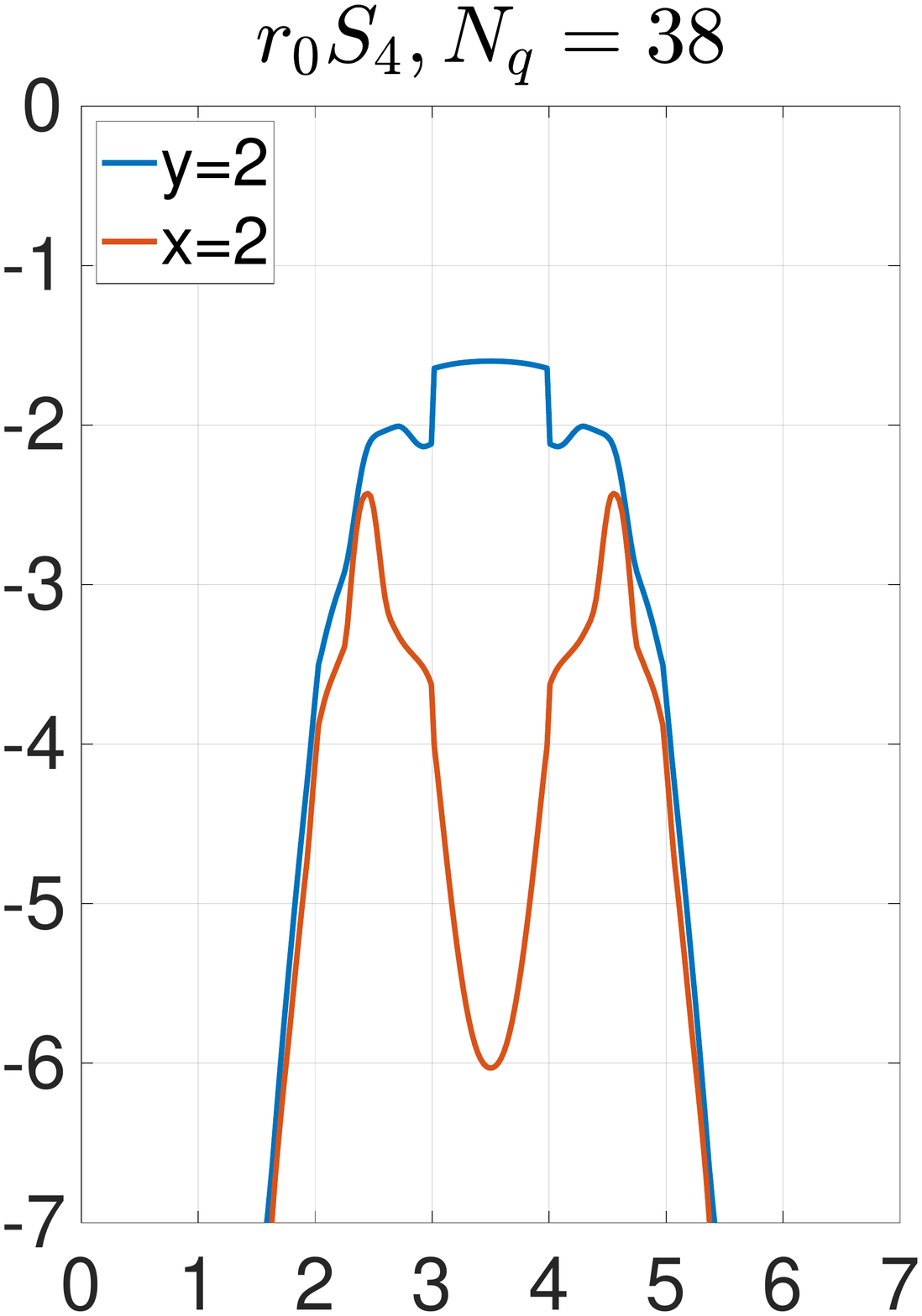}
		
		\label{fig:sub1}
	\end{subfigure}%
	\begin{subfigure}{0.24\linewidth}
		\centering
		\includegraphics[scale=0.17]{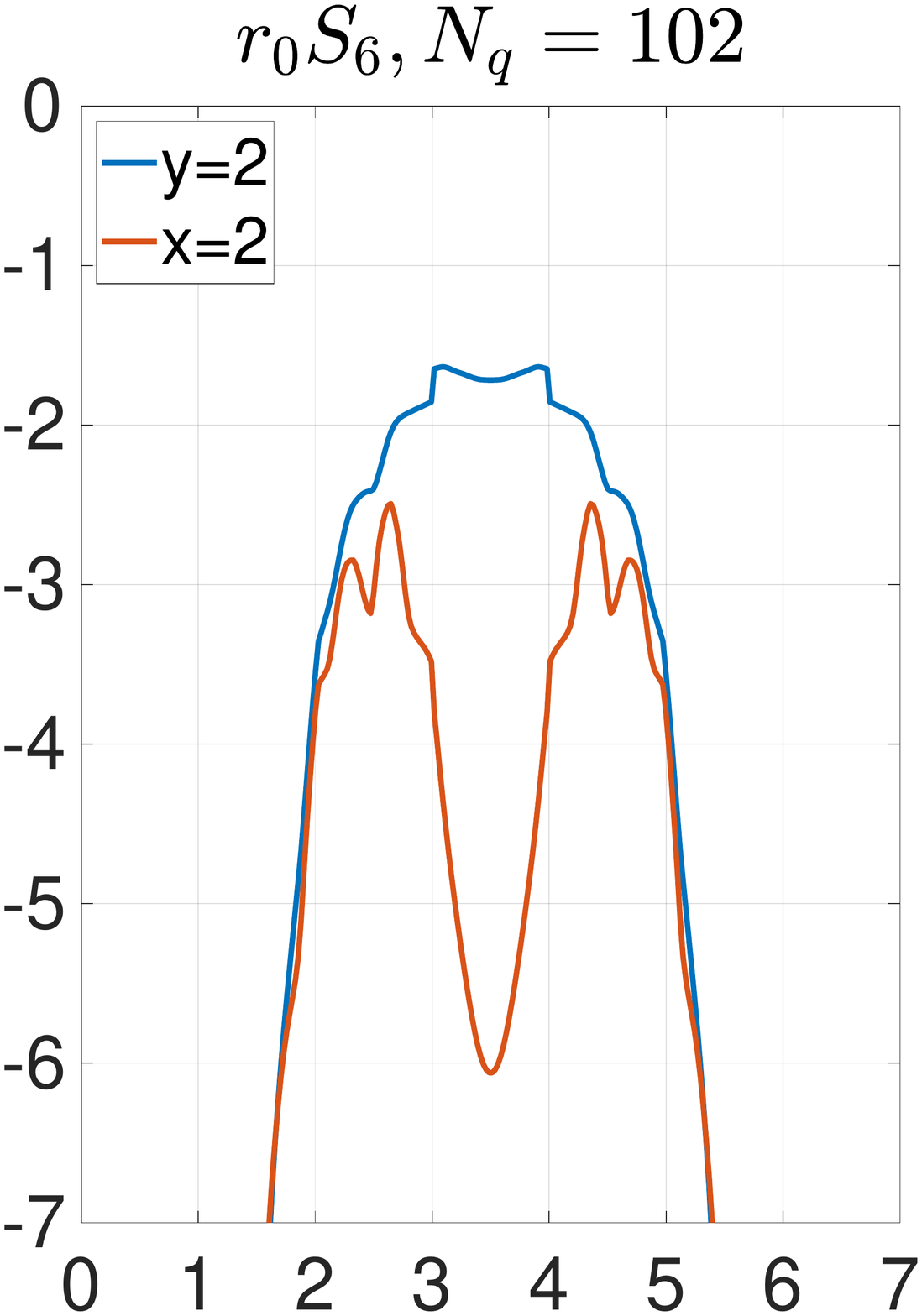}
		
		\label{fig:sub2}
	\end{subfigure}
	\begin{subfigure}{0.24\linewidth}
		\centering
		\includegraphics[scale=0.17]{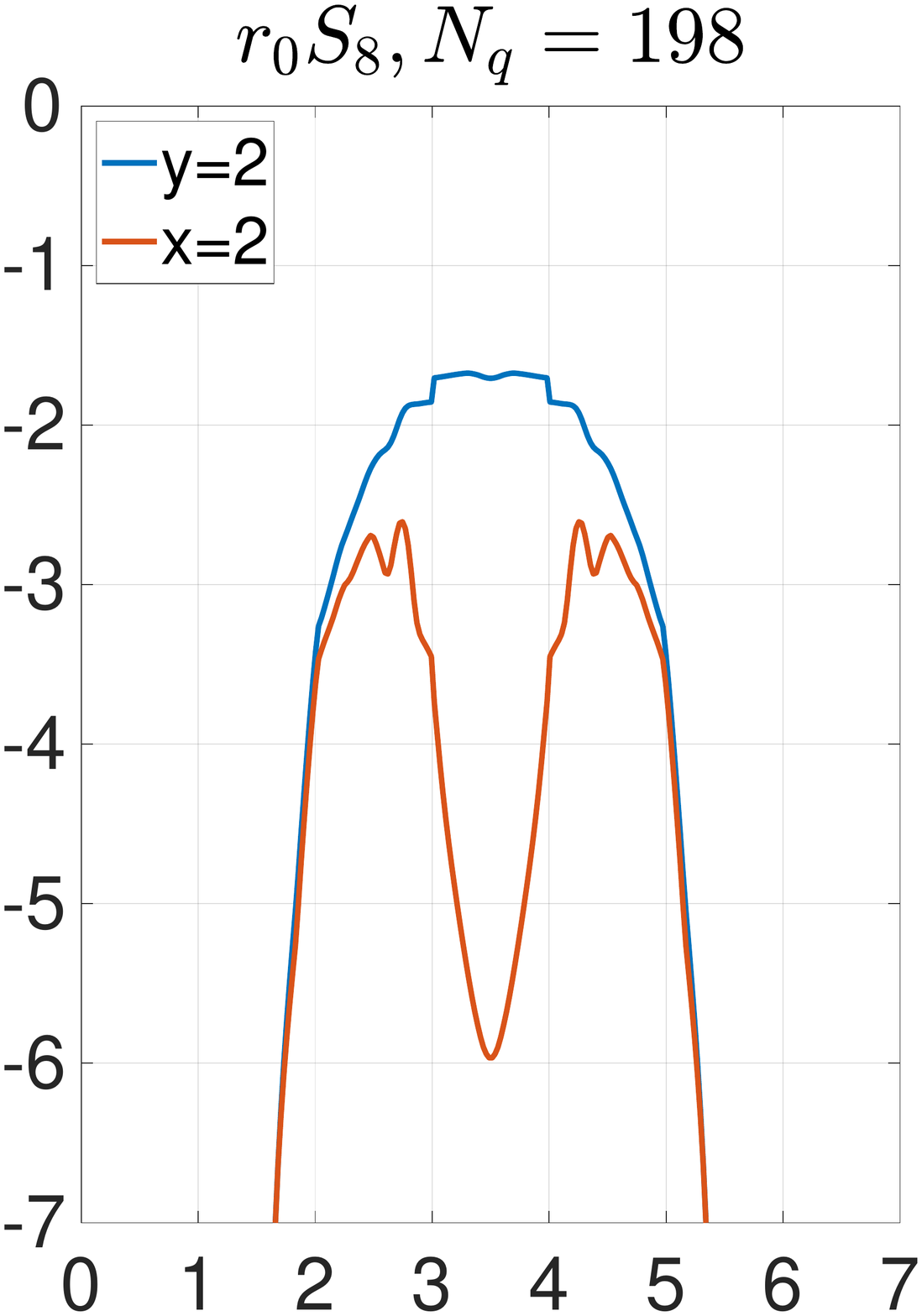}
		
		\label{fig:sub3}
	\end{subfigure}
	\begin{subfigure}{0.24\linewidth}
		\centering
		\includegraphics[scale=0.17]{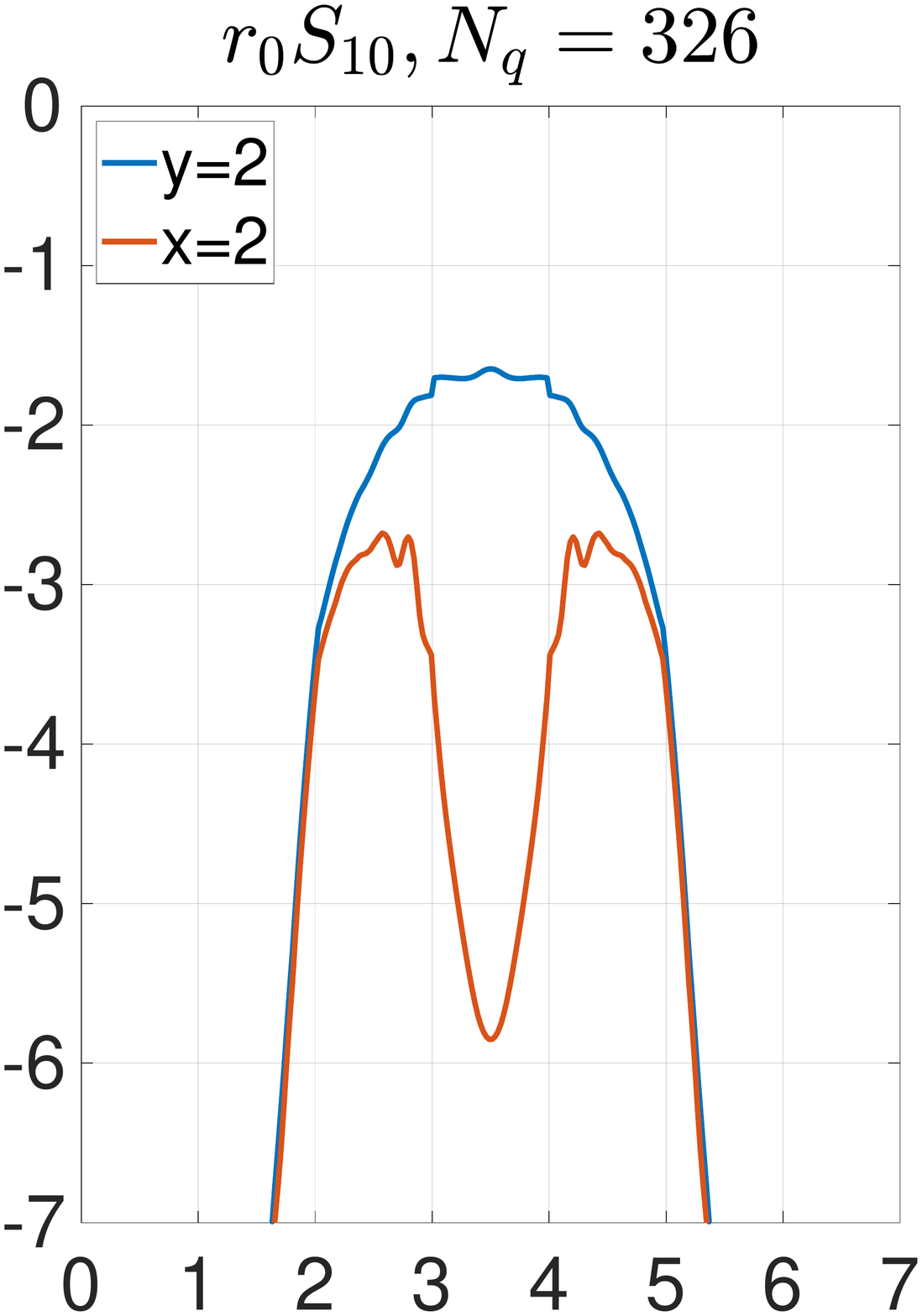}
		
		\label{fig:sub3}
	\end{subfigure}
	\\[-1ex]
	\begin{subfigure}{0.24\linewidth}
		\centering
		\includegraphics[scale=0.17]{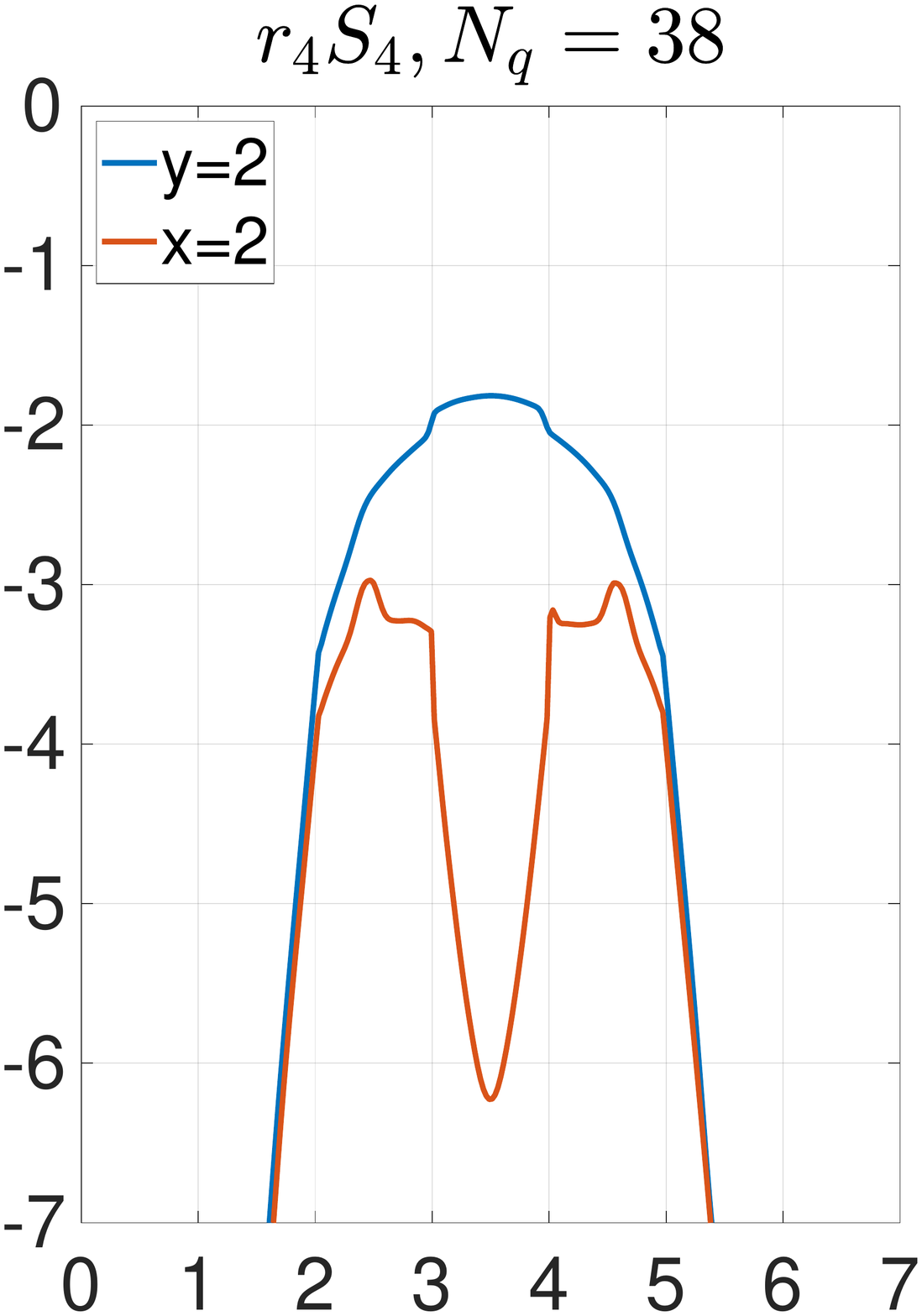}
		
		\label{fig:sub1}
	\end{subfigure}%
	\begin{subfigure}{0.24\linewidth}
		\centering
		\includegraphics[scale=0.17]{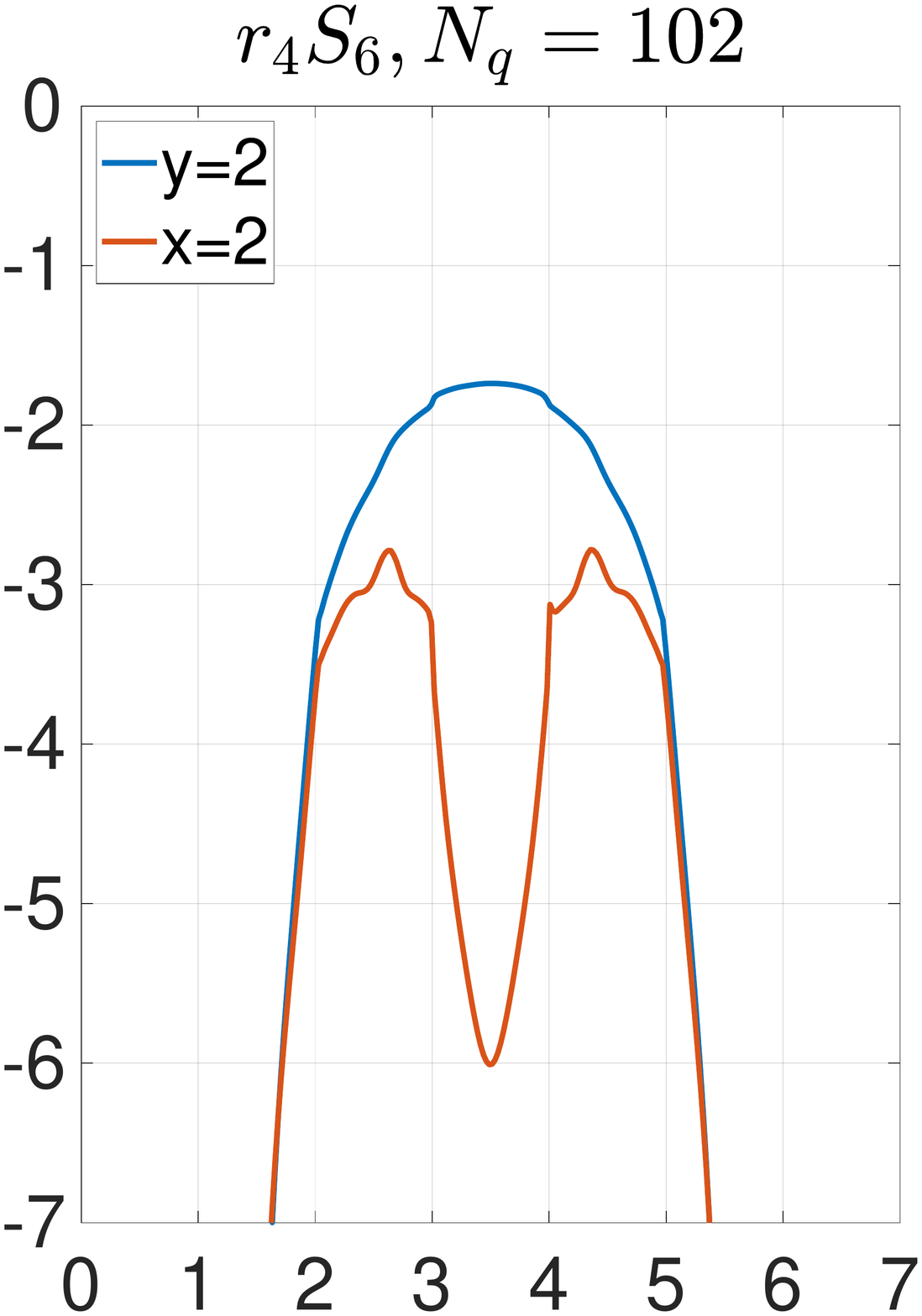}
		
		\label{fig:sub2}
	\end{subfigure}
	\begin{subfigure}{0.24\linewidth}
		\centering
		\includegraphics[scale=0.17]{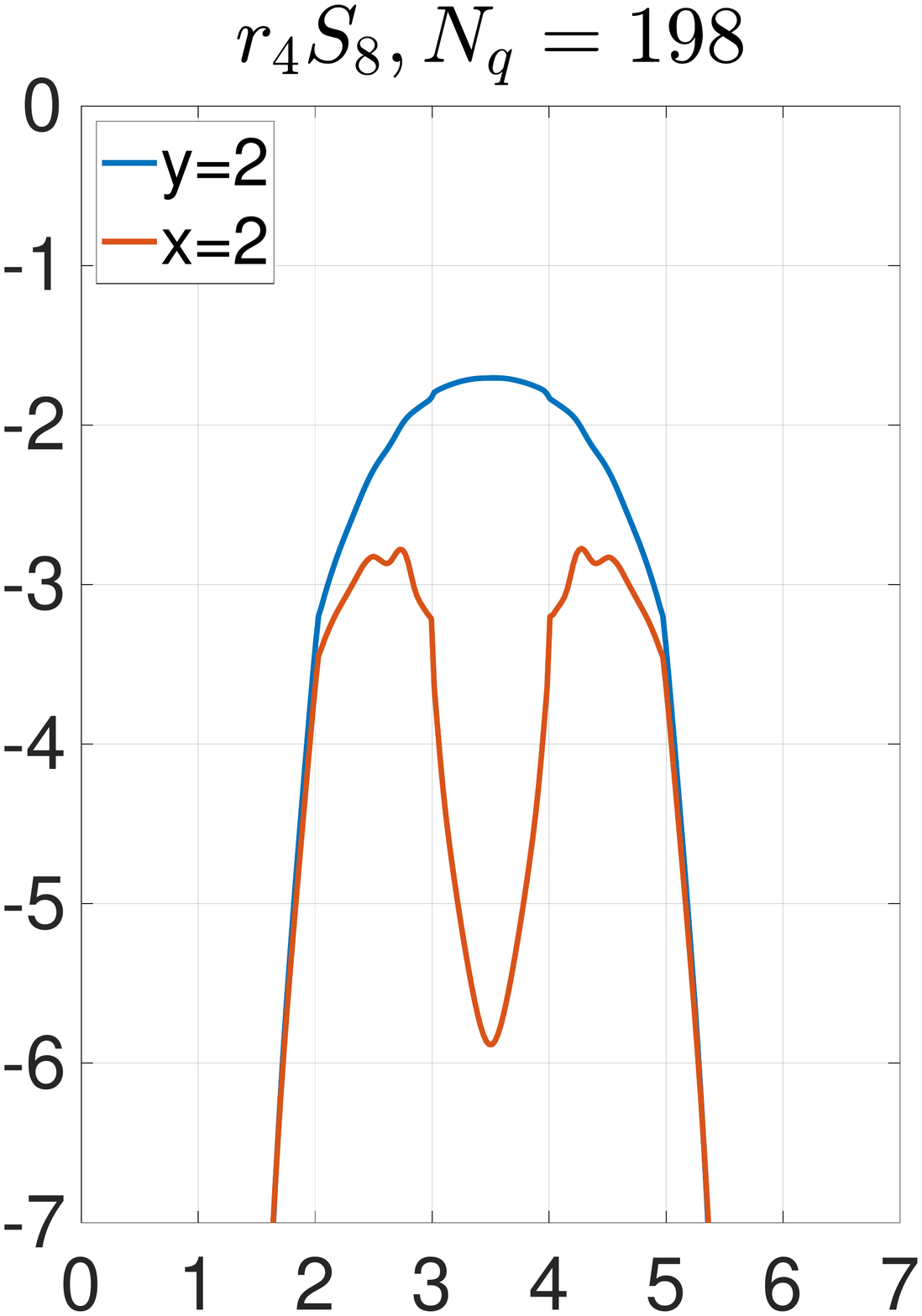}
		
		\label{fig:sub3}
	\end{subfigure}
	\begin{subfigure}{0.24\linewidth}
		\centering
		\includegraphics[scale=0.17]{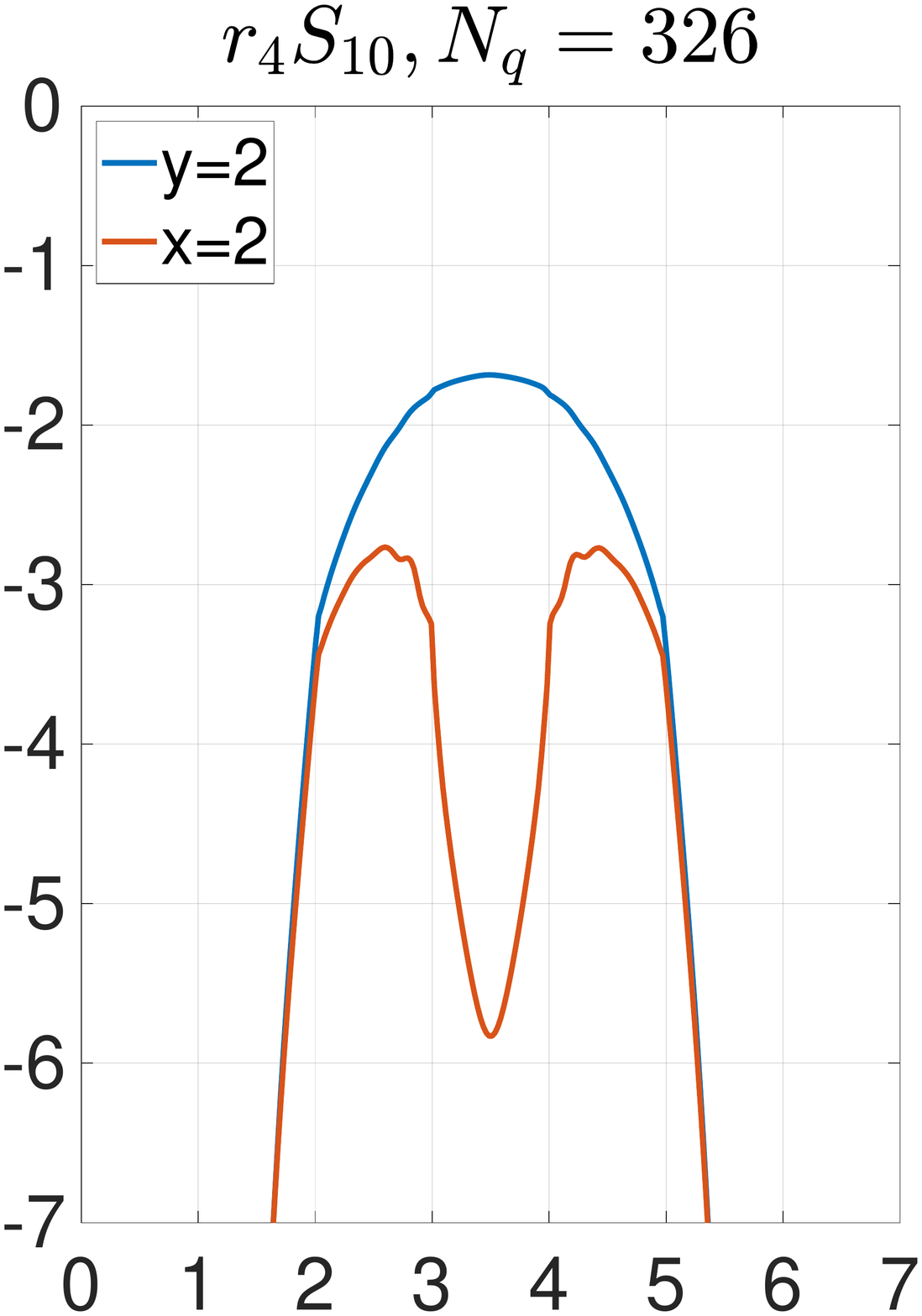}
		
		\label{fig:sub3}
	\end{subfigure}\\[-1ex]
	\begin{subfigure}{0.24\linewidth}
		\centering
		\includegraphics[scale=0.17]{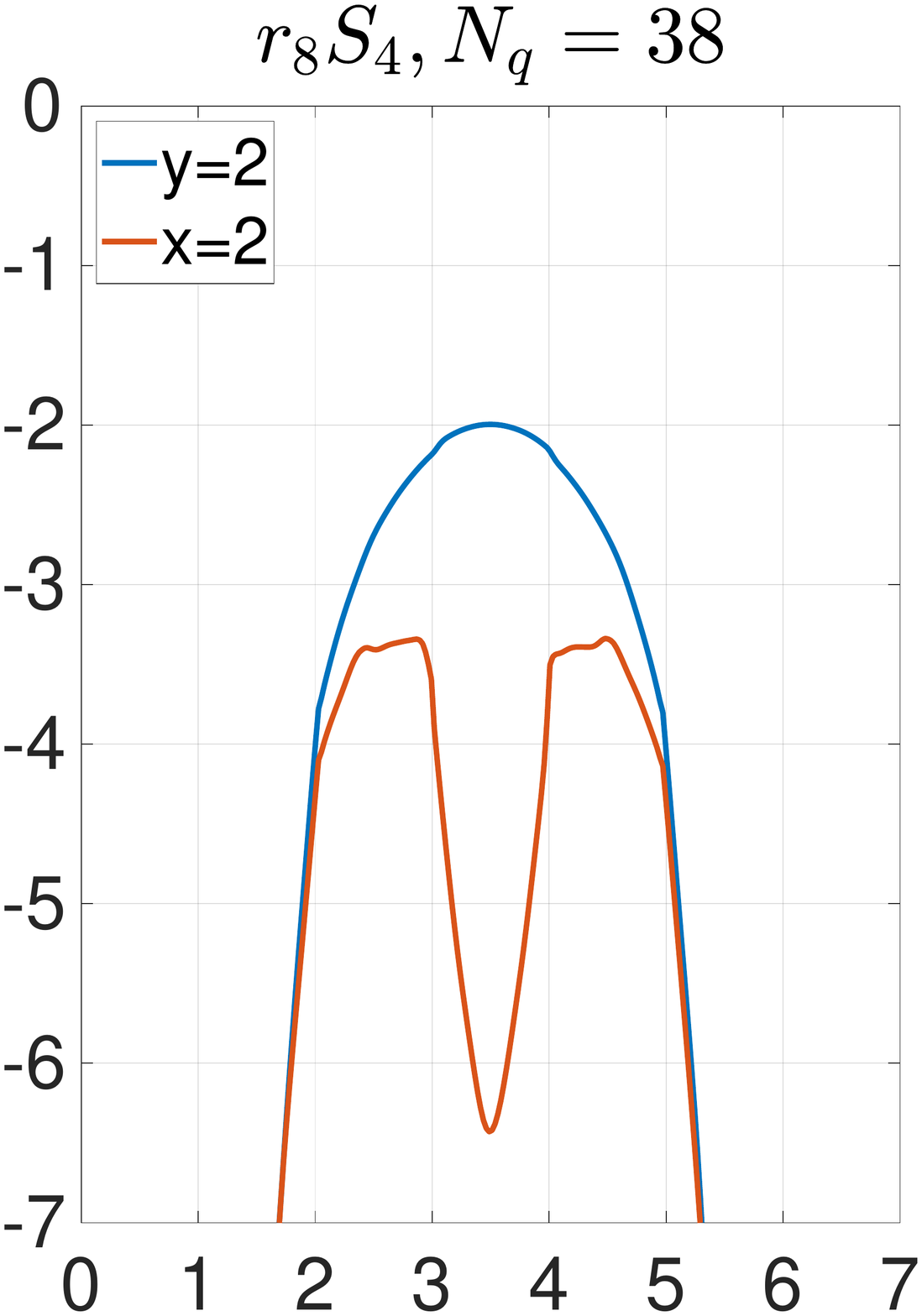}
		
		\label{fig:sub1}
	\end{subfigure}%
	\begin{subfigure}{0.24\linewidth}
		\centering
		\includegraphics[scale=0.17]{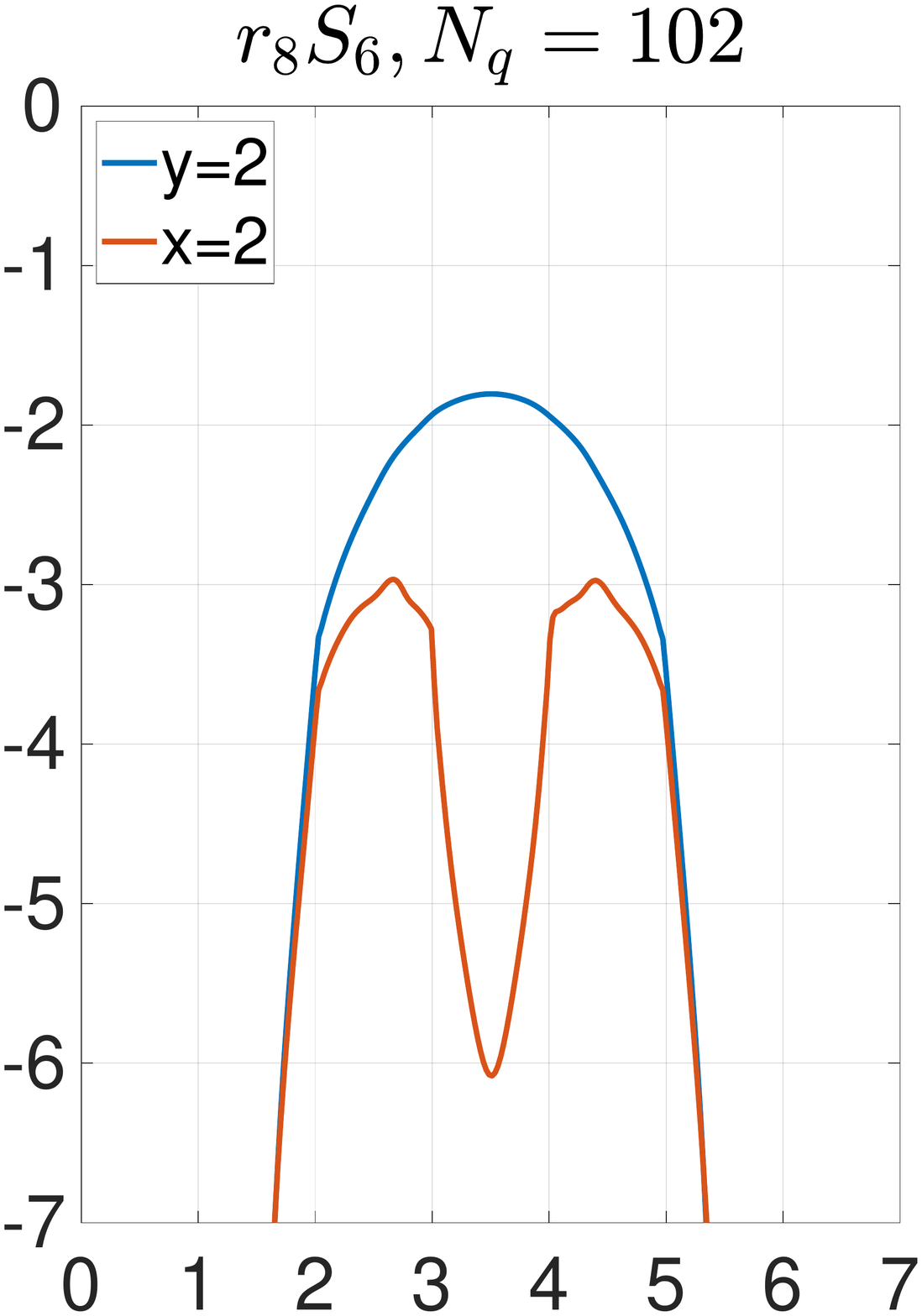}
		
		\label{fig:sub2}
	\end{subfigure}
	\begin{subfigure}{0.24\linewidth}
		\centering
		\includegraphics[scale=0.17]{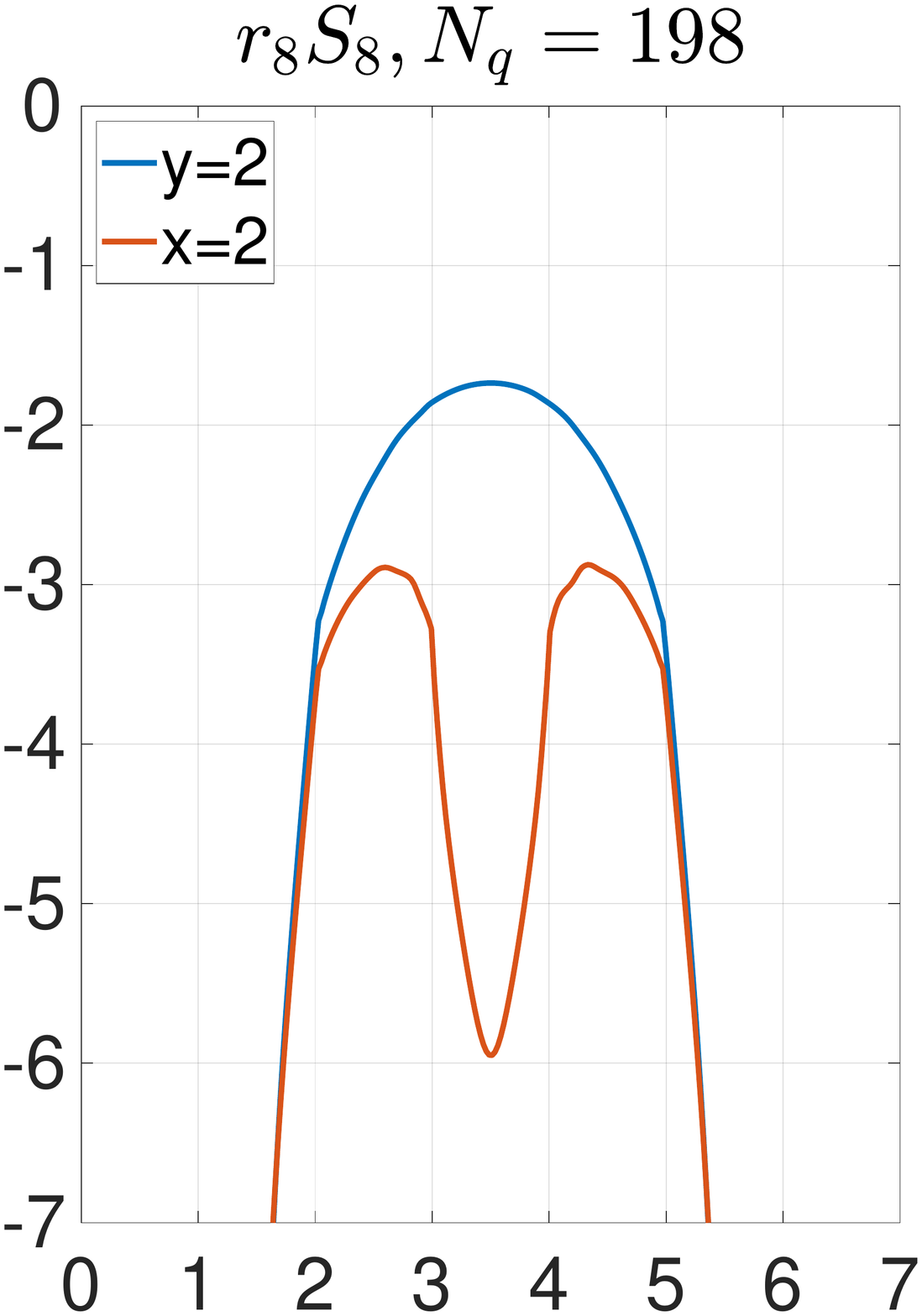}
		
		\label{fig:sub3}
	\end{subfigure}
	\begin{subfigure}{0.24\linewidth}
		\centering
		\includegraphics[scale=0.17]{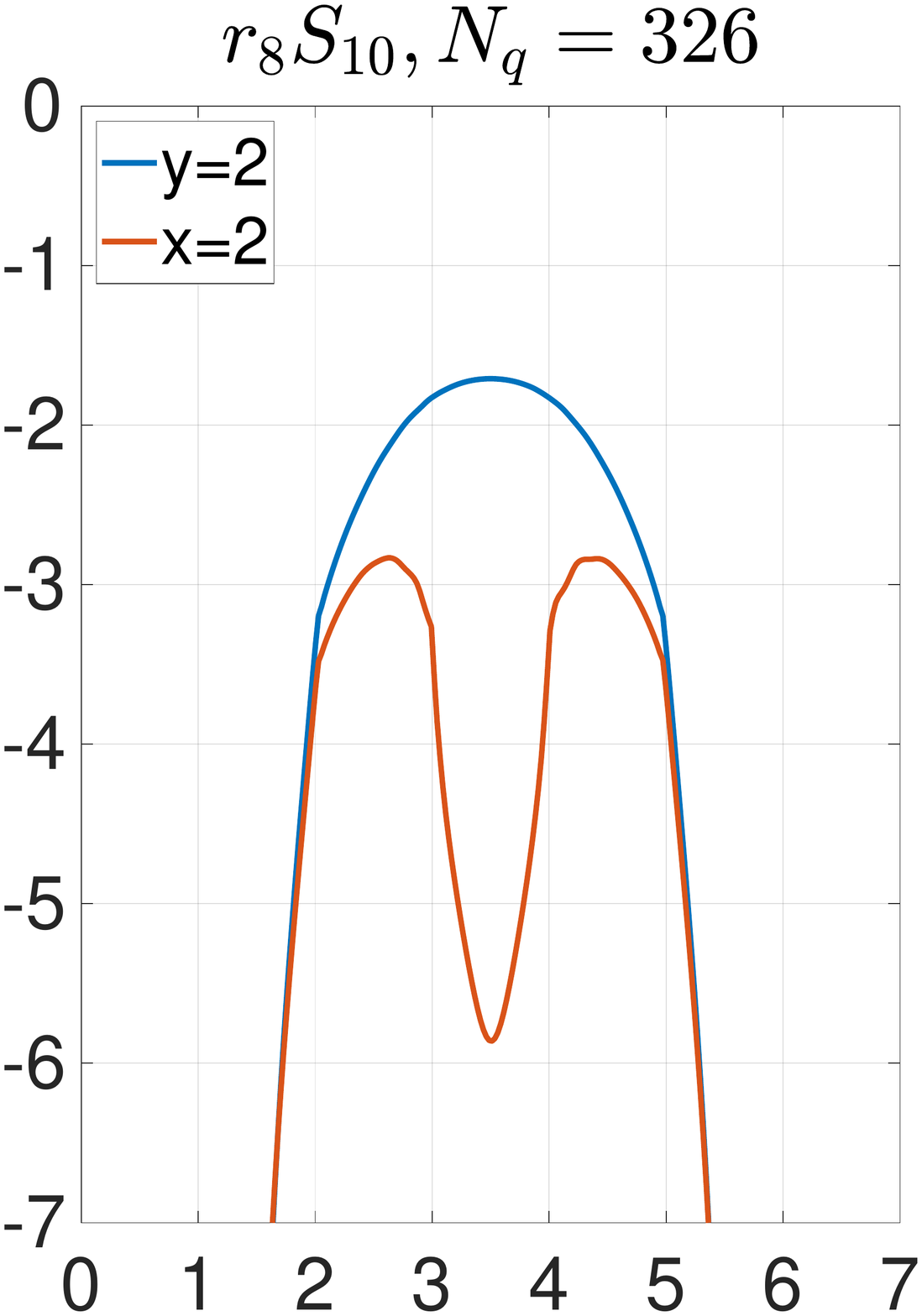}
		\label{fig:sub3}
	\end{subfigure}
	\caption{Cross sections for the lattice problem.}%
	\label{fig:checkerboardcross}
\end{figure}

To demonstrate that the \rSN method does not only yield satisfactory results 
for radially symmetric problems, we consider a lattice problem in the second 
application. The \emph{lattice problem} simulates a source within a heterogeneous material. 
The two-dimensional physical space consists of multiple materials, namely a set of 
squares belonging to a strongly absorbing medium as well as a strongly 
scattering background \cite{brunner2002forms,brunner2005two}. Furthermore, an 
isotropic source is placed into the center of the physical domain. There is no 
initial distribution of particles. Particles solely enter the domain by the 
source term 
isotropically. Again, we investigate the solution's dependency on the number of 
ordinates and the rotation strength.
The used parameters, as well as the underlying layout are shown in Tab. 
\ref{fig:cbSetting} and Fig.~\ref{fig:cblayout}, respectively.
The results are summarized in Fig.~\ref{fig:checkerboard} with the logarithmic 
density 
distribution along the horizontal (blue) and vertical (red) cut shown in 
Fig.~\ref{fig:checkerboardcross}. Similar to the line-source problem, we 
organize the plots such that similar number of quadrature points can be found 
in every column and the rotation strength is fixed within one row. The first 
row shows the result for the standard \SN method, all following rows show the 
solution of the \rSN method with increasing rotation strength.
Up to $N_q=324$ quadrature points we observe ray effects in the \SN method, 
shown in the last column of the first row in Fig.~\ref{fig:checkerboard}. The 
solution along $x=1$ and $y=1$ still shows oscillatory behavior, seen in  
Fig.~\ref{fig:checkerboardcross}. As before, the \rSN method without rotation 
(i.e. $\delta=0$) shows similar ray effects compared to the standard \SN 
method. While these ray effects are 
similar in structure and strength, the direction of the rays is slightly 
different due to the different quadrature set. Activating the rotation with 
$\delta=4$ or $\delta=8$ in the third and fourth row visibly diminishes the 
ray effects. For $\delta=4$, the ray effects seem to be diminished with 
$N_q=326$ and for $\delta=8$ with $N_q=102$, respectively. Remarkable is the 
absence of strong ray effects for all number of quadrature points with 
activated rotation. A convergent behavior of the solution along the horizontal 
and vertical cut can also be observed in  
Fig.~\ref{fig:checkerboardcross}. The last two rows indicate a quality that is 
not matched by the standard \SN method. For the checkerboard case, the rotation 
strength influences the solution in a less significant way then for the 
line-source problem.

%%%%%%%%%%%

%%%%%%%%%%%

%%%%%%%%%%%

%%%%%%%%%%%

\section{Summary and outlook}

To mitigate ray effects of the \SN solution, we have introduced a rotation of the quadrature set. By choosing a mesh-based quadrature rule, we enabled an efficient interpolation step to compute the solution values at the rotated quadrature nodes. In addition to mitigating ray effects, the \rSN method is easy to implement and promises positivity of solution values. It can be shown analytically that in a simplified setting, the rotation and interpolation step adds a diffusive term of the numerical discretization. Furthermore, we have provided a guideline on how to choose the rotation angle depending on time step and angular grid.

We tested our method on the line-source and lattice problems, and observed a mitigation of ray effects.

Future work will focus on a rigorous analysis of the modified equations when performing a general 3D rotation step. Furthermore, different mesh-based quadrature sets should be studied to identify a quadrature rule with desirable properties, such as homogeneity of integration weights. Additionally, different versions on the \rSN method need to be compared. These for example include back and forth rotation before streaming. In an implicit or steady state transport calculation, the ability to perform transport sweeps is of key importance. As the method has been presented, it is not directly compatible with sweeping because the rotation might actually change the face of the boundary from which the solution is propagated. However, one could solve a modified equation directly (without rotations). This will also be investigated in future work.

\section*{Acknowledgements}
The authors acknowledge many fruitful discussions with Cory D.\ Hauck (Oak 
Ridge) and Ryan G.\ McClarren (Notre Dame).

\bibliographystyle{siamplain}

\bibliography{main}

\begin{thebibliography}{10}

\bibitem{abu2001angular}
{\sc I.~Abu-Shumays}, {\em Angular quadratures for improved transport
  computations}, {Transport Theory and Statistical Physics}, 30 (2001),
  pp.~169--204.

\bibitem{brunner2002forms}
{\sc T.~A. Brunner}, {\em Forms of approximate radiation transport}, Sandia
  report,  (2002).

\bibitem{brunner2005two}
{\sc T.~A. Brunner and J.~P. Holloway}, {\em Two-dimensional time dependent
  {Riemann} solvers for neutron transport}, {Journal of Computational Physics},
  210 (2005), pp.~386--399.

\bibitem{JuliaWN}
{\sc T.~Camminady, M.~Frank, K.~Kuepper, and J.~Kusch}, {\em Source code {rSn}
  method}, 2018, \url{https://git.scc.kit.edu/qd4314/rSN}.

\bibitem{case1967linear}
{\sc K.~M. Case and P.~F. Zweifel}, {\em Linear transport theory},  (1967).

\bibitem{fleck1971implicit}
{\sc J.~Fleck~Jr and J.~Cummings~Jr}, {\em An implicit {Monte Carlo} scheme for
  calculating time and frequency dependent nonlinear radiation transport},
  {Journal of Computational Physics}, 8 (1971), pp.~313--342.

\bibitem{fryer2006snsph}
{\sc C.~L. Fryer, G.~Rockefeller, and M.~S. Warren}, {\em {SNSPH}: a parallel
  three-dimensional smoothed particle radiation hydrodynamics code}, The
  {Astrophysical Journal}, 643 (2006), p.~292.

\bibitem{garrett2013comparison}
{\sc C.~K. Garrett and C.~D. Hauck}, {\em A comparison of moment closures for
  linear kinetic transport equations: {The} line source benchmark}, {Transport
  Theory and Statistical Physics}, 42 (2013), pp.~203--235.

\bibitem{jung1972discrete}
{\sc J.~Jung, H.~Chijiwa, K.~Kobayashi, and H.~Nishihara}, {\em Discrete
  ordinate neutron transport equation equivalent to {PL} approximation},
  {Nuclear Science and Engineering}, 49 (1972), pp.~1--9.

\bibitem{lathrop1971remedies}
{\sc K.~Lathrop}, {\em Remedies for ray effects}, {Nuclear Science and
  Engineering}, 45 (1971), pp.~255--268.

\bibitem{lathrop1968ray}
{\sc K.~D. Lathrop}, {\em Ray effects in discrete ordinates equations},
  {Nuclear Science and Engineering}, 32 (1968), pp.~357--369.

\bibitem{levequenumerical}
{\sc R.~J. LeVeque}, {\em Numerical {Methods for Conservation Laws}. 1992},
  Birkha{\O} user Basel.

\bibitem{lewis1984computational}
{\sc E.~E. Lewis and W.~F. Miller}, {\em Computational methods of neutron
  transport},  (1984).

\bibitem{marinak2001three}
{\sc M.~Marinak, G.~Kerbel, N.~Gentile, O.~Jones, D.~Munro, S.~Pollaine,
  T.~Dittrich, and S.~Haan}, {\em Three-dimensional {HYDRA} simulations of
  {National Ignition Facility} targets}, {Physics of Plasmas}, 8 (2001),
  pp.~2275--2280.

\bibitem{mathews1999propagation}
{\sc K.~A. Mathews}, {\em On the propagation of rays in discrete ordinates},
  Nuclear science and engineering, 132 (1999), pp.~155--180.

\bibitem{matzen2005pulsed}
{\sc M.~K. Matzen, M.~Sweeney, R.~Adams, J.~Asay, J.~Bailey, G.~Bennett,
  D.~Bliss, D.~Bloomquist, T.~Brunner, R.~e. Campbell, et~al.}, {\em
  Pulsed-power-driven high energy density physics and inertial confinement
  fusion research}, {Physics of Plasmas}, 12 (2005), p.~055503.

\bibitem{mcclarren2010robust}
{\sc R.~G. McClarren and C.~D. Hauck}, {\em Robust and accurate filtered
  spherical harmonics expansions for radiative transfer}, {Journal of
  Computational Physics}, 229 (2010), pp.~5597--5614.

\bibitem{miller1977ray}
{\sc W.~Miller~Jr and W.~H. Reed}, {\em Ray-effect mitigation methods for
  two-dimensional neutron transport theory}, {Nuclear Science and Engineering},
  62 (1977), pp.~391--411.

\bibitem{morel2003analysis}
{\sc J.~Morel, T.~Wareing, R.~Lowrie, and D.~Parsons}, {\em Analysis of
  ray-effect mitigation techniques}, Nuclear science and engineering, 144
  (2003), pp.~1--22.

\bibitem{pomraning1973equations}
{\sc G.~C. Pomraning}, {\em The equations of radiation hydrodynamics}, Courier
  Corporation, (1973).

\bibitem{reed1972spherical}
{\sc W.~H. Reed}, {\em Spherical harmonic solutions of the neutron transport
  equation from discrete ordinate codes}, {Nuclear Science and Engineering}, 49
  (1972), pp.~10--19.

\bibitem{swesty2009numerical}
{\sc F.~D. Swesty and E.~S. Myra}, {\em A numerical algorithm for modeling
  multigroup neutrino-radiation hydrodynamics in two spatial dimensions}, {The
  Astrophysical Journal Supplement Series}, 181 (2009), p.~1.

\bibitem{tencer2016ray}
{\sc J.~Tencer}, {\em Ray effect mitigation through reference frame rotation},
  Journal of Heat Transfer, 138 (2016), p.~112701.

\bibitem{thurgood1995tn}
{\sc C.~Thurgood, A.~Pollard, and H.~Becker}, {\em The {TN} quadrature set for
  the discrete ordinates method}, {Journal of heat transfer}, 117 (1995),
  pp.~1068--1070.

\end{thebibliography}
\end{document}